\documentclass{amsart}

\usepackage[T1]{fontenc}
\usepackage[utf8]{inputenc}
\usepackage[american]{babel}

\numberwithin{equation}{section}

\theoremstyle{plain}
\newtheorem{thm}{Theorem}[section]
\newtheorem{lem}[thm]{Lemma}
\newtheorem{cor}[thm]{Corollary}

\theoremstyle{remark}
\newtheorem{rem}[thm]{Remark}

\usepackage{amsmath}
\usepackage{amssymb}
\usepackage{mathtools}
\usepackage{enumitem}

\usepackage[hyperfootnotes=false,pdfpagelabels,hidelinks]{hyperref}

\usepackage{lmodern}
\usepackage[scr=rsfs]{mathalfa}
\usepackage[only,rrbracket,llbracket]{stmaryrd}

\newcount\afterconjunctionpenalty
\makeatletter
\@namedef{u8:\detokenize{⋅}}{\penalty\afterconjunctionpenalty\ }
\makeatother
\newcommand*\hairspace{\kern 0.08333em}

\DeclareFontFamily{OT1}{mathc}{}
\DeclareFontShape{OT1}{mathc}{m}{it}{<-> mathc10}{}
\DeclareMathAlphabet{\mathabxcal}{OT1}{mathc}{m}{it}
\ExplSyntaxOn
\NewDocumentCommand\op{m}{%
  \tl_set:Nn\l_tmpa_tl {#1}
  \regex_replace_all:nnN {[a-z]}{\c{mathabxcal}\cB\{\0\cE\}}\l_tmpa_tl
  \regex_replace_all:nnN {[A-Z]}{\c{mathcal}\cB\{\0\cE\}}\l_tmpa_tl
  \l_tmpa_tl
}
\ExplSyntaxOff

\let\set\mathscr
\let\wt\tilde
\let\wb\bar
\let\wh\hat
\newcommand*{\spcheck}[1][A]{\check{\vphantom{#1}}}

\renewcommand{\Re}{\operatorname{Re}}

\providecommand{\rbr}[1]{(#1)}
\providecommand{\bigrbr}[1]{\bigl(#1\bigr)}
\providecommand{\biggrbr}[1]{\biggl(#1\biggr)}
\providecommand{\sbr}[1]{[#1]}
\providecommand{\cbr}[1]{\{#1\}}
\providecommand{\bigcbr}[1]{\bigl\{#1\bigr\}}
\providecommand{\abr}[1]{\langle#1\rangle}
\providecommand{\abs}[1]{\lvert#1\rvert}
\providecommand{\bigabs}[1]{\bigl\lvert#1\bigr\rvert}
\providecommand{\norm}[1]{\lVert#1\rVert}
\providecommand{\Norm}[1]{\mathopen{|\mkern-1.5mu|\mkern-1.5mu|}#1\mathclose{|\mkern-1.5mu|\mkern-1.5mu|}}
\let\seminorm\sbr

\newcommand\restrict[1]{\raisebox{-0.1ex}{$\rvert_{#1}$}}

\DeclareMathOperator{\diam}{diam}
\DeclareMathOperator{\dist}{dist}
\DeclareMathOperator{\Div}{div}
\DeclareMathOperator{\adj}{adj}
\DeclareMathOperator{\supp}{supp}

\newcommand*{\findiff}{\varDelta}
\newcommand*{\join}{{\vee}}
\newcommand*{\meet}{{\wedge}}

\def\B{\mathbf{B}}
\def\C{\mathbb{C}}
\def\N{\mathbb{N}}
\def\R{\mathbb{R}}

\DeclareMathOperator{\VMO}{VMO}

\newcommand*{\osc}{\omega}
\newcommand*{\wc}{{\mkern 2mu\cdot\mkern 2mu}}
\newcommand*{\dd}{\mathop{}\!d}

\makeatletter
\newcommand{\aev}{a.e\@ifnextchar.{}{.}}
\newcommand{\cf}{cf.}
\newcommand{\eg}{e.g.}
\newcommand{\etc}{etc\@ifnextchar.{}{.}}
\newcommand{\ie}{i.e.}
\makeatother

\newcommand*{\dual}{{+}}
\let\hconj\prime
\newcommand*{\dualref}[1]{\hyperref[#1]{\ref*{#1}\textsuperscript{+}}}

\newcommand*{\ball}{B}
\newcommand*{\cube}{Q}
\newcommand*{\cell}{\set{\cube}}
\newcommand*{\cellradius}{r}
\newcommand*{\lattice}{\set{L}}
\newcommand*{\nullset}{\set{N}}

\usepackage{bbm}
\newcommand*{\chf}{\mathbbm{1}}

\usepackage{subcaption}

\usepackage{tikz}
\usetikzlibrary{calc}
\usetikzlibrary{intersections}

\tikzset{
  boundary/.style = {semithick},
  interface/.style = {boundary},
  lattice/.style = {very thin, black!50},
  inclusion/.style = {fill=black!#1, draw=none},
  matrix/.style = {fill=gray!#1},
  leader/.style = {thin},
  crosspoint/.style = {circle, fill, inner sep=0pt, minimum size=3pt},
  supp/.style = {thin, dashed},
  dimension/.style = {thin, latex-latex},
  latticedimmed/.style = {thin, black!25},
}

\begin{document}
\frenchspacing

\title[Homogenization for Piecewise Locally Periodic Problems]
      {Homogenization for \\ Piecewise Locally Periodic Problems on a Domain}
\author{Nikita N. Senik}
\address{Saint~Petersburg State University, Universitetskaya nab.~7/9, Saint~Petersburg 199034, Russia}
\email{nnsenik@gmail.com}
\thanks{Research supported by RSF grant~22-11-00092-P}


\date{\today.}


\begin{abstract}
We study homogenization of a piecewise locally periodic elliptic operator
of the form $-\operatorname{div}A(x,x/\varepsilon)\nabla$ on a Lipschitz domain
divided by interfaces into finitely many subdomains. In the fast variable,
the coefficient $A$ is merely bounded and periodic, with a periodic structure
and a microscopic scale specific to each subdomain. In the slow variable,
it is piecewise uniformly continuous, with a modulus of continuity $\omega$
(or, in general, one for each subdomain), and may jump across the interfaces.
The boundary conditions include Dirichlet, Robin-type and mixed ones.
Assuming that the resolvent is uniformly bounded between
the corresponding Sobolev spaces and that the effective problem
is piecewise regular in the scale of Lipschitz--Besov-type
spaces $\Lambda_p^{1+\omega}$, we obtain operator error estimates in $L_p$
and, with a corrector, in $W_p^1$; for a single microscopic scale $\varepsilon$,
both rates are $\smash[t]{\omega(\varepsilon)^{1/p}}$. If the adjoint
effective problem is regular as well, the rate in $L_p$ improves
to $\omega(\varepsilon)$, as does the rate in $W_p^1$ away from
the boundary and the interfaces. For instance, we show that,
with $p$ close to $2$, these hypotheses hold for strongly elliptic systems
on domains whose boundaries and interfaces are of class $C^{1,\omega}$
and do not meet.
To this end, we prove regularity in the scale of $\Lambda_p^{1+\omega}$-spaces
for elliptic problems with piecewise regular coefficients and data.
\end{abstract}

\maketitle

\section{Introduction}

Differential operators with rapidly oscillating periodic coefficients
arise naturally in models of physical processes in heterogeneous
media.  As⋅the⋅period of the oscillations becomes small, such
processes can often be described by effective models with nonoscillating
coefficients. In⋅mathematical terms, this amounts to approximating
the resolvent of the original differential operator by that of the
effective operator. Quantitative homogenization is concerned with
the accuracy of this approximation.    One direction within this
field, which has grown into an area of its own, focuses on estimates
 in  operator norms, known as operator error estimates.

Operator error estimates are by now well established for purely periodic
operators (see the survey~\cite{ZhPas:2016}, the monograph~\cite{Shen:2018}
and the references therein). In⋅contrast, locally periodic operators
have received less attention. Their coefficients may depend not
only on a rapidly oscillating periodic variable (the⋅fast variable)
but also on a macroscopic variable (the⋅slow variable) describing
the large-scale variation of the medium. This allows, \eg, the
local amplitude of the oscillations to vary gradually in space. In⋅such
a setting, the rate of approximation is affected by the regularity
of the coefficients in the slow variable: if they are H\"older continuous,
the known rates are powers of⋅$\varepsilon^{s}$, with⋅$s$ the H\"older
exponent (see~\cite{Senik:2023-1}). When this condition is relaxed
to mere uniform continuity, one naturally expects⋅$\varepsilon^{s}$
to be replaced by  the modulus of continuity at~$\varepsilon$.
One can go further and require uniform continuity only within individual
phases of the medium, so that the coefficients may jump across the
interfaces (as⋅in bicrystals, layered structures or welded composite
materials). The⋅present paper is devoted  precisely to the case
of piecewise uniformly continuous coefficients.

\begin{figure}
\centering
\begin{subfigure}{0.48\textwidth}
\centering
\begin{tikzpicture}[
  declare function={
    RI=1.6; RII=0.6; 
    DI=0.3; DII=0.2; 
    eI=0.42; eII=0.30; 
    aI=-20; aII=12; 
    dI=0.25*eI; 
    axII=0.32*eII; ayII=0.26*eII; 
  },
]
  \pgfmathsetseed{1024}
  \path[save path=\OuterDomain]
    plot[smooth cycle, variable=\t, samples at={0,30,...,330}] (\t:{RI+DI*rnd});
  \pgfmathsetseed{2048}
  \path[save path=\InnerDomain, name path global=Gamma]
    plot[smooth cycle, variable=\t, samples at={0,30,...,330}] (\t:{RII+DII*rnd});

  \begin{scope}
    \clip[use path=\OuterDomain];
    \shade[left color=gray!10, right color=gray!30, shading angle=90+aI]
      ({-(RI+DI)},{-(RI+DI)}) rectangle ({(RI+DI)},{(RI+DI)});
    \begin{scope}[rotate=aI]
      \pgfmathtruncatemacro{\N}{ceil((RI+DI)/eI)}
      \draw[lattice] ({-\N*eI},{-\N*eI}) grid[step=eI] ({\N*eI},{\N*eI});
      \clip \foreach \i in {-\N,...,\N} {
        \foreach \j in {-\N,...,\N} {
          ({(\i+1/2)*eI},{(\j+1/2)*eI}) circle[radius=dI]
        }
      };
      \shade[left color=black!55, right color=black!85]
        (-{(RI+DI)},-{(RI+DI)}) rectangle ({(RI+DI)},{(RI+DI)});
    \end{scope}
  \end{scope}

  \begin{scope}
    \clip[use path=\InnerDomain];
    \shade[inner color=gray!0, outer color=gray!30]
      (-{(RII+DII)},-{(RII+DII)}) rectangle ({(RII+DII)},{(RII+DII)});
    \begin{scope}[rotate=aII]
      \pgfmathtruncatemacro{\M}{ceil((RII+DII)/eII)}
      \draw[lattice] ({-\M*eII},{-\M*eII}) grid[step=eII] ({\M*eII},{\M*eII});
      \clip \foreach \i in {-\M,...,\M} {
        \foreach \j in {-\M,...,\M} {
          ({(\i+1/2)*eII},{(\j+1/2)*eII}) ellipse[x radius=axII, y radius=ayII]
        }
      };
      \shade[inner color=black!85, outer color=black!55]
        (-{(RII+DII)},-{(RII+DII)}) rectangle ({(RII+DII)},{(RII+DII)});
    \end{scope}
  \end{scope}

  \draw[boundary, use path=\OuterDomain];
  \draw[interface, use path=\InnerDomain];

  \draw[leader] (35:{0.8*RI}) -- (35:{1.3*RI})
    node[anchor=south west, inner sep=1pt] {$\Omega_-$};
  \draw[leader] (215:{0.35*RII}) -- (215:{1.3*RI})
    node[anchor=north east, inner sep=1pt] {$\Omega_+$};
  \path[name path=ray] (0,0) -- (125:{2*RI});
  \path[name intersections={of=Gamma and ray, by=G}];
  \draw[leader] (G) -- (125:{1.3*RI})
    node[anchor=south east, inner sep=1pt] {$\Gamma$};
\end{tikzpicture}\caption{}
\label{subfig: OO | Smooth domains}\end{subfigure}%
\hfill%
\begin{subfigure}{0.48\textwidth}
\centering
\begin{tikzpicture}[
  declare function={
    s=2.8; 
    eI=s/6; eII=s/6; 
    aI=-10; aII=-aI; 
    dI=0.25*eI; 
    dII=0.25*eII; 
  },
]
  \begin{scope}[shift={(-s/2,s/2)}]
    \clip (-s/2,-s/2) rectangle ++(s,s);
    \begin{scope}[rotate=aI]
    \pgfmathsetmacro{\maxdistI}{s*(cos(aI)+abs(sin(aI)))/2}
    \shade[left color=gray!0, right color=gray!50, shading angle=90+aI]
      (-\maxdistI,-\maxdistI) rectangle (\maxdistI,\maxdistI);
      \pgfmathtruncatemacro{\N}{ceil(\maxdistI/eI)}
      \draw[lattice] ({-\N*eI},{-\N*eI}) grid[step=eI] ({\N*eI},{\N*eI});
      \clip \foreach \i in {-\N,...,\N} {
        \foreach \j in {-\N,...,\N} {
          ({(\i+1/2)*eI},{(\j+1/2)*eI}) circle[radius=dI]
        }
      };
      \shade[left color=black!50, right color=black!100]
        (-\maxdistI,-\maxdistI) rectangle (\maxdistI,\maxdistI);
    \end{scope}
  \end{scope}

  \begin{scope}[shift={(s/2,s/2)}]
    \clip (-s/2,-s/2) rectangle ++(s,s);
    \begin{scope}[rotate=aII]
    \pgfmathsetmacro{\maxdistII}{s*(cos(aII)+abs(sin(aII)))/2}
    \shade[left color=gray!50, right color=gray!0, shading angle=90+aII]
      (-\maxdistII,-\maxdistII) rectangle (\maxdistII,\maxdistII);
      \pgfmathtruncatemacro{\N}{ceil(\maxdistII/eII)}
      \draw[lattice] ({-\N*eII},{-\N*eII}) grid[step=eII] ({\N*eII},{\N*eII});
      \clip \foreach \i in {-\N,...,\N} {
        \foreach \j in {-\N,...,\N} {
          ({(\i+1/2)*eII},{(\j+1/2)*eII}) circle[radius=dII]
        }
      };
      \shade[left color=black!100, right color=black!50]
        (-\maxdistII,-\maxdistII) rectangle (\maxdistII,\maxdistII);
    \end{scope}
  \end{scope}

  \draw[boundary] (-s,0) rectangle (s,s);
  \draw[interface] (0,0) -- (0,s);
  \node[crosspoint] at (0,0) {};
  \node[crosspoint] at (0,s) {};

  \node[below, inner sep=2pt] at (-s/2,0) {$\Omega_{\mathrlap{-}}$};
  \node[below, inner sep=2pt] at (s/2,0)  {$\Omega_{\mathrlap{+}}$};
  \draw[leader] (0,0.75*s) -- (0.2*s,1.1*s)
    node[anchor=south west, inner sep=1pt] {$\Gamma$};
\end{tikzpicture}\caption{}
\label{subfig: OO | Lipschitz domains w/ boundary cross-points}

\end{subfigure}
\caption{A domain~$\Omega$ divided by an interface~$\Gamma$ into subdomains~$\Omega_{\pm}$,
each with its own periodic structure: (\subref{subfig: OO | Smooth domains})~a⋅smooth
domain with a compactly contained smooth subdomain; (\subref{subfig: OO | Lipschitz domains w/ boundary cross-points})~a⋅Lipschitz
domain in which the interface meets the boundary at two cross-points.
The shading indicates the dependence of the coefficient on the slow
variable.}
\label{fig: OO | Domains}
\end{figure}
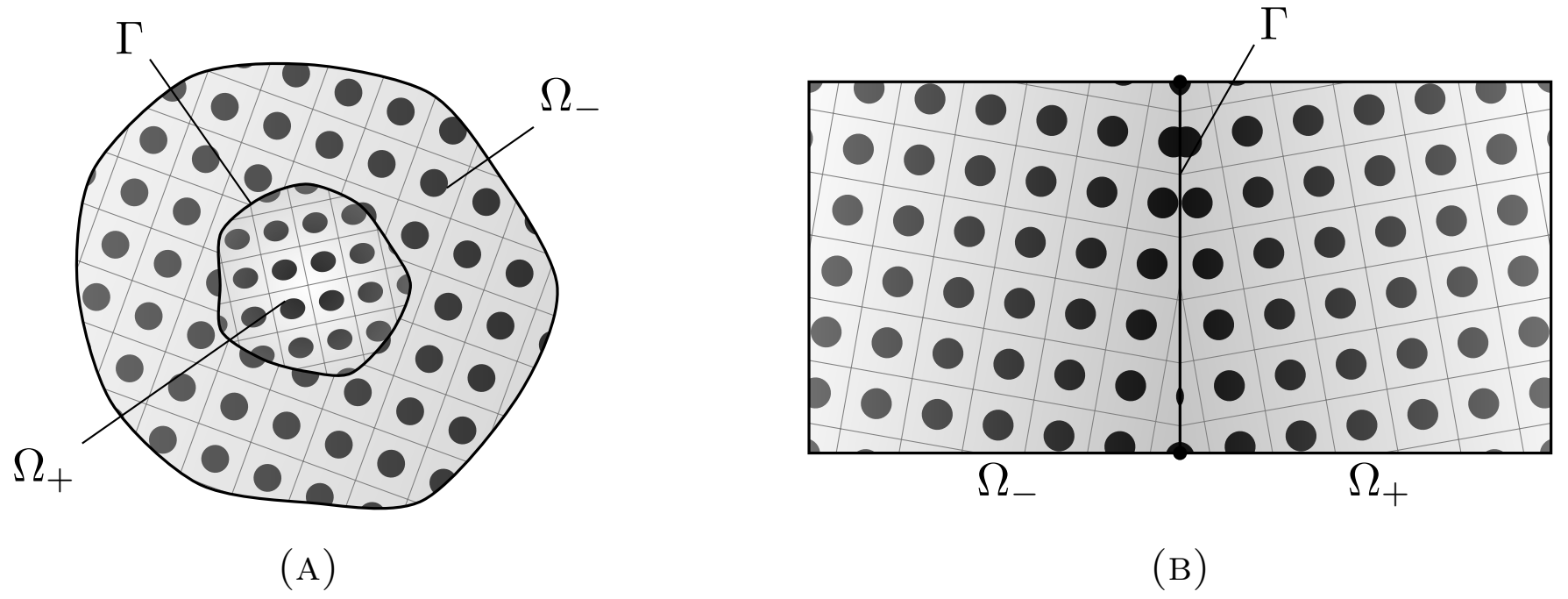

To⋅state our results in a simple setting,  let⋅$\Omega$ be a
Lipschitz domain in⋅$\R^{d}$ divided by a⋅Lipschitz surface~$\Gamma\subset\Omega$
into two  subdomains~$\Omega_{+}$ and~$\Omega_{-}$ (see⋅Figure~\ref{fig: OO | Domains}).
Let⋅$A\colon\Omega\times\nobreak\R^{d}\to\B\rbr{\C^{n\times d}}$
be a⋅bounded function that is uniformly continuous in its first variable
on each⋅$\Omega_{\pm}$ and periodic in its second, possibly with
respect to different lattices on~$\Omega_{\pm}$. We assume that
the continuity on⋅$\Omega_{\pm}$ is controlled by a  regular modulus
of continuity~$\omega_{\pm}$. The microscopic scale on each⋅$\Omega_{\pm}$
is described by⋅$\varepsilon_{\pm}=\varepsilon_{\pm}\rbr{\varepsilon}$,
with⋅$\varepsilon_{\pm}\rbr{\varepsilon}\to0$ as~$\varepsilon\to0$.
We set $A^{\varepsilon}\rbr x=A\rbr{x,x/\varepsilon_{\#}}$, where
$\varepsilon_{\#}=\varepsilon_{\pm}$ on $\Omega_{\pm}$, and consider
an operator on⋅$\Omega$ of the form
\[
\op A^{\varepsilon}=-\Div A^{\varepsilon}\nabla
\]
subject to fairly general boundary conditions. These include the
Dirichlet condition~$u=0$  and a Robin-type condition~$\partial_{\nu_{\op A^{\varepsilon}}}u+\sigma u=0$
 for a suitable operator~$\sigma$ on the boundary, as well as their
combinations on different parts of⋅$\partial\Omega$ or for different
components of⋅$u$ when⋅$u$ is vector-valued.  By⋅$\smash{\set W_{q}^{1}}\rbr{\Omega;\C^{n}}$,
we denote the subspace of the complex Sobolev space~$\smash{W_{q}^{1}}\rbr{\Omega}^{n}$
corresponding to the  boundary conditions.

Throughout, we assume that, for some⋅$\mu\in\C$ and⋅$p\in\rbr{1,\infty}$,
the operator~$\op A^{\varepsilon}-\mu$ is an isomorphism from⋅$\set W_{p}^{1}\rbr{\Omega;\C^{n}}$
onto⋅$\set W_{p^{\prime}}^{1}\rbr{\Omega;\C^{n}}^{*}$, with its inverse
uniformly bounded in~$\varepsilon$. A⋅standard sufficient condition
for this hypothesis, and the one we have in mind, is strong ellipticity,
at least for⋅$p$ close to~$2$. It⋅is not, however, needed for
what follows. On⋅the⋅other hand, the uniform boundedness of the
resolvents is essential even for  convergence in the weak operator
topology. In⋅the⋅classical strongly elliptic setting, the resolvent
of⋅$\op A^{\varepsilon}$ converges to that of an effective operator~$\op A^{0}$
(see~\cite{BLP:1978}, \cite{BP:1984},~\cite{ZhKO:1993}). The⋅latter
has the same form as⋅$\op A^{\varepsilon}$ but with nonoscillating
coefficients, which are obtained from the solution of an auxiliary
problem on the periodicity cell, the cell problem. In⋅the absence
of strong ellipticity, we require that the cell problem be uniquely
solvable for each value of the slow variable, with the solution operator
bounded uniformly in that variable. This implies a weaker ellipticity-type
condition on~$A$, see Remark~\ref{Rem: EO | Weaker ellipticity-type condition}.
We further suppose that the inverse of⋅$\op A^{0}-\mu$ is bounded
as an operator from⋅$L_{p}\rbr{\Omega}^{n}$ to⋅$\smash{\set W_{p}^{1}\rbr{\Omega;\C^{n}}}\cap\smash{\Lambda_{p}^{1+\omega_{\pm}}\rbr{\Omega_{\pm}}^{n}}$,
where $\smash{\Lambda_{p}^{1+\omega_{\pm}}}\rbr{\Omega_{\pm}}$ is
the space of functions in⋅$W_{p}^{1}\rbr{\Omega_{\pm}}$ whose gradients
have $L_{p}$\nobreakdash-moduli of continuity  controlled by~$\omega_{\pm}$.
This is an assumption of elliptic regularity type, reflecting 
properties of the coefficients, the geometry and the boundary conditions.

Under⋅the⋅above hypotheses, we prove that
\begin{align}
\norm{\rbr{\op A^{\varepsilon}-\mu}^{-1}-\rbr{\op A^{0}-\mu}^{-1}}_{L_{p}\rbr{\Omega}^{n}\to L_{p}\rbr{\Omega}^{n}} & \le C\varpi_{\join}\rbr{\varepsilon}^{1/p},\label{est: Intro | Convergence in L_p}\\
\norm{\nabla\rbr{\op A^{\varepsilon}-\mu}^{-1}-\nabla\rbr{\op A^{0}-\mu}^{-1}-\nabla\varepsilon_{\#}\op K_{\mu}^{\varepsilon}}_{L_{p}\rbr{\Omega}^{n}\to L_{p}\rbr{\Omega}^{n\times d}} & \le C\varpi_{\join}\rbr{\varepsilon}^{1/p},\label{est: Intro | Approximation in W_p^1}
\end{align}
where⋅$\varpi_{\join}\rbr{\varepsilon}$ is the larger of⋅$\varpi_{\pm}\rbr{\varepsilon}=\omega_{\pm}\rbr{\varepsilon_{\pm}}$
(see⋅Theorem~\ref{thm: MR | Conv in L_q and approx in W_p^1}). As⋅usual,
the⋅$L_{p}$\nobreakdash-limit captures only  the macroscopic behavior
of the solution, while approximation in the⋅$W_{p}^{1}$\nobreakdash-norm
requires an additional corrector~$\op K_{\mu}^{\varepsilon}$ to
recover  the microscopic   oscillatory structure. If⋅the⋅regularity
assumption on⋅$\op A^{0}$ also holds for its adjoint, then the first
estimate improves⋅to
\begin{equation}
\norm{\rbr{\op A^{\varepsilon}-\mu}^{-1}-\rbr{\op A^{0}-\mu}^{-1}}_{L_{p}\rbr{\Omega}^{n}\to L_{p}\rbr{\Omega}^{n}}\le C\varpi_{\join}\rbr{\varepsilon}\label{est: Intro | Better convergence in L_p}
\end{equation}
(see⋅Theorem~\ref{thm: MR | Conv in L_p with better rate}, which,
 more generally, allows different moduli in the two regularity assumptions).
The⋅second estimate is also improved in this case, but only away
from the boundary and the interface. Namely, for every domain~$\Sigma$
compactly contained in either of⋅$\Omega_{\pm}$, we have
\begin{equation}
\norm{\nabla\rbr{\op A^{\varepsilon}-\mu}^{-1}-\nabla\rbr{\op A^{0}-\mu}^{-1}-\nabla\varepsilon_{\#}\op K_{\mu}^{\varepsilon}}_{L_{p}\rbr{\Omega}^{n}\to L_{p}\rbr{\Sigma}^{n\times d}}\le C\varpi_{\join}\rbr{\varepsilon}\label{est: Intro | Interior approximation in W_p^1}
\end{equation}
provided that a uniform $L_{p}$\nobreakdash-Caccioppoli-type inequality
 holds for the operator~$\op A^{\varepsilon}$ (see⋅Theorem~\ref{thm: MR | Interior approx in W_p^1}).

As⋅our main example,  we consider strongly elliptic operators in
the case where the subdomains are of class~$C^{1,\omega_{\pm}}$
and the coefficients are⋅$C^{0,\omega_{\pm}}$ in the slow variable,
with a mild additional condition near the interface (see⋅Section~\ref{sec: Examples}).
 All the above hypotheses are then satisfied for a range of exponents~$p$
near~$2$ (with⋅$p\ge2$ for the interior $W_{p}^{1}$\nobreakdash-estimate).
In⋅particular, we prove the required  regularity property of the
effective  problem, for which we could not find a reference (see⋅Appendix~\ref{sec: Appendix - Elliptic regularity}).
Some of these results, with Dirichlet or Neumann boundary conditions
and with⋅$p=2$ and⋅$\omega_{\pm}\rbr t=t$, were announced in~\cite{Senik:2023-2}.

Operator error estimates in homogenization were first obtained
by⋅Birman and⋅Suslina~\cite{BSu:2001} for purely periodic problems
in the whole space. They used spectral methods, combining the Floquet–Bloch
decomposition with analytic perturbation theory. Different, nonspectral
approaches, based on a shift technique and on periodic unfolding,
were developed in the works of⋅Zhikov~\cite{Zh:2005} and of⋅Zhikov
and⋅Pastukhova~\cite{ZhPas:2005} and in those of⋅Griso~\cite{Gri:2004},~\cite{Gri:2006},
respectively, extending  the  results of⋅\cite{BSu:2001} to problems
on domains. Since then, many further contributions have been made,
both for purely periodic operators (see,⋅\eg,⋅\cite{Nazarov:2006},
\cite{KLS:2012}, \cite{Su:2013-1} and~\cite{Su:2013-2}) and for
locally periodic ones (see⋅\cite{PasT:2007}, \cite{Bor:2008}, \cite{Se:2017-3},
\cite{Senik:2022} and~\cite{Senik:2023-1}).

A⋅detailed comparison of known operator error estimates for problems
without interfaces,  from purely periodic coefficients to those
that are globally H\"older continuous in the slow variable, can be
found in~\cite{Senik:2023-1}.  While recovering the results of
that work, the present paper goes beyond it in the two directions
outlined above. First, the regularity of the coefficients in the slow
variable and that of the effective  resolvent are measured by general
moduli of continuity rather than by H\"older  moduli. Second, the
coefficients and the gradient of the effective resolvent may have
jumps across an interface, and the periodic structure may differ
on the two sides. Our⋅results thus cover the whole range from coefficients
independent of the slow variable, as in the purely periodic case,
through Lipschitz and H\"older dependence, to coefficients that are
merely piecewise uniformly continuous in~it. This appears to be
the weakest regularity in the slow variable for which operator error
estimates have been obtained.

As⋅for⋅homogenization problems with interfaces, existing works focus
mainly on the purely periodic ones, establishing uniform Lipschitz
estimates and pointwise approximations of solutions. In⋅\cite{Josien:2019},
such  approximations of a suitably normalized solution and its gradient
were obtained for a purely periodic problem on⋅$\R^{d}$ with a flat
interface. The⋅setting, however, differs from ours, since the spectral
parameter there,⋅$\mu=0$, lies in the spectrum. Zhuge⋅\cite{Zhuge:2021}
considered a transmission problem on a bounded Lipschitz domain divided
by a Lipschitz interface. Both the Dirichlet  data and the jump of
the conormal derivative across the interface were allowed to be of
minimal regularity. When these data vanish, as in the present paper,
the results of⋅\cite{Zhuge:2021} imply convergence of the corresponding
resolvent in the operator norm on~$L_{2}$.  Since Zhuge's arguments
 were designed for rough data, the rate is only a power of⋅$\varepsilon$
with a small exponent, dictated by Meyers' estimate on the higher
integrability of gradients rather than by  higher regularity of
the effective solution, which is unavailable in general. This is also
why no approximation for the gradient was given.

 Our approach builds on that of⋅\cite{Senik:2023-1}, which we briefly
describe first. That approach is based on an identity of resolvent
type for the error of a first-order approximation. Besides the effective
resolvent, this approximation involves a corrector suitably regularized
by mollification and Steklov smoothing (such regularizations go back
to⋅\cite{CDG:2002} and~\cite{ZhPas:2005}). The⋅identity   has
the form
\begin{equation}
\rbr{\op A^{\varepsilon}-\mu}^{-1}-\rbr{\op A^{0}-\mu}^{-1}-\varepsilon_{\#}\op K_{\mu}^{\varepsilon}\restrict{L_{p}\rbr{\Omega}^{n}}=\rbr{\op A^{\varepsilon}-\mu}^{-1}\op F_{\mu}^{\varepsilon}\rbr{\op A^{0}-\mu}^{-1},\label{eq: Intro | Resolvent Identity w/o interfaces}
\end{equation}
where⋅$\op F_{\mu}^{\varepsilon}$ collects the defect terms. Among
them, the leading ones cancel, which is the very heart of homogenization.
The⋅remaining ones split into interior and boundary parts, each
shown to be small: the interior part by following the arguments
for the whole-space case (see⋅\cite{Senik:2022}) and the boundary
part by using that it is supported in an⋅$\varepsilon$\nobreakdash-neighborhood
of the boundary, though with a generally larger error.

In⋅the⋅present paper, we extend this machinery from the H\"older
case to any regular modulus of continuity.  To⋅handle interfaces,
we reduce the problem to the subdomains, treating the interface as
part of their boundaries, a classical approach (see,⋅\eg,~\cite{BLP:1978}).
In⋅this way, the global error is represented as a sum of terms associated
with the subdomains, each of the same structure as the right-hand
side of⋅(\ref{eq: Intro | Resolvent Identity w/o interfaces}), with
the defect splitting into interior and boundary parts as before (see~(\ref{eq: Ps | Resolvent identity})).
Combining  the bounds for these parts yields (\ref{est: Intro | Convergence in L_p})
and~(\ref{est: Intro | Approximation in W_p^1}). Under⋅the⋅additional
regularity assumption on the adjoint of the effective operator, the
boundary part admits a better bound in the⋅$L_{p}$\nobreakdash-norm,
leading to~(\ref{est: Intro | Better convergence in L_p}). For⋅the
interior estimate~(\ref{est: Intro | Interior approximation in W_p^1}),
we use a localized version of the resolvent identity (see~(\ref{eq: Localized resolvent identity})),
which has no boundary part and therefore gives the improved rate.

The⋅paper is organized as follows. Section~\ref{sec: Notation} sets
up the basic notation. In⋅Section~\ref{sec: OO}, we describe the
problem geometry and the periodic structures and introduce the operator~$\op A^{\varepsilon}$
along with the first of our assumptions. Section~\ref{sec: EO} is
devoted to the cell problem and the effective operator and states
the remaining assumptions. In⋅Section~\ref{sec: C}, we define a
two-scale corrector on each subdomain, regularize it and assemble
the global corrector  from these pieces. The⋅main results are stated
in Section~\ref{sec: MR}. In⋅Section~\ref{sec: Examples}, we present
our main example, strongly elliptic operators, and verify  that they
satisfy all the hypotheses; the required regularity of the corresponding
effective  problem is proved in Appendix~\ref{sec: Appendix - Elliptic regularity}.
Finally, Section~\ref{sec: Proof of the main results} derives the
global and localized resolvent identities and uses them to prove the
main results.

\section{\label{sec: Notation}Notation}

We use the notation~$\norm{\wc}_{U}$ for the norm on a normed
space~$U$. If⋅$U$ and⋅$V$ are Banach spaces, then $\B\rbr{U,V}$
is the Banach space of bounded linear operators from⋅$U$ to~$V$.
When⋅$U=V$, the space~$\B\rbr U=\B\rbr{U,U}$ becomes a Banach algebra
with the identity map~$\op I$. The⋅norm and the inner product on⋅$\C^{n}$
are denoted by⋅$\abs{\wc}$ and⋅$\abr{\wc,\wc}$. We will  identify⋅$\B\rbr{\C^{n},\C^{m}}$
and~$\C^{m\times n}$.

Let⋅$\Sigma\subset\R^{d}$ be a domain and⋅$U$ a Banach space. Then
$L_{0}\rbr{\Sigma;U}$ is the vector space of all strongly measurable
functions from⋅$\Sigma$ to~$U$. The⋅symbol~$L_{p}\rbr{\Sigma;U}$
stands for the⋅$L_{p}$\nobreakdash-space of  functions in~$L_{0}\rbr{\Sigma;U}$.
 For⋅$s\ge0$, we let⋅$W_{p}^{s}\rbr{\Sigma;U}$ denote the usual
Sobolev space or Sobolev–Slobodetskii space of functions in~$L_{0}\rbr{\Sigma;U}$.
 When⋅$U=\C^{n}$, we write⋅$\norm{\wc}_{p,\Sigma}$ and⋅$\norm{\wc}_{s,p,\Sigma}$
for the norms on⋅$L_{p}\rbr{\Sigma}^{n}=L_{p}\rbr{\Sigma;U}$ and⋅$W_{p}^{s}\rbr{\Sigma}^{n}=W_{p}^{s}\rbr{\Sigma;U}$,
respectively, and⋅$\rbr{\wc,\wc}_{\Sigma}$ for the inner product
on⋅$L_{2}\rbr{\Sigma}^{n}$ and for its extensions. We often omit
the superscripts indicating the dimensions in the notation for spaces
of tensor-valued functions when they are clear from the context. The⋅conjugate
dual space of⋅$W_{p}^{s}\rbr{\Sigma}$ under the pairing~$\rbr{\wc,\wc}_{\Sigma}$
is denoted by⋅$\smash{W_{p}^{s}\rbr{\Sigma}^{*}}$ and its norm by~$\smash{\norm{\wc}_{-s,p^{\hconj},\Sigma}^{*}}$.
Here $p^{\hconj}$ is the exponent conjugate to⋅$p$, that is,⋅$1/p^{\hconj}=1-1/p$.
The⋅closure of⋅$\smash{C_{c}^{\infty}}\rbr{\Sigma}$ in⋅$\smash{W_{p}^{s}}\rbr{\Sigma}$
is⋅$\smash{\mathring{W}_{p}^{s}}\rbr{\Sigma}$, and $\smash{W_{p^{\hconj}}^{-s}}\rbr{\Sigma}$
is its dual with norm~$\norm{\wc}_{-s,p^{\hconj},\Sigma}$. As⋅usual,
$\smash[t]{H^{s}}=\smash[t]{W_{2}^{s}}$, $\smash[t]{H^{-s}}=\smash[t]{W_{2}^{-s}}$,~\etc.

Let⋅$\cell\subset\smash[t]{\R^{d}}$ be an open parallelepiped and~$m\in\N_{0}$.
We denote by⋅$\smash[t]{\wt C^{m}}\rbr{\cell}$ the class of $C^{m}$\nobreakdash-functions
on⋅$\cell$ whose $\cell$\nobreakdash-periodic extensions to⋅$\R^{d}$
are again $C^{m}$\nobreakdash-functions, and by⋅$\smash{\wt W_{p}^{m}}\rbr{\cell}$
the completion of⋅$\smash{\wt C^{m}}\rbr{\cell}$ in the⋅$\smash{W_{p}^{m}}$\nobreakdash-norm.
Notice that $\smash{\wt L_{p}}\rbr{\cell}=\smash{\wt W_{p}^{0}}\rbr{\cell}$
can be regarded as the space of all $\cell$\nobreakdash-periodic
functions in~$L_{p,\text{loc}}\rbr{\R^{d}}$.  The⋅spaces~$\smash{\wt W_{p}^{m}}\rbr{\Sigma\times\cell}$
are defined similarly. The⋅conjugate dual of⋅$\smash{\wt W_{p}^{m}}$
is denoted by~$\smash{\wt W_{p^{\hconj}}^{-m}}$.

In⋅addition to the Sobolev spaces, we  introduce spaces with a prescribed
$L_{p}$-modulus of continuity. A⋅function~$\omega\colon\smash[t]{\wb{\R}_{+}}\to\smash[t]{\wb{\R}_{+}}$
is called a regular modulus of continuity if $\omega\rbr 0=\lim_{t\to0}\omega\rbr t=0$,
$\omega\rbr 1=1$, the function~$t\mapsto\omega\rbr t$ is nondecreasing
and the function~$t\mapsto\omega\rbr t/t$ is nonincreasing (see~\cite{Fef:2009}).
For⋅such an⋅$\omega$ and any⋅$p\in\sbr{1,\infty}$ and⋅$m\in\N_{0}$,
 let⋅$\Lambda_{p}^{m+\omega}\rbr{\Sigma;U}$ denote the space of
functions~$u\in W_{p}^{m}\rbr{\Sigma;U}$ with finite norm
\begin{align*}
\abs u_{\Lambda_{p}^{m+\omega}\rbr{\Sigma;U}} & =\biggrbr{\sum_{\abs{\alpha}=m}\seminorm{D^{\alpha}u}_{\Lambda_{p}^{\omega}\rbr{\Sigma;U}}^{p}+\norm u_{W_{p}^{m}\rbr{\Sigma;U}}^{p}}^{1/p}
\end{align*}
(with⋅the⋅usual convention for~$p=\infty$), where
\[
\seminorm u_{\Lambda_{p}^{\omega}\rbr{\Sigma;U}}=\sup_{t\in\R_{+}}\omega\rbr t^{-1}\osc_{p}\rbr{u;t}
\]
and⋅$\osc_{p}\rbr{u;\wc}$ is the⋅$L_{p}$\nobreakdash-modulus of
continuity of⋅$u$, given⋅by
\[
\osc_{p}\rbr{u;t}=\sup_{\abs h\le t}\norm{\Delta_{h}u}_{L_{p}\rbr{\Delta_{h}\Sigma;U}}
\]
with⋅$\Delta_{h}u=u\rbr{\wc+h}-u$ and⋅$\Delta_{h}\Sigma=\Sigma\cap\rbr{-h+\Sigma}$,
the set on which $\Delta_{h}u$ is defined. Here and below, $D=-i\nabla$
is regarded as a row-vector operator, so that, for any map~$u$,
$Du$ is, up to the factor~$-i$, the Jacobian matrix of~$u$. We
abbreviate the norm on⋅$\smash{\Lambda_{p}^{m+\omega}}\rbr{\Sigma}^{n}=\smash{\Lambda_{p}^{m+\omega}}\rbr{\Sigma;\C^{n}}$
to⋅$\smash[b]{\abs{\wc}_{m+\omega,p,\Sigma}}$ and the seminorm on~$\Lambda_{p}^{\omega}\rbr{\Sigma}^{n}$
to~$\seminorm{\wc}_{\omega,p,\Sigma}$.  We write⋅$\smash{C^{m,\omega}\rbr{\wb{\Sigma};U}}$
for⋅$\smash{\Lambda_{\infty}^{m+\omega}}\rbr{\Sigma;U}$, in line
with the notation for H\"older spaces.  For⋅a⋅power modulus~$\omega\rbr t=t^{r}$,
the space~$\smash{\Lambda_{p}^{m+\omega}}\rbr{\Sigma;U}$ coincides
with the Sobolev space~$\smash{W_{p}^{m+r}}\rbr{\Sigma;U}$ if⋅$r=1$
and⋅$p>1$, and with the Lipschitz–Besov space~$\smash{\Lambda_{p}^{m+r}}\rbr{\Sigma;U}=\smash{B_{p,\infty}^{m+r}}\rbr{\Sigma;U}$
if⋅$r\in\rbr{0,1}$ (see,⋅\eg,~\cite{AF:2003} or⋅\cite{Agranovich:2013}
for its definition). In⋅particular, $\smash{\Lambda_{\infty}^{m+\omega}}\rbr{\Sigma;U}$
is the space~$C^{m,r}\rbr{\smash[t]{\wb{\Sigma}};U}$ of H\"{o}lder
continuous functions.

Let⋅$B\subset\R^{d}$ be the open unit ball  centered at the origin,
and let $B_{\pm}$ be the open unit half-balls with~$\pm x_{d}\in\rbr{0,1}$.
We say that a domain~$\Sigma\subset\R^{d}$  satisfies the uniform
weak Lipschitz condition if it is a⋅$d$\nobreakdash-dimensional
Lipschitz manifold with boundary embedded in~$\R^{d}$. More precisely,
there is a uniformly locally finite open covering~$\cbr{W_{k}}$
of⋅$\partial\Sigma$ and a family of bi-Lipschitz bijections~$w_{k}\colon W_{k}\to B$
such that (1)~$w_{k}\rbr{W_{k}\cap\Sigma}=B_{+}$ and $w_{k}\rbr{W_{k}\cap\partial\Sigma}=\partial B_{+}\!\setminus\partial B$;
(2)~$\sup_{k}\norm{Dw_{k}}_{\infty,W_{k}}$ and⋅$\sup_{k}\norm{Dw_{k}^{-1}}_{\infty,B}$
are finite; and (3)~for a fixed~$\delta>0$, every open ball~$B_{\delta}\rbr x$
with center~$x\in\partial\Sigma$ and radius~$\delta$ is contained
in some coordinate patch~$W_{k}$. The⋅last two conditions are automatically
satisfied provided that the boundary of⋅$\Sigma$ is compact. Notice
that the open set~$\R^{d}\setminus\wb{\Sigma}$, which need not be
connected, satisfies the uniformly weak Lipschitz condition whenever
$\Sigma$~does.

For⋅such~$\Sigma$, there exists a $C^{\infty}$\nobreakdash-partition
of unity~$\cbr{\varphi_{k}}$ subordinate to⋅$\cbr{W_{k}}$ with
the property that $\sup_{k}\norm{D^{\alpha}\varphi_{k}}_{\infty,W_{k}}$
is finite for any~$\alpha$, see~\cite[Chapter~6, Section~3]{St:1970}.
Then  $\Sigma$ admits a linear extension operator   mapping⋅$\smash[b]{W_{p}^{1}}\rbr{\Sigma}$
to⋅$\smash[b]{W_{p}^{1}}\rbr{\R^{d}}$ for every⋅$p\in[1,\infty)$
(see⋅the proof of Theorem~5.22 and Remark~5.23(2) in~\cite{AF:2003}).
It⋅follows that the standard density and embedding results for⋅$\R^{d}$
also hold for~$\Sigma$. In⋅particular, the Sobolev  theorem states
that $W_{p}^{1}\rbr{\Sigma}$ is continuously embedded in⋅$L_{q}\rbr{\Sigma}$
for any⋅$q\in\sbr{p,p^{*}}$. Here, $p^{*}$ is the Sobolev conjugate
to $p$, given by⋅$1/p^{*}=1/p-1/d$ if⋅$p<d$; $p^{*}$ is any finite
number greater than or equal to⋅$p$ if⋅$p=d$; and $p^{*}=\infty$
if~$p>d$.   There is also a bounded linear trace operator from⋅$\smash{W_{p}^{1}}\rbr{\Sigma;U}$
to⋅$\smash{W_{p}^{1-1/p}}\rbr{\partial\Sigma;U}$ provided that~$p>1$.
As⋅for⋅$\smash{\Lambda_{p}^{\omega}}$\nobreakdash-spaces with⋅$p\in\rbr{1,\infty}$,
 any linear extension operator~$\op E$ that maps⋅$L_{p}\rbr{\Sigma}$
and⋅$\smash{W_{p}^{1}}\rbr{\Sigma}$ (or⋅a⋅subspace of the latter)
continuously into the same spaces on⋅$\R^{d}$ satisfies, for⋅$t\in\sbr{0,1}$,
\begin{equation}
\osc_{p}\rbr{u;t}\le\osc_{p}\rbr{\op Eu;t}\lesssim\osc_{p}\rbr{u;t}+t\norm u_{p,\Sigma},\qquad u\in L_{p}\rbr{\Sigma},\label{eq: Notation | E and osc_p}
\end{equation}
with⋅a⋅constant depending on⋅$\op E$ (see~(6) in~\cite{Sharpley:1983}, where the result is stated
for homogeneous Sobolev spaces). In⋅particular, $\op E$ preserves
the  $\Lambda_{p}^{\omega}$\nobreakdash-regularity. For⋅$p=\nobreak\infty$,
the extension by reflection across the boundary has the same property.

Given⋅a⋅set~$G\subset\R^{d}$, we let $G_{\delta}$,⋅$\delta>0$,
denote the $\delta$\nobreakdash-neighborhood of~$G$:
\[
G_{\delta}=\bigcbr{x\in\R^{d}\colon\dist\rbr{x,G}<\delta}.
\]

We will often use the notation~$\alpha\lesssim\beta$  (which is
the same as saying that~$\beta\gtrsim\alpha$) to mean that there
is a positive constant~$C$ depending only on some fixed parameters
(specified in the statements of the results) such that~$\alpha\le C\beta$.
If⋅$\alpha\lesssim\beta$ and⋅$\alpha\gtrsim\beta$ simultaneously,
then⋅$\alpha\sim\beta$. Finally, $\alpha\meet\beta$ and⋅$\alpha\join\beta$
are, respectively, the smaller and the larger of $\alpha$ and~$\beta$.

\section{\label{sec: OO}Original operator}

Let⋅$\Omega$ be a (possibly unbounded) uniformly weakly Lipschitz
domain in⋅$\R^{d}$, and let $\Gamma$ be a⋅$\rbr{d-1}$\nobreakdash-dimensional
surface in⋅$\Omega$,   called the interface, that divides $\Omega$
into finitely many uniformly weakly Lipschitz subdomains~$\Omega_{i}$,⋅$i\in\set I$
(see⋅Figure~\ref{fig: OO | Domains}).

The⋅periodic structure on each⋅$\Omega_{i}$ is described by a
full-rank lattice~$\lattice_{i}$ and a small scale~$\varepsilon_{i}$,
which is a positive nondecreasing function on⋅$\R_{+}$ satisfying⋅$\varepsilon_{i}\rbr{\varepsilon}\to0$
as⋅$\varepsilon\to0$. Below, $\varepsilon_{\join}\rbr{\varepsilon}$
stands for the largest of the scales~$\varepsilon_{i}\rbr{\varepsilon}$.
Let⋅$\cell_{i}$ be an open primitive cell of⋅$\lattice_{i}$ centered
at the origin, and let $\cellradius_{i}=1/2\diam\cell_{i}$; without
loss of generality, $\abs{\cell_{i}}=1$. We
will often work with families~$\cbr{O_{i}}_{i\in\set I}$ of objects,
where $O_{i}$ is associated with⋅$\Omega_{i}$. The corresponding
piecewise object on⋅$\Omega\setminus\Gamma$ (whose precise meaning
will be clear from the context) is denoted by~$O_{\#}$, and relations
involving such objects are understood componentwise. One  example
is the function~$\varepsilon_{\#}$ defined by⋅$\varepsilon_{\#}\rbr{x,\varepsilon}=\varepsilon_{i}\rbr{\varepsilon}$
for every⋅$x\in\Omega_{i}$. To⋅lighten the notation, we suppress
the dependence of⋅$\varepsilon_{i}\rbr{\varepsilon}$ and⋅$\varepsilon_{\#}\rbr{x,\varepsilon}$
on⋅$\varepsilon$ and simply write⋅$\varepsilon_{i}$ and⋅$\varepsilon_{\#}\rbr x$.

Next, we introduce the scaling~$\tau_{\#}^{\varepsilon}$ that turns
a function of the slow and fast variables into the corresponding locally
periodic one. Namely, $\tau_{i}^{\varepsilon}$, $i\in\set I$, sends
$u\colon\Omega_{i}\times\R^{d}\to L_{0}\rbr{\cell_{i}}$ to the
function~$\tau_{i}^{\varepsilon}u\colon\Omega_{i}\to L_{0}\rbr{\cell_{i}}$
given by⋅$\tau_{i}^{\varepsilon}u\rbr x=u\rbr{x,\varepsilon_{i}^{-1}x}$.
We identify⋅$\rbr{u\rbr{x,y}}\rbr z$ with⋅$u\rbr{x,y,z}$, and so
on, choosing suitable representatives where necessary, so that
\begin{equation}
\tau_{\#}^{\varepsilon}u\rbr{x,z}=u\rbr{x,\varepsilon_{\#}^{-1}x,z}\label{def: OO | Operation tau^e}
\end{equation}
for⋅all⋅$x\in\smash[b]{\Omega_{\#}}$ and⋅$z\in\smash[b]{\cell_{\#}}$.
Clearly, $\smash[b]{\tau_{\#}^{\varepsilon}}$ is a homomorphism of
the respective algebras; in other words, for any two functions~$u$
and⋅$v$ from⋅$\Omega_{\#}\times\R^{d}$ to⋅$L_{0}\rbr{\cell_{\#}}$,
one has
\begin{equation}
\tau_{\#}^{\varepsilon}\rbr{u+v}=\tau_{\#}^{\varepsilon}u+\tau_{\#}^{\varepsilon}v,\qquad\tau_{\#}^{\varepsilon}uv=\tau_{\#}^{\varepsilon}u\cdot\tau_{\#}^{\varepsilon}v.\label{eq: OO | Operation tau^e is multiplicative}
\end{equation}
Here the dot~``$\cdot$'' stands for the pointwise product, and
we adopt the convention that an operator acts on everything to its
right up to the dot, \eg, $\smash[b]{\tau_{\#}^{\varepsilon}uv=\tau_{\#}^{\varepsilon}\rbr{uv}}$
and⋅$\tau_{\#}^{\varepsilon}u\cdot v=\rbr{\tau_{\#}^{\varepsilon}u}v$.
We set⋅$u^{\varepsilon}=\tau_{\#}^{\varepsilon}u$ for short.

Let⋅$A=\smash[t]{\cbr{A_{kl}}_{k,l=1}^{d}}$ be a tensor-valued
function whose entries satisfy
\begin{equation}
A_{kl}\in C^{0,\omega_{\#}}\rbr{\wb{\Omega}_{\#};\wt L_{\infty}\rbr{\cell_{\#}}}^{n\times n}.\label{assu: OO| Coefficient A regular}
\end{equation}
Thus, $A$ is continuous on⋅$\Omega$ away from the interface~$\Gamma$,
and the   modulus of continuity~$\omega_{i}$ controls its regularity
on~$\Omega_{i}$. Following⋅\cite[Lemma~5.6]{Al:1992}, we may
choose a representative~$A\colon\rbr{\Omega\setminus\Gamma}\times\R^{d}\to\B\rbr{\C^{n\times d}}$
and null sets~$\nullset_{\#}\subset\cell_{\#}$ such that, for each⋅$i\in\set I$,
 the map~$A\rbr{\wc,y}$ is  continuous on⋅$\Omega_{i}$ with
modulus of continuity~$\omega_{i}$ for every⋅$y\in\cell_{i}\setminus\nullset_{i}$,
while $A\rbr{x,\wc}$ is  $\lattice_{i}$\nobreakdash-periodic and
measurable on⋅$\R^{d}$ for every~$x\in\Omega_{i}$. In⋅particular,
the function~$A^{\varepsilon}\colon\Omega\to\B\rbr{\C^{n\times d}}$
is measurable and uniformly bounded.

We now set up the functional framework to accommodate general boundary
conditions. Fix⋅$p\in\rbr{1,\infty}$, and let $\set W_{p}^{1}\rbr{\Omega;\C^{n}}$
and⋅$\smash{\set W_{p^{\hconj}}^{1}\rbr{\Omega;\C^{n}}}$ be closed
subspaces of $\smash[t]{W_{p}^{1}}\rbr{\Omega}^{n}$ and⋅$\smash{W_{p^{\hconj}}^{1}\rbr{\Omega}^{n}}$,
respectively, that contain⋅$C_{c}^{\infty}\rbr{\Omega}^{n}$, so
that,⋅\eg,
\begin{equation}
\mathring{W}_{p}^{1}\rbr{\Omega}^{n}\subseteq\set W_{p}^{1}\rbr{\Omega;\C^{n}}\subseteq W_{p}^{1}\rbr{\Omega}^{n}.\label{def: OO | Curly W_p^1 space}
\end{equation}
Let⋅$\smash[t]{\set W_{p}^{-1}}\rbr{\Omega;\C^{n}}$ and⋅$\smash[b]{\Norm{\wc}_{-1,p,\Omega}}$
stand for the conjugate dual of⋅$\smash[b]{\set W_{p^{\hconj}}^{1}}\rbr{\Omega;\C^{n}}$
(under the $L_{2}$\nobreakdash-pairing) and the associated norm.
This dual is isometrically isomorphic to the quotient space⋅$\smash[t]{\rbr{W_{p^{\hconj}}^{1}\rbr{\Omega}^{n}}^{*}\!/\rbr{\set W_{p^{\hconj}}^{1}\rbr{\Omega;\C^{n}}}^{\bot}}$,
where $\smash[t]{\rbr{\set W_{p^{\hconj}}^{1}\rbr{\Omega;\C^{n}}}^{\bot}}$
is the subspace of  functionals in⋅$\smash[t]{\rbr{W_{p^{\hconj}}^{1}\rbr{\Omega}^{n}}^{*}}$
vanishing on⋅$\smash[t]{\set W_{p^{\hconj}}^{1}\rbr{\Omega;\C^{n}}}$.
We can thus regard the canonical projection
\begin{equation}
\pi\colon f\mapsto f+\rbr{\set W_{p^{\hconj}}^{1}\rbr{\Omega;\C^{n}}}^{\bot}\label{def: OO | Canonical projection for the dual}
\end{equation}
as⋅an⋅epimorphism of⋅$\rbr{W_{p^{\hconj}}^{1}\rbr{\Omega}^{n}}^{*}$
onto⋅$\set W_{p}^{-1}\rbr{\Omega;\C^{n}}$,⋅with
\begin{equation}
\Norm{\pi f}_{-1,p,\Omega}\le\norm f_{-1,p,\Omega}^{*}.\label{est: OO | Canonical projection bounded}
\end{equation}
We also introduce⋅$\smash{\set W_{p}^{1-1/p}}\rbr{\partial\Omega;\C^{n}}$,
the subspace of⋅$\smash{W_{p}^{1-1/p}}\rbr{\partial\Omega}^{n}$ consisting
of boundary traces of⋅$\smash{\set W_{p}^{1}}\rbr{\Omega;\C^{n}}$.
Similar notation is used for the space of boundary traces of~$\smash{\set W_{p^{\hconj}}^{1}}\rbr{\Omega;\C^{n}}$.
The⋅conjugate dual of⋅the latter trace space (again with respect
to the $L_{2}$\nobreakdash-pairing) is denoted by⋅$\smash[t]{\set W_{p}^{-1/p}}\rbr{\partial\Omega;\C^{n}}$
and  its norm by~$\smash[b]{\Norm{\wc}_{-1/p,p,\partial\Omega}}$.

We are now in a position to define our main object, the operator~$\op A^{\varepsilon}$.
Given a bounded linear map~$\sigma\colon\smash[t]{\set W_{p}^{1-1/p}}\rbr{\partial\Omega;\C^{n}}\to\smash[t]{\set W_{p}^{-1/p}}\rbr{\partial\Omega;\C^{n}}$,
 let $\op A^{\varepsilon}\colon\smash[t]{\set W_{p}^{1}}\rbr{\Omega;\C^{n}}\to\smash[t]{\set W_{p}^{-1}}\rbr{\Omega;\C^{n}}$
be the operator associated with the form
\begin{equation}
\begin{aligned}\rbr{\op A^{\varepsilon}u,u^{\prime}}_{\Omega}=\rbr{A^{\varepsilon}Du,Du^{\prime}}_{\Omega}+\rbr{\sigma\gamma u,\gamma u^{\prime}}_{\partial\Omega},\\
u\in\set W_{p}^{1}\rbr{\Omega;\C^{n}},u^{\prime}\in\set W_{p^{\hconj}}^{1}\rbr{\Omega;\C^{n}},
\end{aligned}
\label{def: OO | Operator Ae}
\end{equation}
where⋅$\gamma$ is the trace operator for~$\Omega$. If⋅$\smash[b]{\set W_{p}^{1}}\rbr{\Omega;\C^{n}}$
and⋅$\smash[b]{\set W_{p^{\hconj}}^{1}}\rbr{\Omega;\C^{n}}$ coincide
with the respective closures of⋅$C_{c}^{\infty}\rbr{\Omega}^{n}$,
then  (\ref{def: OO | Operator Ae}) defines the operator~$D^{*}A^{\varepsilon}D$
on⋅$\Omega$ with Dirichlet boundary condition (in⋅which case $\sigma$
plays no role). On⋅the⋅other hand, if $\smash[b]{\set W_{p}^{1}}\rbr{\Omega;\C^{n}}$
and⋅$\smash[b]{\set W_{p^{\hconj}}^{1}}\rbr{\Omega;\C^{n}}$ are $\smash[b]{W_{p}^{1}}\rbr{\Omega}^{n}$
and⋅$\smash[b]{W_{p^{\hconj}}^{1}}\rbr{\Omega}^{n}$, then  we arrive
at the operator~$D^{*}A^{\varepsilon}D$ on⋅$\Omega$ with the (nonlocal)
Robin-type boundary condition~$\partial_{\nu_{\mkern-1mu \op A^{\varepsilon}}}u+\sigma\gamma u=0$
on⋅$\partial\Omega$, where $\partial_{\nu_{\mkern-1mu \op A^{\varepsilon}}}u$
is the generalized weak conormal derivative of⋅$u$ relative to⋅$\op A^{\varepsilon}$.
 Mixed Dirichlet–Robin-type boundary conditions are allowed as well,
and, moreover, different components of a function in⋅$\smash[b]{\set W_{p}^{1}}\rbr{\Omega;\C^{n}}$
may satisfy different boundary conditions.

Obviously, $\op A^{\varepsilon}$ is bounded uniformly with respect
to⋅$\varepsilon>0$,⋅and
\begin{equation}
\Norm{\op A^{\varepsilon}u}_{-1,p,\Omega}\lesssim\bigrbr{\norm A_{L_{\infty}}+\norm{\sigma}_{\B}}\norm u_{1,p,\Omega},\qquad u\in\set W_{p}^{1}\rbr{\Omega;\C^{n}}.\label{est: OO | Operator Ae bounded}
\end{equation}
We suppose that, for all sufficiently small~$\varepsilon$, the operators~$\op A^{\varepsilon}$
have a common regular point and the resolvents are uniformly bounded
in~$\varepsilon$. That is, there are⋅$\mu\in\C$ and⋅$\varepsilon_{\mu}>0$
such that $\smash[b]{\op A_{\mu}^{\varepsilon}}=\op A^{\varepsilon}-\mu$
is an isomorphism for any⋅$\varepsilon\in\set E_{\mu}=(0,\varepsilon_{\mu}]$⋅and
\begin{equation}
\norm{\rbr{\op A_{\mu}^{\varepsilon}}^{-1}f}_{1,p,\Omega}\lesssim\Norm f_{-1,p,\Omega},\qquad f\in\set W_{p}^{-1}\rbr{\Omega;\C^{n}}.\label{assu: OO | Resolvent of Ae bounded}
\end{equation}
  Throughout the paper, we assume that⋅$\varepsilon_{\join}\rbr{\varepsilon_{\mu}}\le1$.

Let⋅$\smash[b]{\rbr{\op A_{\mu}^{\varepsilon}}^{\dual}}$ denote the
adjoint of~$\smash[b]{\op A_{\mu}^{\varepsilon}}$. The⋅objects
and results related to⋅$\smash[b]{\rbr{\op A_{\mu}^{\varepsilon}}^{\dual}}$
will be marked with~``$\dual$'' as well. Notice that  (\dualref{est: OO | Operator Ae bounded})
and⋅(\dualref{assu: OO | Resolvent of Ae bounded}) hold for⋅$\smash{\rbr{\op A_{\mu}^{\varepsilon}}^{\dual}}$
with $p^{\dual}=p^{\hconj}$ in place of⋅$p$ and, in fact, the same
constants.

For⋅sharper estimates away from the boundary and the interfaces,
we will need the following condition. Let $\chi\in\smash[t]{C^{0,1}\rbr{\wb{\Omega}}}$
be a cutoff function with⋅$\supp\chi\subset\Omega\setminus\nobreak\Gamma$.
We say that $\smash[b]{\op A_{\mu}^{\varepsilon}}$ satisfies an interior
energy estimate if there is $\smash{\chi^{\prime}}\in\smash[t]{C^{0,1}\rbr{\wb{\Omega}}}$
with⋅$\supp\chi\subset\supp\smash{\chi^{\prime}}\subset\Omega\setminus\Gamma$
such that
\begin{equation}
\norm{D\chi u}_{p,\Omega}\lesssim\norm u_{p,\Omega}+\norm{\chi^{\prime}\op A_{\mu}^{\varepsilon}u}_{-1,p,\Omega}^{*},\qquad u\in C_{c}^{\infty}\rbr{\Omega\setminus\Gamma}^{n}.\label{assu: OO | Interior energy estimate}
\end{equation}
For⋅compactly supported cutoff functions, this Caccioppoli-type inequality
enjoys a self-improvement property, as the next lemma shows.
\begin{lem}
\label{lem: OO | Self-improvement of the interior energy estimate}Suppose
that the inverse of⋅$\smash{\op A_{\mu}^{\varepsilon}}$ is also uniformly
bounded from⋅$\smash{\rbr{W_{q^{\hconj}}^{1}\rbr{\Omega}^{n}}^{*}}$
to⋅$\smash{W_{q}^{1}\rbr{\Omega}^{n}}$ for some⋅$q\in\rbr{1,p}$.
Assume further that, given any⋅$\chi\in\smash{C_{c}^{0,1}}\rbr{\Omega\setminus\Gamma}$,
there is $\smash[t]{\chi^{\prime}}\in\smash{C_{c}^{0,1}}\rbr{\Omega\setminus\Gamma}$,
with⋅$\supp\chi\subset\supp\smash[t]{\chi^{\prime}}$, such that
\begin{equation}
\norm{D\chi u}_{q,\supp\chi}\lesssim\norm u_{q,\supp\chi^{\prime}}+\norm{\chi^{\prime}\op A_{\mu}^{\varepsilon}u}_{-1,q,\Omega}^{*},\qquad u\in C_{c}^{\infty}\rbr{\Omega\setminus\Gamma}^{n}.\label{est: OO | Self-improvement of the interior energy estimate}
\end{equation}
Then the same estimate holds for⋅$p$ in place of⋅$q$, possibly
with a different~$\chi^{\prime}$.
\end{lem}
\begin{proof}
We⋅have $\sigma\gamma\chi=0$ and⋅$\rbr{\gamma\chi}^{*}\sigma=0$,
so the proof is no different from that for the case~$\sigma=0$
in~\cite[Lemma~3.2]{Senik:2023-1}.
\end{proof}

\section{\label{sec: EO}Effective operator}

We begin with the cell problem, whose solution is used to define
the effective operator.    For⋅$\zeta\in\smash[t]{\C^{n\times d}}$
and⋅$x\in\Omega_{\#}$, it reads
\begin{equation}
\op A\rbr x\hairspace N_{\zeta}\rbr{x,\wc}=-D_{2}^{*}A\rbr{x,\wc}\hairspace\zeta,\label{def: EO | Cell problem}
\end{equation}
where⋅$\smash[b]{N_{\zeta}}\rbr{x,\wc}\in\smash{\wt W_{p}^{1}}\rbr{\cell_{\#}}^{n}$
is the unknown function normalized to have zero mean over⋅$\cell_{\#}$,
and $\op A\rbr x\colon\smash{\wt W_{p}^{1}}\rbr{\cell_{\#}}^{n}\mkern-2mu /\mkern1mu \C^{n}\to\rbr{\smash{\wt W_{p^{\prime}}^{1}}\rbr{\cell_{\#}}^{n}\mkern-2mu /\mkern1mu \C^{n}}^{*}\simeq\smash{\wt W_{p}^{-1}}\rbr{\cell_{\#}}^{n}\cap\rbr{\C^{n}}^{\bot}$
is the operator given⋅by
\begin{equation}
\op A\rbr x=D_{2}^{*}A\rbr{x,\wc}\hairspace D_{2}.\label{def: EO | Operator A(x)}
\end{equation}
Here and below, $D_{1}$ and⋅$D_{2}$ stand for differentiation
in the first⋅(slow) and the second⋅(fast) variables, respectively.
 It⋅is natural to require that $\op A\rbr x$ be an isomorphism
for any⋅$x\in\Omega_{\#}$, with its inverse uniformly bounded in~$x$:
\begin{equation}
\norm{D_{2}\op A\rbr x^{-1}f}_{p,\cell_{\#}}\lesssim\norm f_{-1,p,\cell_{\#}},\qquad f\in\wt W_{p}^{-1}\rbr{\cell_{\#}}^{n}\cap\rbr{\C^{n}}^{\bot}.\label{assu: EO | Inverse of A(x) bounded}
\end{equation}
Then, for every $\zeta\in\C^{n\times d}$ and⋅$x\in\Omega_{\#}$,
the cell problem~(\ref{def: EO | Cell problem}) has a unique solution.
It⋅follows that the mapping~$\zeta\mapsto N_{\zeta}\rbr{x,\wc}$
is linear, \ie, $N_{\zeta}\rbr{x,\wc}=N\rbr{x,\wc}\zeta$ for a
function~$N\rbr{x,\wc}$. The⋅next lemma establishes the regularity
of~$N$.

\begin{lem}
\label{lm: EO | Function N regular}It⋅holds that
\begin{equation}
N\in C^{0,\omega_{\#}}\rbr{\wb{\Omega}_{\#};\wt W_{p}^{1}\rbr{\cell_{\#}}}.\label{eq: EO | Function N regular}
\end{equation}
\end{lem}
\begin{proof}
We have
\begin{equation}
N_{\zeta}\rbr{x,\wc}+\C^{n}=-\op A\rbr x^{-1}D_{2}^{*}A\rbr{x,\wc}\hairspace\zeta,\label{eq: EO | Function N as a solution}
\end{equation}
whence, by⋅(\ref{assu: EO | Inverse of A(x) bounded}),
\[
\norm{D_{2}N_{\zeta}\rbr{x,\wc}}_{p,\cell_{\#}}\lesssim\norm{D_{2}^{*}A\rbr{x,\wc}\hairspace\zeta}_{-1,p,\cell_{\#}},
\]
and⋅therefore $D_{2}N$ is uniformly bounded on⋅$\Omega_{\#}$⋅with
\[
\norm{D_{2}N}_{L_{\infty}\rbr{\Omega_{\#};\wt L_{p}\rbr{\cell_{\#}}}}\lesssim\norm A_{L_{\infty}}.
\]

To⋅see that $D_{2}N$ is actually uniformly continuous on each⋅$\Omega_{i}$
with modulus of continuity~$\omega_{i}$, we calculate, using⋅(\ref{eq: EO | Function N as a solution}),
that, for any⋅$h\in\R^{d}$ and⋅$x\in\Delta_{h}\Omega_{i}$,
\begin{equation}
\findiff_{h}N_{\zeta}\rbr{x,\wc}+\C^{n}=-\op A\rbr{x+h}^{-1}D_{2}^{*}\rbr{\findiff_{h}A\rbr{x,\wc}\cdot\rbr{D_{2}N_{\zeta}\rbr{x,\wc}+\zeta}},\label{eq: EO | Function Delta_h N as a solution}
\end{equation}
where⋅$\findiff_{h}=\op T_{h}-\op I$ and⋅$\op T_{h}u\rbr{x,y}=u\rbr{x+h,y}$
for all⋅$u\in L_{0}\rbr{\R^{d}\times\R^{d}}$. As⋅a⋅consequence,
\[
\norm{\findiff_{h}D_{2}N_{\zeta}\rbr{x,\wc}}_{p,\cell_{i}}\lesssim\norm{D_{2}^{*}\rbr{\findiff_{h}A\rbr{x,\wc}\cdot\rbr{D_{2}N_{\zeta}\rbr{x,\wc}+\zeta}}}_{-1,p,\cell_{i}},
\]
and⋅so
\[
\osc_{\infty}\rbr{\op r_{i}D_{2}N;\wc}\lesssim\osc_{\infty}\rbr{\op r_{i}A;\wc}\hairspace\norm{I+D_{2}N}_{L_{\infty}\rbr{\Omega_{i};\wt L_{p}\rbr{\cell_{i}}}},
\]
where⋅$\op r_{i}$ is the operator of restriction to~$\Omega_{i}$.
This, together with (\ref{assu: OO| Coefficient A regular}) and the
Poincar\'{e} inequality, yields⋅(\ref{eq: EO | Function N regular}).
\end{proof}

\begin{rem}
\label{Rem: EO | Weaker ellipticity-type condition}The assumption~(\ref{assu: EO | Inverse of A(x) bounded})
implies the following ellipticity-type condition on the coefficient~$A$,
which is weaker than strong ellipticity:
\begin{equation}
\abs{\abr{A\rbr{\cdot,\cdot}\hairspace\eta\otimes\xi,\xi}}\gtrsim\abs{\eta}\hairspace\abs{\xi}^{2},\qquad\eta\in\C^{n},\xi\in\R^{d}.\label{est: EO | Weaker ellipticity-type condition}
\end{equation}
This is a special case of Lemma~\ref{lem: Appendix | Weaker ellipticity condition}
with~$\Sigma=\cell_{\#}$.
\end{rem}
We introduce the effective coefficient~$A^{0}\colon\Omega\setminus\Gamma\to\B\rbr{\C^{n\times d}}$
by setting
\begin{equation}
A^{0}\rbr x=\int_{\cell_{\#}}A\rbr{x,y}\hairspace\rbr{I+D_{2}N\rbr{x,y}}\dd y.\label{def: EO | Coefficient A0}
\end{equation}
It⋅is⋅straightforward to check that
\begin{equation}
\norm{A^{0}}_{L_{\infty}}\le\norm A_{L_{\infty}}\norm{I+D_{2}N}_{L_{\infty}\rbr{\Omega_{\#};\wt L_{p}\rbr{\cell_{\#}}}}\label{est: EO | Coefficient A0 bounded}
\end{equation}
and
\begin{equation}
A^{0}\in C^{0,\omega_{\#}}\rbr{\wb{\Omega}_{\#}}.\label{est: EO | Coefficient A0 regular}
\end{equation}
The⋅effective operator~$\op A^{0}\colon\set W_{p}^{1}\rbr{\Omega;\C^{n}}\to\set W_{p}^{-1}\rbr{\Omega;\C^{n}}$
 is now defined, just as⋅$\op A^{\varepsilon}$, through its form:
\begin{equation}
\begin{aligned}\rbr{\op A^{0}u,u^{\prime}}_{\Omega}=\rbr{A^{0}Du,Du^{\prime}}_{\Omega}+\rbr{\sigma\gamma u,\gamma u^{\prime}}_{\partial\Omega},\\
u\in\set W_{p}^{1}\rbr{\Omega;\C^{n}},u^{\prime}\in\set W_{p^{\hconj}}^{1}\rbr{\Omega;\C^{n}}.
\end{aligned}
\label{def: EO | Operator A0}
\end{equation}
By⋅(\ref{est: EO | Coefficient A0 bounded}), $\op A^{0}$ is 
bounded. We suppose that $\op A_{\mu}^{0}=\op A^{0}-\mu$ is invertible,
with
\begin{equation}
\norm{\rbr{\op A_{\mu}^{0}}^{-1}f}_{1,p,\Omega}\lesssim\Norm f_{-1,p,\Omega},\qquad f\in\set W_{p}^{-1}\rbr{\Omega;\C^{n}},\label{assu: EO | Resolvent of A0 bounded}
\end{equation}
and,⋅moreover,
\begin{equation}
\seminorm{D\rbr{\op A_{\mu}^{0}}^{-1}f}_{\omega_{\#},p,\Omega_{\#}}\lesssim\norm f_{p,\Omega},\qquad f\in L_{p}\rbr{\Omega}^{n}.\label{assu: EO | Resolvent of A0 regular}
\end{equation}

\begin{rem}
 Usually, the isomorphism property~(\ref{assu: EO | Resolvent of A0 bounded})
follows from~(\ref{assu: OO | Resolvent of Ae bounded}). On⋅the⋅other
hand, the assumption~(\ref{assu: EO | Resolvent of A0 regular})
requires sufficient regularity of the coefficients, the domain, the
interface and the boundary conditions. We refer to Section~\ref{sec: Examples}
as well as to⋅\cite[Section~7]{Senik:2023-1} for examples.
\end{rem}
\begin{rem}
The⋅piecewise regularity in⋅(\ref{assu: EO | Resolvent of A0 regular})
can sometimes imply the corresponding global regularity. For⋅example,
if $\omega$ does not vanish too fast at~$0$, as for a power modulus~$\omega\rbr t=t^{s}$
with⋅$s<1/p$, then the extension by zero is bounded in the  Lipschitz–Besov
spaces~$\smash{\Lambda_{p}^{s}}$  (in⋅fact, only for⋅$s<1/p$;
see~\cite[Subsection~4.6.3, Theorem~1]{Runst+Sickel:1996}). Therefore,
piecewise $\smash{\Lambda_{p}^{1+s}}$\nobreakdash-regularity together
with global $W_{p}^{1}$\nobreakdash-regularity is equivalent to
global $\smash{\Lambda_{p}^{1+s}}$\nobreakdash-regularity. If, on
the contrary, $\omega\rbr t/\smash[t]{t^{1/p}}\to0$ as⋅$t\to0$,
such an implication fails in general, since a nontrivial jump of
the gradient across an interface is no more regular than~$\Lambda_{p}^{\smash[b]{1/p}}$.
\end{rem}
The⋅above construction applies equally to⋅$\rbr{\op A^{\varepsilon}}^{\dual}$,
since (\ref{assu: EO | Inverse of A(x) bounded}) holds if and only
if⋅(\dualref{assu: EO | Inverse of A(x) bounded}) does. The⋅operator~$\rbr{\op A^{0}}^{\dual}$
obtained in this way, \ie, by⋅(\dualref{def: EO | Coefficient A0})
and⋅(\dualref{def: EO | Operator A0}), turns out to be the adjoint
of the effective operator~$\op A^{0}$. Hence (\dualref{assu: EO | Resolvent of A0 bounded})
also holds. The⋅regularity property for⋅$\rbr{\op A^{0}}^{\dual}$,
however, is a separate requirement, and the moduli in⋅(\ref{assu: EO | Resolvent of A0 regular})
and⋅(\dualref{assu: EO | Resolvent of A0 regular}) need not
coincide. Indeed,  the regularity is typically limited by a threshold
depending on the Lebesgue exponent, which is $p$ in⋅(\ref{assu: EO | Resolvent of A0 regular})
and⋅$p^{\prime}$ in~(\dualref{assu: EO | Resolvent of A0 regular}).
For⋅instance,  in mixed Dirichlet–Neumann problems, the threshold~$1/2+2/p$
for the smoothness exponent arises near the points where the boundary
condition changes, see~\cite[Theorem~3.1]{Shamir:1968}. Similarly,
in the Dirichlet problem on a Lipschitz domain, the threshold is
$1+1/p$ for⋅$p$ up to some exponent greater than~$2$ depending
only on the Lipschitz character (with a more involved threshold beyond⋅it),
and for all⋅$p$ if the domain is~$C^{1}$, see~\cite[Theorem~1.1]{Jerison+Kenig:1995}.
  Accordingly, whenever (\dualref{assu: EO | Resolvent of A0 regular})
is assumed, it comes with its own modulus~$\smash{\omega_{\#}^{\dual}}$,
which may differ from~$\smash{\omega_{\#}}$. We always take $\smash{\omega_{\#}^{\dual}}$
so that⋅$\smash{\omega_{\#}}\le\smash{\omega_{\#}^{\dual}}$; as the
proofs will show, this entails no loss of generality.

\section{\label{sec: C}Corrector}

The corrector is constructed in two stages. We first define and analyze
two-scale correctors associated with each⋅$\Omega_{i}$ and then glue
them together.

\subsection{Two-scale corrector}

For⋅a⋅fixed $i\in\set I$, let $A_{i}$ and⋅$N_{i}$ be extensions
of⋅$\op r_{i}A$ and⋅$\op r_{i}N$ to⋅$\R^{d}$, \eg,  by reflection
in⋅$\partial\Omega_{i}$, such that
\begin{equation}
\begin{aligned}A_{i} & \in C^{0,\omega_{i}}\rbr{\wb{\R}^{d};\wt L_{\infty}\rbr{\cell_{i}}},\\
N_{i} & \in C^{0,\omega_{i}}\rbr{\wb{\R}^{d};\wt W_{p}^{1}\rbr{\cell_{i}}}.
\end{aligned}
\label{eq: C | Functions Ai and Ni regular}
\end{equation}
We also choose a bounded extension operator~$\op e_{i}$ that
maps⋅$L_{p}\rbr{\Omega_{i}}$ to⋅$L_{p}\rbr{\R^{d}}$ and⋅$\smash[b]{W_{p}^{1}}\rbr{\Omega_{i}}$
to⋅$\smash[b]{W_{p}^{1}}\rbr{\R^{d}}$, and define⋅$\op R_{\mu,i}^{0}\colon L_{p}\rbr{\Omega}^{n}\to W_{p}^{1}\rbr{\Omega_{i}}^{n}$⋅as
\begin{equation}
\op R_{\mu,i}^{0}=\op r_{i}\rbr{\op A_{\mu}^{0}}^{-1}.\label{def: C |Resolvent R0_(i,mu)}
\end{equation}
The⋅two-scale corrector~$\op K_{\mu,i}\colon L_{p}\rbr{\Omega}^{n}\to L_{p}\rbr{\R^{d};\wt W_{p}^{1}\rbr{\cell_{i}}}^{n}$
is then given⋅by
\begin{equation}
\op K_{\mu,i}=N_{i}\op e_{i}D_{1}\op R_{\mu,i}^{0}.\label{def: C | Two-scale corrector K_(i,mu)}
\end{equation}
Both $\op R_{\mu,i}^{0}$ and⋅$\op K_{\mu,i}$ are bounded, as can
be seen from⋅(\ref{eq: EO | Function N regular}) and~(\ref{assu: EO | Resolvent of A0 bounded}).
\begin{rem}
\label{rem: Cor | Anisotropic norms}Throughout, we will use the notation~$\smash{W_{p}^{1_{1}}}\rbr{\R^{d}\times\cell_{i}}$
and⋅$\smash{W_{p}^{1_{2}}}\rbr{\R^{d}\times\cell_{i}}$, \etc,
to denote the anisotropic Sobolev spaces~$W_{p}^{1}\rbr{\R^{d};L_{p}\rbr{\cell_{i}}}$
and⋅$W_{p}^{1}\rbr{\cell_{i};L_{p}\rbr{\R^{d}}}$, respectively, and
write⋅$\norm{\wc}_{1_{1},p,\R^{d}\times\cell_{i}}$ and⋅$\norm{\wc}_{1_{2},p,\R^{d}\times\cell_{i}}$
for the associated norms. In⋅addition, $\seminorm{\wc}_{\omega,p,\R^{d}\times\cell_{i}}$
stands for the seminorm on~$\Lambda_{p}^{\omega}\rbr{\R^{d};L_{p}\rbr{\cell_{i}}}$.
\end{rem}
\begin{lem}
\label{lem: C | Operators R0_(mu,i) and K_(mu,i) bounded}For⋅any⋅$f\in L_{p}\rbr{\Omega}^{n}$,
we have
\begin{align}
\norm{\op R_{\mu,i}^{0}f}_{1,p,\Omega_{i}} & \lesssim\norm f_{p,\Omega},\label{est: C | Operator R0_(mu,i) bounded}\\
\norm{\op K_{\mu,i}f}_{1_{2},p,\R^{d}\times\cell_{i}} & \lesssim\norm f_{p,\Omega}.\label{est: C | Operator K_(mu,i) bounded}
\end{align}
\end{lem}
What is more, we can estimate the operator moduli of continuity of⋅$\op R_{\mu,i}^{0}$
and~$\op K_{\mu,i}$.
\begin{lem}
\label{lem: C | Operators R0_(i,mu) and K_(i,mu) regular}For⋅any⋅$f\in L_{p}\rbr{\Omega}^{n}$,
we have
\begin{align}
\seminorm{\op e_{i}D\op R_{\mu,i}^{0}f}_{\omega_{i},p,\R^{d}} & \lesssim\norm f_{p,\Omega},\label{est: C | Operator R0_(i,mu) regular}\\
\seminorm{D_{2}\op K_{\mu,i}f}_{\omega_{i},p,\R^{d}\times\cell_{i}}+\seminorm{\op K_{\mu,i}f}_{\omega_{i},p,\R^{d}\times\cell_{i}} & \lesssim\norm f_{p,\Omega},\label{est: C | Operator K_(i,mu) regular}\\
\seminorm{A_{i}\rbr{\op e_{i}D_{1}\op R_{\mu,i}^{0}+D_{2}\op K_{\mu,i}}f}_{\omega_{i},p,\R^{d}\times\cell_{i}} & \lesssim\norm f_{p,\Omega}.\label{est: C | Operator F_(i,mu) regular}
\end{align}
\end{lem}
\begin{proof}
If⋅$t>1$, then $\omega_{i}\rbr t\ge1$ and the $L_{p}$\nobreakdash-modulus
of continuity is at most twice the $L_{p}$\nobreakdash-norm, so
the estimates follow from Lemma~\ref{lem: C | Operators R0_(mu,i) and K_(mu,i) bounded}.
For⋅$t\le1$, recall from⋅(\ref{eq: Notation | E and osc_p}) that
\[
\osc_{p}\rbr{\op e_{i}u;t}\lesssim\osc_{p}\rbr{u;t}+t\norm u_{p,\Omega_{i}},\qquad u\in L_{p}\rbr{\Omega_{i}}.
\]
The⋅first estimate then follows from⋅(\ref{assu: EO | Resolvent of A0 bounded})
and~(\ref{assu: EO | Resolvent of A0 regular}). To⋅obtain the second,
we write
\begin{equation}
\findiff_{h}\op K_{\mu,i}=\findiff_{h}N_{i}\cdot\op e_{i}D_{1}\op R_{\mu,i}^{0}+\op T_{h}N_{i}\cdot\findiff_{h}\op e_{i}D_{1}\op R_{\mu,i}^{0},\label{eq: Cor | Operator Delta_h K_(i,mu)}
\end{equation}
whence
\[
\osc_{p}\rbr{\op K_{\mu,i}f;t}\le\osc_{\infty}\rbr{N_{i};t}\hairspace\norm{\op e_{i}D\op R_{\mu,i}^{0}f}_{p,\R^{d}}+\norm{N_{i}}_{L_{\infty}\rbr{\R^{d};\wt L_{p}\rbr{\cell_{i}}}}\osc_{p}\rbr{\op e_{i}D\op R_{\mu,i}^{0}f;t}.
\]
A⋅similar bound for⋅$D_{2}\op K_{\mu,i}f$ is obtained by differentiating⋅(\ref{eq: Cor | Operator Delta_h K_(i,mu)}).
Now⋅(\ref{eq: C | Functions Ai and Ni regular}), (\ref{est: C | Operator R0_(mu,i) bounded})
and⋅(\ref{est: C | Operator R0_(i,mu) regular}) yield~(\ref{est: C | Operator K_(i,mu) regular}).
The⋅same reasoning gives~(\ref{est: C | Operator F_(i,mu) regular}).
\end{proof}
\begin{rem}
In⋅the⋅definition of⋅$\smash[b]{\op K_{\mu,i}}$, we apply the extension
operator~$\op e_{i}$ to the derivative of⋅$\smash[b]{\op R_{\mu,i}^{0}}$
to ensure that the extension  does not affect the corresponding
modulus of continuity. In⋅certain cases, one may instead extend⋅$\smash[b]{\op R_{\mu,i}^{0}}$,
 provided that $\op e_{i}$ also maps $\smash[b]{W_{p}^{2}}\rbr{\Omega_{i}}$
to⋅$\smash[b]{W_{p}^{2}}\rbr{\R^{d}}$. For⋅instance, this is the
case if $\omega_{i}$ has a special form, such as⋅$\omega_{i}\rbr t=t^{s}$,
when our arguments can be easily reformulated in terms of the related
Sobolev seminorms, see~\cite{Senik:2023-1}.  We also mention here
the paper~\cite{Fef:2009}, where an extension operator preserving
the modulus of continuity of the highest-order derivatives was constructed
for $C^{m}$-spaces.
\end{rem}
The⋅regularity of the two-scale corrector~$\op K_{\mu,i}$ is not
sufficient for our purposes, since the corrector we need  must
map into⋅$W_{p}^{1}\rbr{\R^{d}}^{n}$. To⋅this end,  we will mollify⋅$\op K_{\mu,i}$
on each⋅$\Omega_{i}$, apply Steklov smoothing and then cut it off
near⋅$\partial\Omega_{i}\cap\Gamma$ to eliminate discontinuities
across the interface.

\subsection{\label{subsec: Mollification}Mollification}

For⋅each⋅$i\in\set I$,  let $J_{i}\in C_{c}^{0,1}\rbr{\cell_{i}}$
be a nonnegative even function satisfying⋅$\smash{\int_{\R^{d}}}J_{i}\rbr z\dd z=1$.
By⋅$\op J_{i}^{\varepsilon}$, we denote the mollification in the
first variable restricted to⋅$L_{1}\rbr{\R^{d}\times\cell_{i}}+L_{\infty}\rbr{\R^{d}\times\cell_{i}}$:
\begin{equation}
\op J_{i}^{\varepsilon}u\rbr{x,y}=\int_{\R^{d}}J_{i}\rbr{\varepsilon_{i}^{-1}x-w}\hairspace u\rbr{\varepsilon_{i}w,y}\dd w=\int_{\cell_{i}}J_{i}\rbr z\hairspace u\rbr{x+\varepsilon_{i}z,y}\dd z,\label{def: C | Mollification Je}
\end{equation}
where⋅$\rbr{x,y}\in\R^{d}\times\cell_{i}$. As⋅is well known, 
for⋅$u\in L_{p}\rbr{\R^{d}}$, the function~$\op J_{i}^{\varepsilon}u$
belongs to⋅$\smash[b]{W_{p}^{1}}\rbr{\R^{d}}$, yet its $\smash[b]{W_{p}^{1}}$\nobreakdash-norm
in general blows up as~$\varepsilon\to0$. At⋅the⋅same time, $\op J_{i}^{\varepsilon}u$
converges to⋅$u$ in~$L_{p}\rbr{\R^{d}}$. The⋅following two lemmas
bound the rates of blow-up and convergence in terms of the $L_{p}$\nobreakdash-modulus
of continuity of~$u$.
\begin{lem}
\label{lem: C | Mollification Je and regularity}Let⋅$q\in\sbr{1,\infty}$.
Then, for any⋅$\varepsilon>0$ and any⋅$u\in L_{q}\rbr{\R^{d}\times\cell_{i}}$,
\begin{equation}
\norm{D_{1}\op J_{i}^{\varepsilon}u}_{q,\R^{d}\times\cell_{i}}\lesssim\varepsilon_{i}^{-1}\osc_{q}\rbr{u;\varepsilon_{i}},\label{est: C | Bound for D Je}
\end{equation}
where⋅$\osc_{q}\rbr{u;\wc}$ is the $L_{q}$\nobreakdash-modulus
of continuity of⋅$u$ regarded as a function on⋅$\R^{d}$ with values
in~$L_{q}\rbr{\cell_{i}}$.
\end{lem}
\begin{proof}
Notice that
\[
\begin{aligned}D_{1}\op J_{i}^{\varepsilon}u\rbr{x,y} & =\varepsilon_{i}^{-1}\int_{\R^{d}}\rbr{DJ_{i}}\rbr{\varepsilon_{i}^{-1}x-w}\hairspace\findiff_{-x+\varepsilon_{i}w}u\rbr{x,y}\dd w\\
 & =-\varepsilon_{i}^{-1}\int_{\cell_{i}}\rbr{DJ_{i}}\rbr z\hairspace\findiff_{\varepsilon_{i}z}u\rbr{x,y}\dd z.
\end{aligned}
\]
Applying H\"{o}lder's inequality, we see that
\[
\norm{D_{1}\op J_{i}^{\varepsilon}u\rbr{\wc,y}}_{q,\R^{d}}\le\varepsilon_{i}^{-1}\norm{DJ_{i}}_{q^{\hconj},\cell_{i}}\biggrbr{\int_{\cell_{i}}\norm{\findiff_{\varepsilon_{i}z}u\rbr{\wc,y}}_{q,\R^{d}}^{q}\dd z}^{1/q},
\]
and⋅so, by taking the $L_{q}\rbr{\cell_{i}}$-norm in~$y$,
\[
\norm{D_{1}\op J_{i}^{\varepsilon}u}_{q,\R^{d}\times\cell_{i}}\lesssim\varepsilon_{i}^{-1}\biggrbr{\int_{\cell_{i}}\norm{\findiff_{\varepsilon_{i}z}u}_{q,\R^{d}\times\cell_{i}}^{q}\dd z}^{1/q}.
\]
For⋅$z\in\cell_{i}$, the shift satisfies $\abs{\varepsilon_{i}z}\le\varepsilon_{i}\cellradius_{i}$,
and therefore
\[
\norm{\findiff_{\varepsilon_{i}z}u}_{q,\R^{d}\times\cell_{i}}\le\osc_{q}\rbr{u;\varepsilon_{i}\cellradius_{i}}.
\]
If⋅$\cellradius_{i}\le1$, the lemma follows by monotonicity of~$\osc_{q}\rbr{u;\wc}$.
Otherwise, subadditivity of⋅$\osc_{q}\rbr{u;\wc}$ over⋅$\R^{d}$
gives $\osc_{q}\rbr{u;\lambda t}\le\rbr{\lambda+1}\osc_{q}\rbr{u;t}$
for $t\ge0$ and $\lambda\ge1$, which completes the proof.
\end{proof}
\begin{lem}
\label{lem: C | Bound for (Je - I)}Let⋅$q\in\sbr{1,\infty}$. Then,
for any⋅$\varepsilon>0$ and any⋅$u\in L_{q}\rbr{\R^{d}\times\cell_{i}}$,
\begin{equation}
\norm{\rbr{\op J_{i}^{\varepsilon}-\op I}u}_{q,\R^{d}\times\cell_{i}}\lesssim\osc_{q}\rbr{u;\varepsilon_{i}},\label{est: C | Bound for (Je - I)}
\end{equation}
where⋅$\osc_{q}\rbr{u;\wc}$ is the $L_{q}$\nobreakdash-modulus
of continuity of⋅$u$ regarded as a function on⋅$\R^{d}$ with values
in~$L_{q}\rbr{\cell_{i}}$.
\end{lem}
\begin{proof}
We write
\[
\rbr{\op J_{i}^{\varepsilon}-\op I}u\rbr{x,y}=\int_{\cell_{i}}J_{i}\rbr z\hairspace\findiff_{\varepsilon_{i}z}u\rbr{x,y}\dd z
\]
and⋅then proceed as in the proof of Lemma~\ref{lem: C | Mollification Je and regularity}.
\end{proof}
We now replace the operators~$D\op R_{\mu,i}^{0}$ and⋅$\op K_{\mu,i}$
by their mollified versions. The⋅two lemmas below collect the properties
of these regularizations that we need. Recall the notation for anisotropic
norms introduced in Remark~\ref{rem: Cor | Anisotropic norms}.
\begin{lem}
\label{lem: C | Operators Je_i R0_(mu,i) and Je_i K_(mu,i) bounded}For⋅any⋅$\varepsilon\in\set E_{\mu}$
and⋅$f\in L_{p}\rbr{\Omega}^{n}$, we have
\begin{align}
\norm{D\op J_{i}^{\varepsilon}\op e_{i}D\op R_{\mu,i}^{0}f}_{p,\R^{d}} & \lesssim\varepsilon_{i}^{-1}\omega_{i}\rbr{\varepsilon_{i}}\hairspace\norm f_{p,\Omega},\label{est: C | Operator Je_i R0_(mu,i) bounded}\\
\norm{D_{1}\op J_{i}^{\varepsilon}\op K_{\mu,i}f}_{1_{2},p,\R^{d}\times\cell_{i}} & \lesssim\varepsilon_{i}^{-1}\omega_{i}\rbr{\varepsilon_{i}}\hairspace\norm f_{p,\Omega},\label{est: C |Operator Je_i K_(mu,i) bounded}\\
\norm{D_{1}\op J_{i}^{\varepsilon}A_{i}\rbr{\op e_{i}D_{1}\op R_{\mu,i}^{0}+D_{2}\op K_{\mu,i}}f}_{p,\R^{d}\times\cell_{i}} & \lesssim\varepsilon_{i}^{-1}\omega_{i}\rbr{\varepsilon_{i}}\hairspace\norm f_{p,\Omega}.
\end{align}

\end{lem}
\begin{lem}
\label{lem: C | Operators Je_i R0_(i,mu) and Je_i K_(i,mu) converge}For⋅any⋅$\varepsilon\in\set E_{\mu}$
and⋅$f\in L_{p}\rbr{\Omega}^{n}$, we have
\begin{align}
\norm{\rbr{\op J_{i}^{\varepsilon}-\op I}\op e_{i}D\op R_{\mu,i}^{0}f}_{p,\R^{d}} & \lesssim\omega_{i}\rbr{\varepsilon_{i}}\hairspace\norm f_{p,\Omega},\label{est: C | Operator Je_i R0_(i,mu) converges}\\
\norm{\rbr{\op J_{i}^{\varepsilon}-\op I}\op K_{\mu,i}f}_{1_{2},p,\R^{d}\times\cell_{i}} & \lesssim\omega_{i}\rbr{\varepsilon_{i}}\hairspace\norm f_{p,\Omega},\label{est: C | Operator Je_i K_(i,mu) converges}\\
\norm{\rbr{\op J_{i}^{\varepsilon}-\op I}A_{i}\rbr{\op e_{i}D_{1}\op R_{\mu,i}^{0}+D_{2}\op K_{\mu,i}}f}_{p,\R^{d}\times\cell_{i}} & \lesssim\omega_{i}\rbr{\varepsilon_{i}}\hairspace\norm f_{p,\Omega}.
\end{align}
\end{lem}
Both lemmas are immediate consequences of Lemmas~\ref{lem: C | Mollification Je and regularity},
\ref{lem: C | Bound for (Je - I)} and~\ref{lem: C | Operators R0_(i,mu) and K_(i,mu) regular}.

A⋅further regularization is needed to   separate the fast and slow
variables in the corrector. For this purpose, we use the Steklov
smoothing operator.

\subsection{Steklov smoothing}

Let $\op T_{i}^{\varepsilon}\colon L_{0}\rbr{\R^{d}\times\cell_{i}}\to L_{0}\rbr{\R^{d}\times\cell_{i};L_{0}\rbr{\cell_{i}}}$,
$i\in\set I$, be the translation operator
\begin{equation}
\op T_{i}^{\varepsilon}u\rbr{x,y,z}=u\rbr{x+\varepsilon_{i}z,y},\label{def: C | Operator Te}
\end{equation}
where⋅$\rbr{x,y}\in\R^{d}\times\cell_{i}$ and~$z\in\cell_{i}$.
 Obviously, $\op T_{i}^{\varepsilon}\rbr{u+v}=\op T_{i}^{\varepsilon}u+\op T_{i}^{\varepsilon}v$
and⋅$\op T_{i}^{\varepsilon}uv=\op T_{i}^{\varepsilon}u\cdot\op T_{i}^{\varepsilon}v$,
so $\op T_{i}^{\varepsilon}$ is an algebra homomorphism. Next, the
formal adjoint of⋅$\op T_{i}^{\varepsilon}$ with respect to the
$L_{2}$\nobreakdash-pairing is given by the formula
\[
\rbr{\op T_{i}^{\varepsilon}}^{*}v\rbr{x,y}=\int_{\cell_{i}}v\rbr{x-\varepsilon_{i}z,y,z}\dd z.
\]
 The⋅Steklov smoothing operator~$\op S_{i}^{\varepsilon}$ associated
with⋅$\Omega_{i}$ is defined to be the restriction of⋅$\rbr{\op T_{i}^{\varepsilon}}^{*}$
to⋅$L_{1}\rbr{\R^{d}\times\cell_{i}}+L_{\infty}\rbr{\R^{d}\times\cell_{i}}$.
Since $\cell_{i}$ is symmetric about the origin, we may write
\begin{equation}
\op S_{i}^{\varepsilon}u\rbr{x,y}=\int_{\cell_{i}}\op T_{i}^{\varepsilon}u\rbr{x,y,z}\dd z.\label{def: C | Steklov smoothing}
\end{equation}
For⋅the⋅same reason, $\op S_{i}^{\varepsilon}$ is  formally self-adjoint.
\begin{rem}
\label{rem: C | Mollification Je_i other expression}It⋅will be useful
to express the operator~$\op J_{i}^{\varepsilon}$ similarly to~$\op S_{i}^{\varepsilon}$:
\[
\op J_{i}^{\varepsilon}u\rbr{x,y}=\int_{\cell_{i}}J_{i}\rbr z\hairspace\op T_{i}^{\varepsilon}u\rbr{x,y,z}\dd z.
\]
From⋅this point of view, the Steklov smoothing operator is a particular
case of~$\op J_{i}^{\varepsilon}$, with⋅$J_{i}$ being the  characteristic
function of~$\cell_{i}$. Unlike a mollifier, this kernel is discontinuous.
\end{rem}
We recall some properties of⋅$\op T_{i}^{\varepsilon}$ and⋅$\op S_{i}^{\varepsilon}$;
see~\cite[Subsection~5.2]{Senik:2023-1} for the proofs. The⋅first
three concern the operators~$\tau_{i}^{\varepsilon}\op T_{i}^{\varepsilon}$
and~$\tau_{i}^{\varepsilon}\op S_{i}^{\varepsilon}$.
\begin{lem}
\label{lem: C | Operator =0003C4e Te isometry on Rd}For⋅any⋅$q\in\sbr{1,\infty}$
and⋅$\varepsilon>0$, $\tau_{i}^{\varepsilon}\op T_{i}^{\varepsilon}$
is an isometry of⋅$\wt L_{q}\rbr{\R^{d}\times\cell_{i}}$ into~$L_{q}\rbr{\R^{d};L_{q}\rbr{\cell_{i}}}$.
\end{lem}

Let $u$ be a function on a neighborhood of a domain~$\Sigma\subset\R^{d}$.
Extending it by zero to⋅$\R^{d}$, applying Lemma~\ref{lem: C | Operator =0003C4e Te isometry on Rd}
and using that $\Sigma+\varepsilon_{i}\cell_{i}\subset\Sigma_{\varepsilon_{i}\cellradius_{i}}$,
we obtain the following results.
\begin{lem}
\label{lem: C | Operator =0003C4e Te bounded on domains}Let⋅$\Sigma$
be a domain in⋅$\R^{d}$, and let⋅$q\in\sbr{1,\infty}$. Then, for
any⋅$\varepsilon>0$ and⋅$u\in\wt L_{q}\rbr{\Sigma_{\varepsilon_{i}\cellradius_{i}}\times\cell_{i}}$,
\begin{equation}
\norm{\tau_{i}^{\varepsilon}\op T_{i}^{\varepsilon}u}_{q,\Sigma\times\cell_{i}}\le\norm u_{q,\Sigma_{\varepsilon_{i}\cellradius_{i}}\times\cell_{i}}.
\end{equation}
\end{lem}
\begin{lem}
\label{lem: C | Operator =0003C4e Se bounded on domains}Let⋅$\Sigma$
be a domain in⋅$\R^{d}$, and let⋅$q\in\sbr{1,\infty}$. Then, for
any⋅$\varepsilon>0$ and⋅$u\in\wt L_{q}\rbr{\Sigma_{\varepsilon_{i}\cellradius_{i}}\times\cell_{i}}$,
\begin{equation}
\norm{\tau_{i}^{\varepsilon}\op S_{i}^{\varepsilon}u}_{q,\Sigma}\le\norm u_{q,\Sigma_{\varepsilon_{i}\cellradius_{i}}\times\cell_{i}}.
\end{equation}
\end{lem}

Both⋅$\op T_{i}^{\varepsilon}$ and⋅$\op S_{i}^{\varepsilon}$ converge
to the identity operator in the uniform operator topologies when restricted
to suitable Sobolev spaces. This is made precise in the next two lemmas.
\begin{lem}
\label{lem: C | Bound for (Te - I)}Let⋅$\Sigma$ be a domain in⋅$\R^{d}$,
and let⋅$q\in\sbr{1,\infty}$. Then, for any⋅$\varepsilon>0$ and⋅$u\in W_{q}^{1_{1}}\rbr{\Sigma_{\varepsilon_{i}\cellradius_{i}}\times\cell_{i}}$,
we⋅have
\begin{equation}
\norm{\rbr{\op T_{i}^{\varepsilon}-\op I}u}_{q,\Sigma\times\cell_{i}\times\cell_{i}}\lesssim\varepsilon_{i}\norm{D_{1}u}_{q,\Sigma_{\varepsilon_{i}\cellradius_{i}}\times\cell_{i}}.\label{est: C | Bound for (Te - I)}
\end{equation}
\end{lem}

\begin{lem}
\label{lem: C | Bound for (Se - I)}Let⋅$\Sigma$ be a domain in⋅$\R^{d}$,
and let⋅$q\in\sbr{1,\infty}$. Then, for any⋅$\varepsilon>0$ and⋅$u\in W_{q}^{1_{1}}\rbr{\Sigma_{\varepsilon_{i}\cellradius_{i}}\times\cell_{i}}$,
we⋅have
\begin{equation}
\norm{\rbr{\op S_{i}^{\varepsilon}-\op I}u}_{q,\Sigma\times\cell_{i}}\lesssim\varepsilon_{i}\norm{D_{1}u}_{q,\Sigma_{\varepsilon_{i}\cellradius_{i}}\times\cell_{i}}.\label{est: C | Bound for (Se - I)}
\end{equation}
\end{lem}

\subsection{Corrector}

For⋅each⋅$i\in\set I$, let⋅$\op K_{\mu,i}^{\varepsilon}\colon L_{p}\rbr{\Omega}^{n}\to W_{p}^{1}\rbr{\Omega_{i}}^{n}$
be given⋅by
\begin{equation}
\op K_{\mu,i}^{\varepsilon}=\op r_{i}\tau_{i}^{\varepsilon}\op S_{i}^{\varepsilon}\op J_{i}^{\varepsilon}\op K_{\mu,i}.\label{def: C | Corrector Ke_(i,mu)}
\end{equation}
Choose a  function~$\rho_{\varepsilon}\in C^{0,1}\rbr{\wb{\Omega}}$
with support in the two-sided (nonsymmetric) $3\varepsilon_{\#}\cellradius_{\#}$\nobreakdash-neighborhood
of⋅$\Gamma$ (see⋅Figure~\ref{fig: C | Rho_eps}) such that $\rho_{\varepsilon}$
equals~$1$ on the⋅$2\varepsilon_{\#}\cellradius_{\#}$-neigh\-borhood
of⋅$\Gamma$ and⋅$\norm{D\rho_{\varepsilon}}_{\infty,\Omega_{\#}}\lesssim\smash[b]{\varepsilon_{\#}^{-1}}$.
We set⋅$\chi_{\varepsilon}=1-\rho_{\varepsilon}$ and define the corrector~$\op K_{\mu}^{\varepsilon}\colon L_{p}\rbr{\Omega}^{n}\to W_{p}^{1}\rbr{\Omega}^{n}$⋅as
\begin{equation}
\op K_{\mu}^{\varepsilon}=\chi_{\varepsilon}\op K_{\mu,\#}^{\varepsilon}.\label{def: C | Corrector Ke}
\end{equation}
More explicitly, for⋅$f\in L_{p}\rbr{\Omega}^{n}$ and⋅$x\in\Omega_{i}$,
$i\in\set I$,
\begin{equation}
\op K_{\mu}^{\varepsilon}f\rbr x=\chi_{\varepsilon}\rbr x\int_{\cell_{i}}\int_{\cell_{i}}J_{i}\rbr w\hairspace N_{i}\rbr{x+\varepsilon_{i}z+\varepsilon_{i}w,\varepsilon_{i}^{-1}x}\hairspace\op e_{i}D\op r_{i}\rbr{\op A_{\mu}^{0}}^{-1}f\rbr{x+\varepsilon_{i}z+\varepsilon_{i}w}\dd z\dd w.\label{eq: C | Corrector Ke explicitly}
\end{equation}
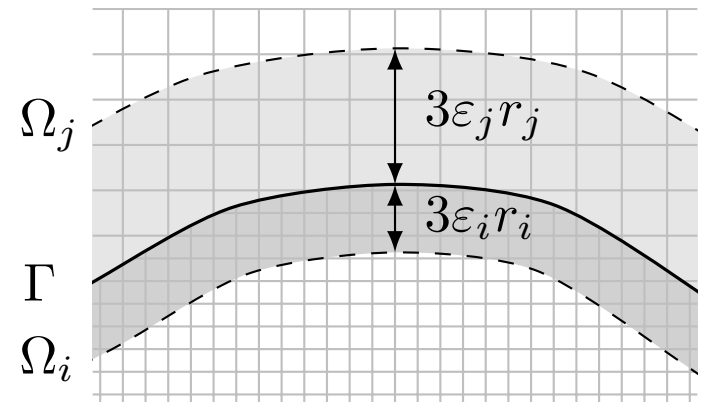
\begin{figure}
\centering
\begin{tikzpicture}[
  declare function={
    r=4; d=0.25; 
    eI=0.15; eII=0.30; 
    h=3; 
    m=d; 
    wx=2; ylo=r-h*eI-d-2*m; yhi=r+h*eII+d+m; 
  },
]
  \def\Seed{4096}
  \def\Samples{30,45,...,150} 
  \def\Window{(-wx,ylo) rectangle (wx,yhi)}

  \pgfmathsetseed{\Seed}
  \path[save path=\Interface, name path global=Iface]
    plot[smooth, variable=\t, samples at=\Samples] (\t:{r+d*rnd});
  \pgfmathsetseed{\Seed}
  \path[save path=\SuppI, name path global=SI]
    plot[smooth, variable=\t, samples at=\Samples] (\t:{r-h*eI+d*rnd});
  \pgfmathsetseed{\Seed}
  \path[save path=\SuppII, name path global=SII]
    plot[smooth, variable=\t, samples at=\Samples] (\t:{r+h*eII+d*rnd});

  \begin{scope}
    \clip \Window;
    \pgfmathsetseed{\Seed}
    \clip plot[smooth, variable=\t, samples at=\Samples] (\t:{r-h*eI+d*rnd})
          -- (150:{2*r}) -- (90:{2*r}) -- (30:{2*r}) -- cycle;
    \pgfmathsetseed{\Seed}
    \fill[black!18] plot[smooth, variable=\t, samples at=\Samples] (\t:{r+d*rnd})
          -- (0,0) -- cycle;
  \end{scope}
  \begin{scope}
    \clip \Window;
    \pgfmathsetseed{\Seed}
    \clip plot[smooth, variable=\t, samples at=\Samples] (\t:{r+d*rnd})
          -- (150:{2*r}) -- (90:{2*r}) -- (30:{2*r}) -- cycle;
    \pgfmathsetseed{\Seed}
    \fill[black!10] plot[smooth, variable=\t, samples at=\Samples] (\t:{r+h*eII+d*rnd})
          -- (0,0) -- cycle;
  \end{scope}

  \begin{scope}
    \clip \Window;
    \pgfmathsetseed{\Seed}
    \clip plot[smooth, variable=\t, samples at=\Samples] (\t:{r+d*rnd})
          -- (0,0) -- cycle;
    \draw[latticedimmed] (-wx,ylo) grid[step=eI] (wx,yhi);
  \end{scope}
  \begin{scope}
    \clip \Window;
    \pgfmathsetseed{\Seed}
    \clip plot[smooth, variable=\t, samples at=\Samples] (\t:{r+d*rnd})
          -- (150:{2*r}) -- (90:{2*r}) -- (30:{2*r}) -- cycle;
    \draw[latticedimmed] (-wx,ylo) grid[step=eII] (wx,yhi);
  \end{scope}

  \begin{scope}
    \clip \Window;
    \draw[supp, use path=\SuppI];
    \draw[supp, use path=\SuppII];
    \draw[interface, use path=\Interface];
  \end{scope}

  \path[name path=ray] (0,0) -- (90:{2*r});
  \path[name intersections={of=SI and ray, by=A}];
  \path[name intersections={of=Iface and ray, by=B}];
  \path[name intersections={of=SII and ray, by=C}];
  \draw[dimension] (A) -- (B) node[midway, right=2pt] {$3\varepsilon_i r_i$};
  \draw[dimension] (B) -- (C) node[midway, right=2pt] {$3\varepsilon_j r_j$};

  \path[name path=leftedge] (-wx,ylo) -- (-wx,yhi);
  \path[name intersections={of=SII and leftedge, by=J}];
  \path[name intersections={of=Iface and leftedge, by=G}];
  \path[name intersections={of=SI and leftedge, by=I}];
  \node[anchor=center] at ({-wx-0.35},0 |- G) {$\Gamma$};
  \node[anchor=center] (lblOmegaUp) at ({-wx-0.35},0 |- J) {$\Omega_{\mathrlap{j}}$};
  \node[anchor=center] (lblOmegaLow) at ({-wx-0.35},0 |- I) {$\Omega_{\mathrlap{i}}$};

  \pgfresetboundingbox
  \path[use as bounding box]
    (-wx,ylo) rectangle (wx,yhi);
\end{tikzpicture}\caption{The⋅$3\varepsilon_{\#}r_{\#}$\protect\nobreakdash-neighborhood of
a part of the interface~$\Gamma$, with different widths on the two
sides, containing the support of~$\rho_{\varepsilon}$}
\label{fig: C | Rho_eps}
\end{figure}

\begin{rem}
The⋅operators~$\op S_{i}^{\varepsilon}$ and⋅$\op J_{i}^{\varepsilon}$
serve two different purposes. The⋅first one separates the  different
scales, represented by the fast and slow variables, while the second
compensates for the limited regularity of the coefficients and the
effective operator resolvent. In⋅fact, $\op J_{i}^{\varepsilon}$
can be dropped from the definition of⋅$\op K_{\mu,i}^{\varepsilon}$
as long as⋅$\omega_{i}\rbr t=t$, \ie, the function~$\op r_{i}A$
is Lipschitz in the first variable and the resolvent of the effective
operator is continuous from⋅$L_{p}\rbr{\Omega}^{n}$ to⋅$W_{p}^{2}\rbr{\Omega_{i}}^{n}$.
If⋅this is the case for every⋅$i\in\set I$, then (\ref{eq: C | Corrector Ke explicitly})
becomes
\[
\op K_{\mu}^{\varepsilon}f\rbr x=\chi_{\varepsilon}\rbr x\int_{\cell_{\#}}N_{\#}\rbr{x+\varepsilon_{\#}z,\varepsilon_{\#}^{-1}x}\hairspace\op e_{\#}D\op r_{\#}\rbr{\op A_{\mu}^{0}}^{-1}f\rbr{x+\varepsilon_{\#}z}\dd z.
\]
For⋅locally periodic problems, the corrector~$\op K_{\mu}^{\varepsilon}$
in this form first appeared  in~\cite{PasT:2007}. We note, however,
that even in the most general case all⋅$\op J_{i}^{\varepsilon}$
(as⋅well as the cutoff function~$\chi_{\varepsilon}$) can be removed
from our results (but not from~$\op K_{\mu}^{\varepsilon}$) simply
by reformulating the statements, see Remark~\ref{rem: MR | Simpler corrector}.
\end{rem}
\begin{rem}
Note also that a mollification can be built from the Steklov smoothing
alone: $\rbr{\op S_{i}^{\varepsilon}}^{2}$ has the same form as⋅$\op J_{i}^{\varepsilon}$,
where $J_{i}$ is  the convolution of the characteristic function
of⋅$\cell_{i}$ with itself. As⋅observed, \eg, in⋅\cite{Pastukhova:2020-1},
this⋅$J_{i}$ is a nonnegative even Lipschitz function of integral~$1$.
It⋅vanishes outside⋅$2\cell_{i}$ rather than⋅$\cell_{i}$, but this
is not essential and only changes constants.
\end{rem}

Let us check that $\op K_{\mu}^{\varepsilon}$ is bounded. We write
\[
\norm{\varepsilon_{\#}D\op K_{\mu}^{\varepsilon}f}_{p,\Omega}\le\norm{\varepsilon_{\#}\chi_{\varepsilon}D\op K_{\mu,\#}^{\varepsilon}f}_{p,\Omega}+\norm{\varepsilon_{\#}\op K_{\mu,\#}^{\varepsilon}f\cdot D\chi_{\varepsilon}}_{p,\Omega}.
\]
Notice that⋅$\varepsilon_{\#}D\tau_{\#}^{\varepsilon}\op S_{\#}^{\varepsilon}=\varepsilon_{\#}\tau_{\#}^{\varepsilon}\op S_{\#}^{\varepsilon}D_{1}+\tau_{\#}^{\varepsilon}\op S_{\#}^{\varepsilon}D_{2}$.
It⋅then follows from Lemma~\ref{lem: C | Operator =0003C4e Se bounded on domains}
that
\[
\norm{\varepsilon_{\#}D\op K_{\mu}^{\varepsilon}f}_{p,\Omega}+\norm{\op K_{\mu}^{\varepsilon}f}_{p,\Omega}\lesssim\sum_{i\in\set I}\bigrbr{\varepsilon_{i}\norm{D_{1}\op J_{i}^{\varepsilon}\op K_{\mu,i}f}_{p,\R^{d}\times\cell_{i}}+\norm{\op J_{i}^{\varepsilon}\op K_{\mu,i}f}_{1_{2},p,\R^{d}\times\cell_{i}}},
\]
and⋅so, by⋅Lemmas~\ref{lem: C | Operators R0_(mu,i) and K_(mu,i) bounded}
and~\ref{lem: C | Operators Je_i R0_(mu,i) and Je_i K_(mu,i) bounded},
for any⋅$\varepsilon\in\set E_{\mu}$,
\begin{equation}
\norm{\varepsilon_{\#}D\op K_{\mu}^{\varepsilon}f}_{p,\Omega}+\norm{\op K_{\mu}^{\varepsilon}f}_{p,\Omega}\lesssim\norm f_{p,\Omega}.\label{est: C | Corrector Ke bounded}
\end{equation}
In⋅a⋅similar manner, applying the Sobolev embedding theorem in each
variable, we can show that, for every⋅$q\in\sbr{p,p^{*}}$,
\begin{equation}
\norm{\varepsilon_{\#}\omega_{\#}\rbr{\varepsilon_{\#}}^{-1}\op K_{\mu}^{\varepsilon}f}_{q,\Omega}\lesssim\norm f_{p,\Omega}.\label{est: C | Corrector Ke bounded in L_q}
\end{equation}

\section{\label{sec: MR}Main results}

We first fix some notation used the results below. For⋅each⋅$\varepsilon>0$,
let  $\varpi_{\join}\rbr{\varepsilon}$ be the largest of⋅$\varpi_{i}\rbr{\varepsilon}=\omega_{i}\rbr{\varepsilon_{i}}$.
By⋅the⋅problem parameters, we mean the numbers~$d$, $p$,  $n$,
$\mu$; the Lipschitz structures of the domain~$\Omega$ and of
the subdomains~$\Omega_{i}$; the periodicity cells~$\cell_{i}$;
the $C^{0,\omega_{i}}$\nobreakdash-norms of⋅$\smash[b]{A\restrict{\Omega_{i}}}$;
and the operator norm of~$\sigma$.

Our first result gives approximations under minimal regularity assumptions.
\begin{thm}
\label{thm: MR | Conv in L_q and approx in W_p^1}Suppose that⋅\textup{(\ref{assu: OO| Coefficient A regular})},
\textup{(\ref{assu: OO | Resolvent of Ae bounded})}, \textup{(\ref{assu: EO | Inverse of A(x) bounded})}
and⋅\textup{(\ref{assu: EO | Resolvent of A0 bounded})}–\textup{(\ref{assu: EO | Resolvent of A0 regular})}
hold. Then, for any⋅$\varepsilon\in\set E_{\mu}$ and⋅$f\in L_{p}\rbr{\Omega}^{n}$,
we⋅have
\begin{align}
\norm{\rbr{\op A_{\mu}^{\varepsilon}}^{-1}f-\rbr{\op A_{\mu}^{0}}^{-1}f}_{q,\Omega} & \lesssim\varpi_{\join}\rbr{\varepsilon}^{1/p}\norm f_{p,\Omega},\label{est: MR | Conv in L_q}\\
\norm{D\rbr{\op A_{\mu}^{\varepsilon}}^{-1}f-D\rbr{\op A_{\mu}^{0}}^{-1}f-D\varepsilon_{\#}\op K_{\mu}^{\varepsilon}f}_{p,\Omega} & \lesssim\varpi_{\join}\rbr{\varepsilon}^{1/p}\norm f_{p,\Omega},\label{est: MR | Approx in W_p^1}
\end{align}
where~$q\in\sbr{p,p^{*}}$. The⋅implicit constants  depend only
on the problem parameters, the exponent~$q$ and the constants in
the bounds~\textup{(\ref{assu: OO | Resolvent of Ae bounded})},
 \textup{(\ref{assu: EO | Inverse of A(x) bounded})} and~\textup{(\ref{assu: EO | Resolvent of A0 regular})}.
\end{thm}

If,⋅in addition, the adjoint effective problem satisfies the regularity
property~(\dualref{assu: EO | Resolvent of A0 regular}), then
the rate of convergence in⋅(\ref{est: MR | Conv in L_q}) for⋅$q=p$
can be  improved. Recall that, by our convention, $\smash{\omega_{\#}}\le\smash{\omega_{\#}^{\dual}}$,
and hence $\smash{\varpi_{\join}}\le\smash{\varpi_{\join}^{\dual}}$.
\begin{thm}
\label{thm: MR | Conv in L_p with better rate}Under⋅the⋅hypotheses
of Theorem~\textup{\ref{thm: MR | Conv in L_q and approx in W_p^1}},
suppose furthermore that \textup{(\dualref{assu: EO | Resolvent of A0 regular})}
holds. Then, for any⋅$\varepsilon\in\set E_{\mu}$ and⋅$f\in L_{p}\rbr{\Omega}^{n}$,
we⋅have
\begin{equation}
\norm{\rbr{\op A_{\mu}^{\varepsilon}}^{-1}f-\rbr{\op A_{\mu}^{0}}^{-1}f}_{p,\Omega}\lesssim\varpi_{\join}\rbr{\varepsilon}^{1/p}\varpi_{\join}^{\dual}\rbr{\varepsilon}^{1/p^{\dual}}\norm f_{p,\Omega}.\label{est: MR | Conv in L_p with better rate}
\end{equation}
The⋅implicit constant  depends only on the problem parameters and
the constants in the bounds~\textup{(\ref{assu: OO | Resolvent of Ae bounded})},
 \textup{(\ref{assu: EO | Inverse of A(x) bounded})}, \textup{(\ref{assu: EO | Resolvent of A0 regular})}
and  \textup{(\dualref{assu: EO | Resolvent of A0 regular})}.
\end{thm}

The⋅other estimate in Theorem~\ref{thm: MR | Conv in L_q and approx in W_p^1}
can also be improved, but only away from the boundary and the interface,
provided that the corresponding interior energy estimate is available.
\begin{thm}
\label{thm: MR | Interior approx in W_p^1}Under⋅the⋅hypotheses of
Theorem~\textup{\ref{thm: MR | Conv in L_q and approx in W_p^1}},
suppose furthermore that \textup{(\dualref{assu: EO | Resolvent of A0 regular})}
holds and that  the interior energy estimate~\textup{(\ref{assu: OO | Interior energy estimate})}
is valid for some⋅$\chi$ and⋅$\chi^{\prime}$ whose supports lie
at a positive distance from⋅$\partial\Omega\cup\Gamma$. Then, for
any⋅$\varepsilon\in\set E_{\mu}$ and⋅$f\in L_{p}\rbr{\Omega}^{n}$,
we⋅have
\begin{equation}
\norm{\chi\rbr{D\rbr{\op A_{\mu}^{\varepsilon}}^{-1}f-D\rbr{\op A_{\mu}^{0}}^{-1}f-D\varepsilon_{\#}\op K_{\mu}^{\varepsilon}f}}_{p,\Omega}\lesssim\varpi_{\join}\rbr{\varepsilon}^{1/p}\varpi_{\join}^{\dual}\rbr{\varepsilon}^{1/p^{\dual}}\norm f_{p,\Omega}.\label{est: MR | Interior approx in W_p^1}
\end{equation}
The⋅implicit constant  depends only on the problem parameters, the
functions~$\chi$ and⋅$\chi^{\prime}$ and the constants in the 
bounds~\textup{(\ref{assu: OO | Resolvent of Ae bounded})}, \textup{(\ref{assu: OO | Interior energy estimate})},
\textup{(\ref{assu: EO | Inverse of A(x) bounded})}, \textup{(\ref{assu: EO | Resolvent of A0 regular})}
and~\textup{(\dualref{assu: EO | Resolvent of A0 regular})}.
\end{thm}

\begin{rem}
\label{rem: MR | Constants}Although the bound~(\ref{assu: EO | Resolvent of A0 bounded})
is used in the proofs, the implicit constants in the above theorems
do not depend~it. This is possible because (\ref{assu: EO | Resolvent of A0 bounded})
holds with the constant from~(\ref{assu: OO | Resolvent of Ae bounded}).
Indeed, once (\ref{est: MR | Conv in L_q}) and⋅(\ref{est: MR | Approx in W_p^1})
are proved, we conclude that, for every⋅$f\in L_{p}\rbr{\Omega}^{n}$,
the function~$\smash{\rbr{\op A_{\mu}^{\varepsilon}}^{-1}}f$ converges
to⋅$\smash{\rbr{\op A_{\mu}^{0}}^{-1}}f$ weakly in⋅$\smash{W_{p}^{1}}\rbr{\Omega}^{n}$,
since the term~$\varepsilon_{\#}\smash{\op K_{\mu}^{\varepsilon}}f$
converges weakly to~$0$. (For⋅strongly elliptic operators, this
also follows from the classical homogenization results.) By⋅the⋅weak
lower semicontinuity of the norm and by the density of⋅$L_{p}\rbr{\Omega}^{n}$
in⋅$\set W_{p}^{-1}\rbr{\Omega;\C^{n}}$, (\ref{assu: OO | Resolvent of Ae bounded})
now implies⋅(\ref{assu: EO | Resolvent of A0 bounded}) with the same
constant.
\end{rem}
\begin{rem}
\label{rem: MR | Simpler corrector}The⋅estimates~(\ref{est: MR | Approx in W_p^1})
and⋅(\ref{est: MR | Interior approx in W_p^1}) provide approximation
for⋅$\rbr{\op A_{\mu}^{\varepsilon}}^{-1}$ in the Sobolev space seminorm,
which is the traditional way of expressing homogenization results.
We can get  a simpler approximation for the operator~$D\rbr{\op A_{\mu}^{\varepsilon}}^{-1}$
by dropping small terms that appear after differentiating~$\varepsilon_{\#}\op K_{\mu}^{\varepsilon}f$.
The⋅identity~$\smash{\varepsilon_{\#}D\tau_{\#}^{\varepsilon}\op S_{\#}^{\varepsilon}}=\smash{\varepsilon_{\#}\tau_{\#}^{\varepsilon}\op S_{\#}^{\varepsilon}}D_{1}+\smash{\tau_{\#}^{\varepsilon}\op S_{\#}^{\varepsilon}}D_{2}$
yields
\[
\begin{aligned}D\varepsilon_{\#}\op K_{\mu}^{\varepsilon} & =\tau_{\#}^{\varepsilon}\op S_{\#}^{\varepsilon}D_{2}\op K_{\mu,\#}+\tau_{\#}^{\varepsilon}\op S_{\#}^{\varepsilon}\rbr{\op J_{\#}^{\varepsilon}-\op I}D_{2}\op K_{\mu,\#}+\varepsilon_{\#}\chi_{\varepsilon}\tau_{\#}^{\varepsilon}\op S_{\#}^{\varepsilon}D_{1}\op J_{\#}^{\varepsilon}\op K_{\mu,\#}\\
 & \quad-\rho_{\varepsilon}\tau_{\#}^{\varepsilon}\op S_{\#}^{\varepsilon}D_{2}\op J_{\#}^{\varepsilon}\op K_{\mu,\#}-\varepsilon_{\#}\op K_{\mu,\#}^{\varepsilon}\cdot D\rho_{\varepsilon}.
\end{aligned}
\]
Arguing as in the proof of⋅(\ref{est: C | Corrector Ke bounded}),
one can show that
\[
\norm{\tau_{\#}^{\varepsilon}\op S_{\#}^{\varepsilon}\rbr{\op I-\op J_{\#}^{\varepsilon}}D_{2}\op K_{\mu,\#}f}_{p,\Omega}\lesssim\varpi_{\join}\rbr{\varepsilon}\hairspace\norm f_{p,\Omega}
\]
and
\[
\norm{\varepsilon_{\#}\chi_{\varepsilon}\tau_{\#}^{\varepsilon}\op S_{\#}^{\varepsilon}D_{1}\op J_{\#}^{\varepsilon}\op K_{\mu,\#}f}_{p,\Omega}\lesssim\varpi_{\join}\rbr{\varepsilon}\hairspace\norm f_{p,\Omega}.
\]
The⋅last two terms  are treated similarly, except that we first
apply Lemma~\ref{lem: Ps | Norm over a thin layer} below:
\[
\norm{\rho_{\varepsilon}\tau_{\#}^{\varepsilon}\op S_{\#}^{\varepsilon}D_{2}\op J_{\#}^{\varepsilon}\op K_{\mu,\#}f}_{p,\Omega}+\norm{\varepsilon_{\#}\op K_{\mu,\#}^{\varepsilon}f\cdot D\rho_{\varepsilon}}_{p,\Omega}\lesssim\varpi_{\join}\rbr{\varepsilon}^{1/p}\norm f_{p,\Omega}.
\]
Thus, $D\varepsilon_{\#}\op K_{\mu}^{\varepsilon}f$ in⋅(\ref{est: MR | Approx in W_p^1})
can be replaced by⋅$\tau_{\#}^{\varepsilon}\op S_{\#}^{\varepsilon}D_{2}\op K_{\mu,\#}f$.
The same is true for⋅(\ref{est: MR | Interior approx in W_p^1}),
because the last two terms, those with the slowest decay, vanish
on the support of⋅$\chi$ for small⋅$\varepsilon$, which is enough
(see the proof of Theorem~\ref{thm: MR | Interior approx in W_p^1}).
\end{rem}
\begin{rem}
Recall that the interior energy estimate~(\ref{assu: OO | Interior energy estimate}),
which is assumed in Theorem~\ref{thm: MR | Interior approx in W_p^1},
follows, for any compactly supported~$\chi$, from its analog with
a smaller exponent~$q>1$ provided that the resolvent of⋅$\op A^{\varepsilon}$
is uniformly bounded from⋅$\rbr{W_{q^{\hconj}}^{1}\rbr{\Omega}^{n}}^{*}$
to⋅$W_{q}^{1}\rbr{\Omega}^{n}$, see Lemma~\ref{lem: OO | Self-improvement of the interior energy estimate}.
\end{rem}
\begin{rem}
\label{rem: Boundary layer term}A⋅comparison of⋅(\ref{est: MR | Approx in W_p^1})
and⋅(\ref{est: MR | Interior approx in W_p^1})  suggests that the
estimate for⋅$D\rbr{\op A_{\mu}^{\varepsilon}}^{-1}$ is weaker
only near⋅$\partial\Omega$ and~$\Gamma$. Indeed, with a suitable
boundary-layer correction term~$\op B_{\mu}^{\varepsilon}$, one
obtains, for any⋅$\varepsilon\in\set E_{\mu}$ and⋅$f\in L_{p}\rbr{\Omega}^{n}$,
\begin{equation}
\norm{\rbr{\op A_{\mu}^{\varepsilon}}^{-1}f-\rbr{\op A_{\mu}^{0}}^{-1}f-\varepsilon_{\#}\op K_{\mu}^{\varepsilon}f-\op B_{\mu}^{\varepsilon}f}_{1,p,\Omega}\lesssim\varpi_{\join}\rbr{\varepsilon}\hairspace\norm f_{p,\Omega},\label{est: Error with boundary-layer correction term}
\end{equation}
see⋅(\ref{eq: Ps | Resolvent identity}) and Lemmas~\ref{lem: Ps | Operator Ie_(i,mu) bounded}
and~\ref{lem: Ps | Operator De_(i,mu) bounded}.  In⋅\cite{ZhPas:2005},
as well as in⋅\cite{PasT:2007}, \cite{PSu:2012}, \cite{Su:2013-1}
and~\cite{Su:2013-2}, an estimate of this kind was  a first step
 toward extending the analysis from the whole space to domains.
Constructing⋅$\smash[b]{\op B_{\mu}^{\varepsilon}}$ explicitly,
however, is no easier than solving the original problem, since it
essentially amounts to inverting~$\smash[b]{\op A_{\mu}^{\varepsilon}}$.
Instead, the boundary-layer approach
 requires only bounds on the norms of~$\smash[b]{\op B_{\mu}^{\varepsilon}}$.
 We may also read (\ref{est: Error with boundary-layer correction term})
in the opposite direction, from domains back to the whole space. If⋅there
is no⋅$\Gamma$ and $\Omega$ is a flat manifold without boundary,
such as the whole space or the torus, then⋅$\smash[b]{\op B_{\mu}^{\varepsilon}}$
vanishes, and  the rates in⋅(\ref{est: MR | Conv in L_q})–(\ref{est: MR | Approx in W_p^1})
improve to~$\varpi_{\join}\rbr{\varepsilon}$. At⋅least when⋅$\varpi_{\join}\rbr{\varepsilon}=\varepsilon$,
this rate is known to be sharp.
\end{rem}

\begin{rem}
\label{rem: MR | Regularity of the effective operator}Even⋅if⋅(\ref{assu: EO | Resolvent of A0 regular})
fails to hold, we may still recover some of our results in a weaker
form. Namely, inspection of the proof of Theorem~\ref{thm: MR | Conv in L_q and approx in W_p^1}
shows that its conclusion remains true for any fixed⋅$f\in L_{p}\rbr{\Omega}^{n}$,
provided that  $\rbr{\op A_{\mu}^{0}}^{-1}f\in\smash{\Lambda_{p}^{1+\omega_{i}}}\rbr{\Omega_{i}}^{n}$
for all⋅$i\in\set I$.  Of⋅course, in that case, the implicit constants
 may depend on the function~$f$, and its norm needs to be removed
from the right-hand sides of the estimates~(\ref{est: MR | Conv in L_q})–(\ref{est: MR | Approx in W_p^1}).
\end{rem}

\section{\label{sec: Examples}Example: Strongly elliptic operators}

In⋅this section, we assume that⋅$\Omega$ is a bounded domain divided
into a finite number of⋅$C^{1,\omega_{i}}$-regular subdomains~$\Omega_{i}$,
$i\in\set I$. Notice that the interface~$\Gamma_{ij}$ between
two adjacent subdomains~$\Omega_{i}$ and⋅$\Omega_{j}$, \ie, the
common part of their boundaries, is actually $C^{1,\omega_{ij}}$-regular
with~$\omega_{ij}=\omega_{i}\meet\omega_{j}$. We suppose that the
interfaces have no cross-points, either interior (more than two subdomains
meet) or boundary (an interface touches the boundary). Figure~\ref{subfig: OO | Smooth domains}
shows a typical geometry.

Other assumptions about the regularity of⋅$\partial\Omega_{i}$ are
also possible, including mere Lipschitz regularity or geometric conditions
such as convexity, but we will not discuss them  here. Nor will
we consider boundary conditions other than the Dirichlet or Robin
ones on the whole boundary~$\partial\Omega$, although mixed boundary
conditions can be studied similarly. (See⋅\cite[Section~7]{Senik:2023-1}
for examples of problems on Lipschitz domains and of mixed Dirichlet–Neumann
problems.)  Thus, from now on, $\smash{\set W_{p}^{1}}\rbr{\Omega;\C^{n}}$
stands for either⋅$\smash{\mathring{W}_{p}^{1}}\rbr{\Omega}^{n}$
or~$\smash{W_{p}^{1}}\rbr{\Omega}^{n}$.
\begin{rem}
\label{rem: Ex | Cross-points}Excluding cross-points rules out typical
sources of lower regularity of the effective solution. Near the points
where the interfaces meet each other or the boundary, its regularity
may drop, much as near the points where the boundary condition changes
in mixed problems or near corners of Lipschitz domains. At⋅interior
cross-points, the solution can be only slightly more regular than⋅$W_{p}^{1}$
(see the example in⋅\cite[Subsection~5.3]{Morin+Nochetto+Siebert:2000},
based on~\cite{Kellogg:1974}); at boundary cross-points, the boundary
condition matters as well, and the regularity available away from
such points extends to them only under compatibility conditions (see⋅\cite{Savare:1998}
and~\cite{Druet:2013}). In⋅the⋅absence of cross-points, piecewise
regularity in Sobolev and H\"older spaces up to the interfaces is
classical for smooth geometry (see,⋅\eg,~\cite[Theorem~4.20]{McL:2000}).
\end{rem}
Our argument is  perturbative  around⋅$\smash[b]{p}=2$, so we
begin with  that case. As⋅usual, $\smash[t]{\set H^{1}}\rbr{\Omega;\C^{n}}$
denotes the space~$\smash{\set W_{2}^{1}}\rbr{\Omega;\C^{n}}$.

Let⋅$\set P_{\sigma}$ be the set of all exponents~$p\in\rbr{1,\infty}$
for which
\begin{equation}
\Norm{\sigma\varrho}_{-1/p,p,\partial\Omega}\lesssim\norm{\varrho}_{1-1/p,p,\partial\Omega},\qquad\varrho\in\set W_{p}^{1-1/p}\rbr{\partial\Omega;\C^{n}}.\label{est: Ex | sigma bounded}
\end{equation}
We also introduce the set~$\smash[t]{\wh{\set P}_{\sigma}}$ of
all⋅$p\in\rbr{1,\infty}$ such that
\begin{equation}
\norm{\sigma\varrho}_{1-1/p,p,\partial\Omega}\lesssim\norm{\varrho}_{1-1/p,p,\partial\Omega},\qquad\varrho\in\set W_{p}^{1-1/p}\rbr{\partial\Omega;\C^{n}}.\label{assu: Ex | sigma better bounded}
\end{equation}
Both  are convex by interpolation, and clearly⋅$\smash[t]{\wh{\set P}_{\sigma}}\subset\set P_{\sigma}$.
The⋅estimate~(\ref{est: Ex | sigma bounded}) implies⋅(\ref{est: OO | Operator Ae bounded})
for all⋅$p\in\set P_{\sigma}$ and~$\varepsilon>0$. The⋅analogous
sets for the adjoint of⋅$\sigma$,⋅$\set P_{\sigma^{*}}$ and⋅$\smash[t]{\wh{\set P}_{\sigma^{*}}}$,
consist of those⋅$p^{\prime}$ for which⋅(\dualref{est: Ex | sigma bounded})
and⋅(\dualref{assu: Ex | sigma better bounded}) hold, respectively.
By⋅duality, $\set P_{\sigma^{*}}$ is the reflection of⋅$\set P_{\sigma}$
under~$p\mapsto p^{\prime}$, whereas $\smash[t]{\wh{\set P}_{\sigma^{*}}}$
may differ substantially from~$\smash[t]{\wh{\set P}_{\sigma}}$
(see  the remark below). In⋅the⋅Dirichlet case, all these sets are
taken to be~$\rbr{1,\infty}$.
\begin{rem}
The⋅following two examples illustrate the sets~$\set P_{\sigma}$,
$\smash[t]{\wh{\set P}_{\sigma}}$ and~$\smash[t]{\wh{\set P}_{\sigma^{*}}}$.
In⋅the⋅first, $\sigma$ is the operator of multiplication by a function~$\theta\in C^{0,s}\rbr{\partial\Omega}^{n\times n}$
with⋅$s\in\nobreak(0,1]$. Then $\set P_{\sigma}=\rbr{1,\infty}$,
and both⋅$\smash[t]{\wh{\set P}_{\sigma}}$ and⋅$\smash[t]{\wh{\set P}_{\sigma^{*}}}$
contain the interval⋅$\rbr{1,1/(1-\nobreak s)}$, since $C^{0,s}$\nobreakdash-functions
are multipliers on⋅$\smash{W_{p}^{1-1/p}}\rbr{\partial\Omega}$ whenever⋅$1-\nobreak1/p<\nobreak s$.
In⋅the⋅second example, $\sigma$ is a rank-one operator given by⋅$\sigma\varrho=\rbr{\varrho,\beta}_{\partial\Omega}\alpha$
with a nonzero⋅$\alpha\in C^{0,s}\rbr{\partial\Omega}^{n}$ and⋅$\beta=\chf_{E}\eta$,
where⋅$\eta\in\C^{n}$ is a nonzero vector and⋅$E\subset\partial\Omega$
is a set whose relative boundary is nonempty and Lipschitz. Then
$\set P_{\sigma}=\rbr{1,\infty}$ and $\smash[t]{\wh{\set P}_{\sigma}}$
contains⋅$\rbr{1,1/(1-\nobreak s)}$ as above, while  for the adjoint,
given by⋅$\sigma^{*}\varrho=\rbr{\varrho,\alpha}_{\partial\Omega}\beta$,
we have⋅$\smash[t]{\wh{\set P}_{\sigma^{*}}}=\rbr{1,2}$, because
$\chf_{E}\in\smash{W_{q}^{1-1/q}}\rbr{\partial\Omega}$ only for⋅$q<2$.
In⋅both examples, $\sigma$ is compact for every~$p\in\rbr{1,\infty}$.
\end{rem}
Let us turn to the estimate~(\ref{assu: OO | Resolvent of Ae bounded}).
To⋅ensure that it holds for⋅$p=2$, we impose the following uniform
coercivity condition: there are⋅$\varepsilon_{0}\in(0,1]$ with⋅$\varepsilon_{\join}\rbr{\varepsilon_{0}}\le1$,
$c_{A}>0$ and⋅$C_{A}<\infty$ such that, for all⋅$\varepsilon\in\set E=(0,\varepsilon_{0}]$,
\begin{equation}
\Re\rbr{A^{\varepsilon}Du,Du}_{\Omega}+C_{A}\norm u_{2,\Omega}^{2}\ge c_{A}\norm{Du}_{2,\Omega}^{2},\qquad u\in\set H^{1}\rbr{\Omega;\C^{n}}.\label{assu: Ex | Operator D^* Ae D coercive}
\end{equation}
It⋅can be shown that the operator~$\op A\rbr x$ is then strongly
coercive for⋅all⋅$x\in\Omega\setminus\Gamma$ (see~\cite[Lemma~7.2]{Senik:2023-1}),
and therefore the function~$A$ satisfies the Legendre–Ha\-da\-mard
condition
\begin{equation}
\Re\abr{A\rbr{\wc,\wc}\hairspace\eta\otimes\xi,\eta\otimes\xi}\ge c_{A}\abs{\eta}^{2}\abs{\xi}^{2},\qquad\eta\in\C^{n},\xi\in\R^{d}\label{est: Ex | Legendre-Hadamard}
\end{equation}
(\cf~(\ref{est: EO | Weaker ellipticity-type condition})). We also
assume that $\sigma$ is compact from⋅$\set H^{1/2}\rbr{\partial\Omega;\C^{n}}$
to⋅$\set H^{-1/2}\rbr{\partial\Omega;\C^{n}}$ (in⋅particular,⋅$2\in\set P_{\sigma}$).
Then, for any⋅$\delta>0$, there exists⋅$C_{\sigma}\rbr{\delta}<\infty$
such that
\begin{equation}
\Norm{\gamma^{*}\sigma\gamma u}_{-1,2,\Omega}\le\delta\norm{Du}_{2,\Omega}+C_{\sigma}\rbr{\delta}\hairspace\norm u_{2,\Omega},\qquad u\in\set H^{1}\rbr{\Omega;\C^{n}}.\label{est: Ex | sigma bounded for p=00003D2, with delta}
\end{equation}
Now (\ref{assu: Ex | Operator D^* Ae D coercive}) and⋅(\ref{est: Ex | sigma bounded for p=00003D2, with delta})
together with⋅(\ref{est: OO | Operator Ae bounded}) for⋅$p=2$ imply
that the operator~$\op A^{\varepsilon}$ is $m$-sectorial for all⋅$\varepsilon\in\set E$,
with a sector~$\set S$ depending on the constants in these estimates
but not on~$\varepsilon$. Consequently, $\op A_{\mu}^{\varepsilon}$
satisfies⋅(\ref{assu: OO | Resolvent of Ae bounded}) for⋅$p=2$ and
any⋅$\varepsilon\in\set E$ provided that~$\mu\notin\set S$.

If⋅the⋅exponent~2 is not an interior point of⋅$\set P_{\sigma}$,
we set⋅$\set P=\cbr 2$. Otherwise, we can go beyond the case~$p=2$.
For⋅simplicity of notation,    fix⋅$p_{\sigma}>2$ such that $\sbr{\smash{p_{\sigma}^{\hconj}},p_{\sigma}}\subset\set P_{\sigma}$.
 We  invoke Shneiberg's stability theorem  (see⋅\cite[Corollary~2]{Shneiberg:1974}
or~\cite[Section~17.2]{Agranovich:2013}). Let
\[
\kappa\rbr p=\tan\frac{\pi}{4}\frac{1/2-1/p}{1/2-1/p_{\mathrlap{\sigma}}}.
\]
Then, for any⋅$p\in\rbr{p_{\sigma}^{\hconj},p_{\sigma}}$ satisfying
\begin{equation}
\abs{\kappa\rbr p}<\kappa_{\mu}^{\varepsilon}=\min_{q\in\cbr{p_{\sigma}^{\hconj},p_{\sigma}}}\norm{\op A_{\mu}^{\varepsilon}}_{\B\rbr{\set W_{q}^{1},\set W_{q}^{-1}}}^{-1}\norm{\rbr{\op A_{\mu}^{\varepsilon}}^{-1}}_{\B\rbr{\set H^{-1},\set H^{1}}}^{-1},\label{est: Ex | Shneiberg's critical exponent}
\end{equation}
the⋅operator~$\smash{\op A_{\mu}^{\varepsilon}}$ viewed as acting
from⋅$\set W_{p}^{1}\rbr{\Omega;\C^{n}}$ to⋅$\set W_{p}^{-1}\rbr{\Omega;\C^{n}}$
has a bounded inverse and, moreover,
\begin{equation}
\norm{\rbr{\op A_{\mu}^{\varepsilon}}^{-1}}_{\B\rbr{\set W_{p}^{-1},\set W_{p}^{1}}}\le\rbr{1-\abs{\kappa\rbr p}/\kappa_{\mu}^{\varepsilon}}^{-1}\norm{\rbr{\op A_{\mu}^{\varepsilon}}^{-1}}_{\B\rbr{\set H^{-1},\set H^{1}}}.\label{est: Ex | Shneiberg's estimate for the norm of the inverse}
\end{equation}
We shrink the range of⋅$p$ by estimating⋅$\kappa_{\mu}^{\varepsilon}$
from below by a positive constant~$\kappa_{\mu}$ that depends only
on the constants in⋅(\ref{est: OO | Operator Ae bounded}) for⋅$p\in\cbr{\smash[t]{p_{\sigma}^{\hconj}},p_{\sigma}}$
and in⋅(\ref{assu: OO | Resolvent of Ae bounded}) for~$p=2$. Denote
by⋅$p_{\mu}$ the smaller of⋅$p_{\sigma}$ and the solution to⋅$\kappa\rbr p=\kappa_{\mu}$
and define~$\set P_{\mu}=\rbr{p_{\mu}^{\hconj},p_{\mu}}$. Now,
to conclude (\ref{assu: OO | Resolvent of Ae bounded}) from⋅(\ref{est: Ex | Shneiberg's estimate for the norm of the inverse})
for each⋅$p\in\set P_{\mu}$, it remains to bound the right-hand
side by replacing⋅$\kappa_{\mu}^{\varepsilon}$ with⋅$\kappa_{\mu}$
 and then use⋅(\ref{assu: OO | Resolvent of Ae bounded}) for~$p=2$.
  Finally, let⋅$\set P$ be the union of the intervals~$\set P_{\mu}$
with~$\mu\notin\set S$. The⋅resolvent identity and the Sobolev
embedding theorem then show that⋅(\ref{assu: OO | Resolvent of Ae bounded})
holds for any⋅$p\in\set P$, $\varepsilon\in\set E$ and~$\mu\notin\set S$.
\begin{rem}
The⋅coercivity condition~(\ref{assu: Ex | Operator D^* Ae D coercive})
depends on⋅$\varepsilon$ and has to be verified uniformly in~$\varepsilon\in\set E$.
Here we formulate simpler  conditions that do not involve⋅$\varepsilon$
and still imply  coercivity.

Fix~$m\ge n$. Let⋅$b\rbr D$ be a matrix first-order differential
operator with symbol
\[
\xi\mapsto b\rbr{\xi}=\sum_{k=1}^{d}b_{k}\xi_{k}
\]
where~$b_{k}\in\C^{m\times n}$. Suppose that $g$ is  a $\C^{m\times m}$-valued
function on⋅$\Omega\times\R^{d}$ such that $g\in\smash[t]{C^{0,\omega_{\#}}\rbr{\wb{\Omega}_{\#};\wt L_{\infty}\rbr{\cell_{\#}}}}$
and $\Re g$ is uniformly positive definite. Now define $A_{kl}=b_{k}^{*}gb_{l}$,
so that we have $D^{*}A^{\varepsilon}D=b\rbr D^{*}g^{\varepsilon}b\rbr D$.
Many operators of mathematical physics can actually be expressed in
this form (see, for instance, the discussion in~\cite{BSu:2003}).

For⋅the⋅Dirichlet problem, we assume that the symbol~$b\rbr{\xi}$
 has the maximal rank~$n$ for every⋅$\xi\in\R^{d}\setminus\cbr 0$,
or, equivalently, that there is⋅$c_{b}>0$ such that
\[
b\rbr{\xi}^{*}b\rbr{\xi}\ge c_{b}\abs{\xi}^{2},\qquad\xi\in\R^{d}.
\]
By⋅using the Fourier transform, one sees that
\begin{equation}
\norm{b\rbr D\hairspace u}_{2,\Omega}^{2}\ge c_{b}\norm{Du}_{2,\Omega}^{2},\qquad u\in\mathring{H}^{1}\rbr{\Omega}^{n},\label{est: Coercivity of b(D)* b(D) on Dirichlet H^1}
\end{equation}
and therefore
\[
\Re\rbr{A^{\varepsilon}Du,Du}_{\Omega}\ge c_{b}\norm{\rbr{\Re g}^{-1}}_{L_{\infty}}^{-1}\norm{Du}_{2,\Omega}^{2},\qquad u\in\mathring{H}^{1}\rbr{\Omega}^{n}.
\]

The⋅coercivity on⋅$H^{1}\rbr{\Omega}^{n}$ for the Robin or Neumann
problems requires a stronger condition that $b\rbr{\xi}$ be of maximal
rank~$n$ for any⋅$\xi\in\C^{d}\setminus\cbr 0$, not just⋅$\xi\in\R^{d}\setminus\cbr 0$.
This ensures that
\begin{equation}
\norm{b\rbr D\hairspace u}_{2,\Omega}^{2}\ge c_{b}\norm{Du}_{2,\Omega}^{2}-C_{b}\norm u_{2,\Omega}^{2},\qquad u\in H^{1}\rbr{\Omega}^{n},\label{est: Coercivity of b(D)* b(D) on H^1}
\end{equation}
with⋅$c_{b}>0$ and⋅$C_{b}<\infty$, see⋅\cite[Theorem~7.8 in Section~3.7]{Necas:2012}.
Hence
\[
\Re\rbr{A^{\varepsilon}Du,Du}_{\Omega}\ge\norm{\rbr{\Re g}^{-1}}_{L_{\infty}}^{-1}\bigrbr{c_{b}\norm{Du}_{2,\Omega}^{2}-C_{b}\norm u_{2,\Omega}^{2}},\qquad u\in H^{1}\rbr{\Omega}^{n}.
\]

Purely periodic homogenization problems for  operators of this
form were studied in a  series of papers by Birman and Suslina and
their students, beginning with~\cite{BSu:2003}. Among these, we
mention the works~\cite{PSu:2012} and⋅\cite{Su:2013-1} on the Dirichlet
problem and the work~\cite{Su:2013-2} on the Neumann problem.
\end{rem}
As⋅already noted, the uniform weak coercivity condition~(\ref{assu: Ex | Operator D^* Ae D coercive})
implies the strong coercivity of⋅$\op A\rbr x$ for all⋅$x\in\Omega\setminus\Gamma$.
More precisely, we have
\begin{equation}
\Re\rbr{\op A\rbr x\hairspace u,u}_{\cell_{\#}}\ge c_{A}\norm{Du}_{2,\cell_{\#}}^{2},\qquad u\in\wt H^{1}\rbr{\cell_{\#}}^{n},\label{est: Operator A(x) coercive}
\end{equation}
with⋅the⋅same⋅$c_{A}$ as in⋅(\ref{assu: Ex | Operator D^* Ae D coercive}),
see~\cite[Lemma~7.2]{Senik:2023-1}. It⋅follows that $\op A\rbr x$
satisfies⋅(\ref{assu: EO | Inverse of A(x) bounded}) for~$p=2$.
 The⋅boundedness and coercivity constants of⋅$\op A\rbr x$ are
no worse than those of⋅$\op A^{\varepsilon}$. Hence, by the above
argument, $\op A\rbr x-\mu$ is invertible for all⋅$p\in\set P$ and⋅$\mu\notin\set S$,
uniformly in⋅$x\in\Omega\setminus\Gamma$, and so is $\op A\rbr x$,
which gives~(\ref{assu: EO | Inverse of A(x) bounded}) for these~$p$.

We now turn to the effective operator. For⋅$p=2$ and⋅$\mu\notin\set S$,
the bound~(\ref{assu: OO | Resolvent of Ae bounded}) and classical
homogenization results show that the inverse of⋅$\smash{\op A_{\mu}^{\varepsilon}}$
converges in the weak operator topology to the inverse of⋅$\smash{\op A_{\mu}^{0}}$
(see,~\eg, \cite{BLP:1978},~\cite{Al:1992}); the limit is an
isomorphism by coercivity (see~\cite[Lemma~6.2]{Tartar:2010}). For⋅general⋅$p\in\set P$,
(\ref{assu: OO | Resolvent of Ae bounded}) then implies⋅(\ref{assu: EO | Resolvent of A0 bounded}),
as in Remark~\ref{rem: MR | Constants}.

 To⋅verify⋅(\ref{assu: EO | Resolvent of A0 regular}) as it stands,
we need $A$ to be more regular near the interface than  assumed
in~(\ref{assu: OO| Coefficient A regular}). Namely, we suppose that,
for each⋅$\Gamma_{ij}$, there is an open neighborhood~$G_{ij}$
such that $A$ has the modulus of continuity~$\omega_{ij}$ on both
sides of⋅$\Gamma_{ij}$ within⋅$G_{ij}$, that⋅is,
\begin{equation}
A\in C^{0,\omega_{ij}}\rbr{\overline{\Omega_{i}\cap G_{ij}};\wt L_{\infty}\rbr{\cell_{i}}}.\label{assu: Ex | Coefficient A additional regularity}
\end{equation}
This  assumption implies the same regularity for⋅$A^{0}$ (as⋅in
the proof of Lemma~\ref{lm: EO | Function N regular}) and allows
us to obtain⋅(\ref{assu: EO | Resolvent of A0 regular}) from⋅(\ref{assu: EO | Resolvent of A0 bounded})
for every⋅$p\in\set P\cap\smash[t]{\wh{\set P}_{\sigma}}$ via 
a difference quotient argument combined with suitable regularization.
The⋅proof is somewhat technical, and we defer it to the appendix,
see Theorem~\ref{thm: Appendix | Elliptic regularity}, which covers
the Dirichlet and Neumann problems. In⋅the⋅Robin case, we treat⋅$\varrho=\sigma\gamma u_{0}$,
where⋅$u_{0}=\rbr{\op A_{\mu}^{0}}^{-1}f$, as Neumann data. This
requires the effective Neumann operator to be invertible, which is
the case for all⋅$p\in\set P$ and~$\mu\notin\set S$. Indeed, replacing⋅$\sigma$
by~$0$ improves the bound~(\ref{est: OO | Operator Ae bounded})
and reduces the coercivity bound, obtained from⋅(\ref{assu: Ex | Operator D^* Ae D coercive})
and⋅(\ref{est: Ex | sigma bounded for p=00003D2, with delta}), to⋅(\ref{assu: Ex | Operator D^* Ae D coercive})
alone. These two bounds determine the sector, the estimate for the
inverse and, through~$\kappa_{\mu}$, the range of exponents. Hence,
the Neumann operator has a smaller sector and a larger range than
the Robin one.  If,⋅moreover, $p\in\smash[t]{\wh{\set P}_{\sigma}}$,
then $\varrho$ belongs to⋅$\smash{W_{p}^{1-1/p}}\rbr{\partial\Omega}^{n}$,
with norm controlled by⋅$\norm f_{p,\Omega}$, so Lemma~\ref{lem: Appendix | Regular Neumann data fit}
and Theorem~\ref{thm: Appendix | Elliptic regularity} apply.

\begin{rem}
\label{rem: Ex | On coefficient A additional regularity}The⋅assumption~(\ref{assu: Ex | Coefficient A additional regularity})
 is natural. Since the effective solution is continuous across each
interface~$\Gamma_{ij}$, its traces from both sides coincide, so
its regularity near⋅$\Gamma_{ij}$ is shared by⋅$\Omega_{i}$ and~$\Omega_{j}$.
If⋅the⋅coefficients on one side, say⋅$\Omega_{j}$, are rougher, they
can limit the regularity of the effective solution near⋅$\Gamma_{ij}$
on the other side as well, in which case it fails to belong to⋅$\smash{\Lambda_{p}^{1+\omega_{i}}}$
up to~$\Gamma_{ij}$. The⋅assumption~(\ref{assu: Ex | Coefficient A additional regularity})
excludes exactly this situation, and it concerns only the rougher
side and only a neighborhood of the interface.

However, (\ref{assu: Ex | Coefficient A additional regularity})
can be dropped at the cost of enlarging the moduli of continuity.
Let⋅$\smash[t]{\wt{\omega}_{i}}$ be the largest of⋅$\omega_{i}$
and the moduli~$\omega_{j}$ for the subdomains~$\Omega_{j}$ adjacent
to~$\Omega_{i}$; it is again a regular modulus of continuity. Since
$\omega_{i}\le\smash[t]{\wt{\omega}_{i}}\meet\smash[t]{\wt{\omega}_{j}}$
for each such~$j$, the assumption~(\ref{assu: Ex | Coefficient A additional regularity})
holds with⋅$\wt{\omega}_{i}$ in place of⋅$\omega_{i}$ for every⋅$A$
satisfying~(\ref{assu: OO| Coefficient A regular}), and hence (\ref{assu: EO | Resolvent of A0 regular})
holds with these moduli. Accordingly, $\varpi_{\join}\rbr{\varepsilon}$
in the main results becomes⋅$\smash[t]{\wt{\varpi}}_{\join}\rbr{\varepsilon}=\max_{i\in\set I}\smash[t]{\wt{\omega}_{i}}\rbr{\varepsilon_{i}}$.
If⋅the⋅scales~$\varepsilon_{i}$ coincide, then $\smash[t]{\wt{\varpi}_{\join}}\rbr{\varepsilon}=\varpi_{\join}\rbr{\varepsilon}$.
Otherwise, the additional assumption yields the rate~$\varpi_{\join}\rbr{\varepsilon}$,
which may be better.
\end{rem}

Of⋅course, the above results carry over to the dual counterparts,
since the adjoint of⋅$\smash{\op A_{\mu}^{\varepsilon}}$ has the
same boundedness and coercivity constants as~$\smash{\op A_{\mu}^{\varepsilon}}$;
the only change is that $\smash[t]{\wh{\set P}_{\sigma}}$ is replaced
by~$\smash[t]{\wh{\set P}_{\sigma^{*}}}$.

In⋅the⋅present setting, Theorems~\ref{thm: MR | Conv in L_q and approx in W_p^1}
and~\ref{thm: MR | Conv in L_p with better rate} specialize to
the following corollary.
\begin{cor}
\label{cor: Ex | Conv in L_p with better rate and approx in W_p^1}Suppose
that  $A$ satisfies⋅\textup{(\ref{assu: OO| Coefficient A regular})}
and⋅\textup{(\ref{assu: Ex | Coefficient A additional regularity})},
that the coercivity condition~\textup{(\ref{assu: Ex | Operator D^* Ae D coercive})}
holds and that $\sigma$ is compact from⋅$\set H^{1/2}\rbr{\partial\Omega;\C^{n}}$
to⋅$\set H^{-1/2}\rbr{\partial\Omega;\C^{n}}$. Then, for any⋅$p\in\set P\cap\smash[t]{\wh{\set P}_{\sigma}}$,
$\varepsilon\in\set E$, $\mu\notin\set S$ and⋅$f\in L_{p}\rbr{\Omega}^{n}$,
we⋅have
\begin{align}
\norm{\rbr{\op A_{\mu}^{\varepsilon}}^{-1}f-\rbr{\op A_{\mu}^{0}}^{-1}f}_{p,\Omega} & \lesssim\varpi_{\join}\rbr{\varepsilon}^{1/p}\norm f_{p,\Omega},\label{est: Ex | Conv in L_p with better rate}\\
\norm{D\rbr{\op A_{\mu}^{\varepsilon}}^{-1}f-D\rbr{\op A_{\mu}^{0}}^{-1}f-D\varepsilon_{\#}\op K_{\mu}^{\varepsilon}f}_{p,\Omega} & \lesssim\varpi_{\join}\rbr{\varepsilon}^{1/p}\norm f_{p,\Omega}.\label{est: Ex | Approx in W_p^1}
\end{align}
Moreover, if⋅$p^{\prime}\in\smash[t]{\wh{\set P}_{\sigma^{*}}}$,
then the rate in⋅\textup{(\ref{est: Ex | Conv in L_p with better rate})}
improves to~$\varpi_{\join}\rbr{\varepsilon}$. The⋅implicit constants
  depend only on the problem parameters, the $C^{1,\omega_{i}}$\nobreakdash-structures
of⋅$\Omega_{i}$,   the $C^{0,\omega_{ij}}$\nobreakdash-norms
of⋅$\smash[b]{A\restrict{\Omega_{i}\cap G_{ij}}}$,  the sets~$G_{ij}$
and the constants in the bounds~\textup{(\ref{est: Ex | sigma bounded})}–\textup{(\ref{est: Ex | sigma bounded for p=00003D2, with delta})}
and, for the improved rate,~\textup{(\dualref{assu: Ex | sigma better bounded})}.
\end{cor}
It⋅remains to discuss the interior energy estimate~(\ref{assu: OO | Interior energy estimate}).
 Let⋅$\Sigma$ be a domain compactly contained in⋅$\Omega_{i}$ for
some~$i\in\set I$. Choose a cutoff function~$\chi$  with⋅$\supp\chi\subset\nobreak\Omega_{i}$
that equals~$1$ near~$\Sigma$. Testing $\op A_{\mu}^{\varepsilon}u$
against⋅$\abs{\chi}^{2}u$ and using⋅(\ref{assu: Ex | Operator D^* Ae D coercive}),
we arrive at the classical Caccioppoli inequality:
\begin{equation}
\norm{\chi Du}_{2,\supp\chi}\lesssim\norm u_{2,\supp\chi}+\norm{\chi\op A_{\mu}^{\varepsilon}u}_{-1,2,\Omega}^{*},\qquad u\in\set H^{1}\rbr{\Omega;\C^{n}}.\label{est: Ex: Caccioppoli's}
\end{equation}
Then Lemma~\ref{lem: OO | Self-improvement of the interior energy estimate}
with⋅$q=2$ yields⋅(\ref{assu: OO | Interior energy estimate}) for
all exponents~$p\in\set P\cap[2,\infty)$.
\begin{cor}
\label{cor: Ex | Interior approx in W_p^1}Under the hypotheses of
Corollary~\textup{\ref{cor: Ex | Conv in L_p with better rate and approx in W_p^1}},
let $\Sigma$ be a domain compactly contained in some~$\Omega_{i}$.
Then, for any⋅$p\in\set P\cap\smash[t]{\wh{\set P}_{\sigma}}\cap[2,\infty)$
such that⋅$p^{\prime}\in\smash[t]{\wh{\set P}_{\sigma^{*}}}$, and
any⋅$\varepsilon\in\set E$, $\mu\notin\set S$ and⋅$f\in L_{p}\rbr{\Omega}^{n}$,
we⋅have
\begin{equation}
\norm{D\rbr{\op A_{\mu}^{\varepsilon}}^{-1}f-D\rbr{\op A_{\mu}^{0}}^{-1}f-D\varepsilon_{\#}\op K_{\mu}^{\varepsilon}f}_{p,\Sigma}\lesssim\varpi_{\join}\rbr{\varepsilon}\hairspace\norm f_{p,\Omega}.\label{est: Ex | Interior approx in W_p^1}
\end{equation}
The⋅implicit constant   depends only on the problem parameters,
the $C^{1,\omega_{i}}$\nobreakdash-structures of⋅$\Omega_{i}$,
  the $C^{0,\omega_{ij}}$\nobreakdash-norms of⋅$\smash[b]{A\restrict{\Omega_{i}\cap G_{ij}}}$,
 the sets~$G_{ij}$, the domain~$\Sigma$ and the constants in
the bounds⋅\textup{(\ref{est: Ex | sigma bounded})}–\textup{(\ref{est: Ex | sigma bounded for p=00003D2, with delta})}
and~\textup{(\dualref{assu: Ex | sigma better bounded})}.
\end{cor}
\begin{rem}
In⋅our construction, the interval~$\set P$ shrinks toward~$\cbr 2$
as the ellipticity weakens. Moreover, the examples in⋅\cite{Mey:1963}
show that (\ref{assu: OO | Resolvent of Ae bounded}) need not hold
for large exponents. Additional assumptions on the function~$A$
may, however, widen the range. For⋅instance, in the absence of interfaces,
 (\ref{assu: OO | Resolvent of Ae bounded}) is valid for all~$p\in\rbr{1,\infty}$
if⋅$A$ is $\VMO$ (\ie,⋅of vanishing mean oscillation) in the fast
variable uniformly in the slow one (see~\cite[Subsection~7.2]{Senik:2023-1}).
\end{rem}

\section{\label{sec: Proof of the main results}Proof of the  main results}

We start with two  identities, a resolvent identity and its localized
version.

\subsection{The resolvent identity}

We would like to obtain a  representation of the form
\begin{equation}
\rbr{\op A_{\mu}^{\varepsilon}}^{-1}-\rbr{\op A_{\mu}^{0}}^{-1}-\varepsilon_{\#}\op K_{\mu}^{\varepsilon}\restrict{L_{p}\rbr{\Omega}^{n}}=\sum_{i\in\set I}\bigrbr{\op I_{\mu,i}^{\varepsilon}+\op D_{\mu,i}^{\varepsilon}+\op B_{\mu,i}^{\varepsilon}}.\label{eq: Ps | Resolvent identity}
\end{equation}
Here $\op I_{\mu,i}^{\varepsilon}$ is an interior part, where homogenization
actually comes into play, $\op D_{\mu,i}^{\varepsilon}$ is another
interior part  involving mollification errors, and $\op B_{\mu,i}^{\varepsilon}$
is a boundary part dealing with effects arising from the boundary
of~$\Omega_{i}$.

For⋅$f\in\smash{L_{p}\rbr{\Omega}^{n}}$ and⋅$f^{\dual}\in\smash{\rbr{W_{p}^{1}\rbr{\Omega}^{n}}^{*}}$,
 let $u_{\varepsilon}=\smash{\rbr{\op A_{\mu}^{\varepsilon}}^{-1}}f$,
$u_{0,i}=\smash{\op R_{\mu,i}^{0}}f$, $U_{i}=\smash[b]{\op K_{\mu,i}}f$
and⋅$U_{\varepsilon,i}=\smash{\op K_{\mu,i}^{\varepsilon}}f=\smash{\tau_{i}^{\varepsilon}\op S_{i}^{\varepsilon}\op J_{i}^{\varepsilon}U_{i}}$.
As⋅before, the superscript~``$\dual$'' marks the dual counterparts,
\eg,⋅$u_{\varepsilon}^{\dual}=\smash{\rbr{\rbr{\op A_{\mu}^{\varepsilon}}^{\dual}}^{-1}}f^{\dual}$,
and the subscript~``$\#$'' in place of⋅$i$ refers to the corresponding
piecewise objects.

Applying the left-hand side of⋅(\ref{eq: Ps | Resolvent identity})
to⋅$f$ and then pairing the result with⋅$f^{\dual}$ gives 
\[
\rbr{u_{\varepsilon}-u_{0,\#}-\varepsilon_{\#}\chi_{\varepsilon}U_{\varepsilon,\#},f^{\dual}}_{\Omega}=\rbr{f,u_{\varepsilon}^{\dual}}_{\Omega}-\rbr{u_{0,\#},f^{\dual}}_{\Omega}-\rbr{\varepsilon_{\#}\chi_{\varepsilon}U_{\varepsilon,\#},f^{\dual}}_{\Omega}.
\]
Using the definition of⋅$u_{0,i}$ and⋅$u_{\varepsilon}^{\dual}$
again, for the first two terms on the right we have
\[
\rbr{f,u_{\varepsilon}^{\dual}}_{\Omega}-\rbr{u_{0,\#},f^{\dual}}_{\Omega}=\rbr{A^{0}Du_{0,\#},Du_{\varepsilon}^{\dual}}_{\Omega}-\rbr{A^{\varepsilon}Du_{0,\#},Du_{\varepsilon}^{\dual}}_{\Omega}.
\]
As⋅for⋅the⋅last term,  choose a function~$\xi_{\varepsilon}\in C^{0,1}\rbr{\smash[t]{\wb{\Omega}}}$
that  vanishes outside the $3\varepsilon_{\#}r_{\#}$\nobreakdash-neighborhood
of⋅$\partial\Omega_{\#}$ and equals~$1$ on the $2\varepsilon_{\#}r_{\#}$\nobreakdash-neighborhood
of⋅$\partial\Omega_{\#}$, with the natural bound⋅$\smash[b]{\norm{D\xi_{\varepsilon}}_{\infty,\Omega_{\#}}}\lesssim\smash{\varepsilon_{\#}^{-1}}$.
Taking $\xi_{\varepsilon}$ equal to⋅$\rho_{\varepsilon}$ outside
the⋅$3\varepsilon_{\#}r_{\#}$\nobreakdash-neighborhood of⋅$\partial\Omega$,
we see that $\rho_{\varepsilon}^{\prime}=\xi_{\varepsilon}-\rho_{\varepsilon}$
is supported in that neighborhood and vanishes in the $2\varepsilon_{\#}r_{\#}$\nobreakdash-neighborhood
of~$\Gamma$.  Set⋅$\eta_{\varepsilon}=1-\xi_{\varepsilon}$. Then
$\varepsilon_{\#}\eta_{\varepsilon}U_{\varepsilon,\#}$ belongs to⋅$\set W_{p}^{1}\rbr{\Omega;\C^{n}}$,
and we conclude that
\[
\rbr{\varepsilon_{\#}\eta_{\varepsilon}U_{\varepsilon,\#},f^{\dual}}_{\Omega}=\rbr{A^{\varepsilon}D\varepsilon_{\#}\eta_{\varepsilon}U_{\varepsilon,\#},Du_{\varepsilon}^{\dual}}_{\Omega}-\mu\rbr{\varepsilon_{\#}\eta_{\varepsilon}U_{\varepsilon,\#},u_{\varepsilon}^{\dual}}_{\Omega}.
\]
As⋅a⋅result,
\begin{equation}
\begin{aligned}\hspace{1em} & \hspace{-1em}\rbr{u_{\varepsilon}-u_{0,\#}-\varepsilon_{\#}\chi_{\varepsilon}U_{\varepsilon,\#},f^{\dual}}_{\Omega}\\
 & =\rbr{\eta_{\varepsilon}A^{0}Du_{0,\#},Du_{\varepsilon}^{\dual}}_{\Omega}-\rbr{\eta_{\varepsilon}A^{\varepsilon}D\rbr{u_{0,\#}+\varepsilon_{\#}U_{\varepsilon,\#}},Du_{\varepsilon}^{\dual}}_{\Omega}+\mu\rbr{\varepsilon_{\#}\eta_{\varepsilon}U_{\varepsilon,\#},u_{\varepsilon}^{\dual}}_{\Omega}\\
 & \quad+\rbr{\xi_{\varepsilon}\rbr{A^{0}-A^{\varepsilon}}Du_{0,\#},Du_{\varepsilon}^{\dual}}_{\Omega}+\rbr{A^{\varepsilon}\varepsilon_{\#}U_{\varepsilon,\#}D\xi_{\varepsilon},Du_{\varepsilon}^{\dual}}_{\Omega}-\rbr{\varepsilon_{\#}\rho_{\varepsilon}^{\prime}U_{\varepsilon,\#},f^{\dual}}_{\Omega}.
\end{aligned}
\label{eq: Ps | First step}
\end{equation}
Introducing the notation~$\eta_{\varepsilon,i}=\op r_{i}\eta_{\varepsilon}$,
$\xi_{\varepsilon,i}=\op r_{i}\xi_{\varepsilon}$, $\rho_{\varepsilon,i}^{\prime}=\op r_{i}\rho_{\varepsilon}^{\prime}$
and⋅$A_{i}^{0}=\op r_{i}A^{0}$, we rewrite  the right-hand side
as the sum over⋅$i\in\set I$⋅of
\[
\begin{aligned} & \rbr{\eta_{\varepsilon,i}A_{i}^{0}Du_{0,i},Du_{\varepsilon}^{\dual}}_{\Omega_{i}}-\rbr{\eta_{\varepsilon,i}A_{i}^{\varepsilon}D\rbr{u_{0,i}+\varepsilon_{i}U_{\varepsilon,i}},Du_{\varepsilon}^{\dual}}_{\Omega_{i}}+\varepsilon_{i}\mu\rbr{\eta_{\varepsilon,i}U_{\varepsilon,i},u_{\varepsilon}^{\dual}}_{\Omega_{i}}\\
 & \quad+\rbr{\xi_{\varepsilon,i}\rbr{A_{i}^{0}-A_{i}^{\varepsilon}}Du_{0,i},Du_{\varepsilon}^{\dual}}_{\Omega_{i}}+\varepsilon_{i}\rbr{A_{i}^{\varepsilon}U_{\varepsilon,i}D\xi_{\varepsilon,i},Du_{\varepsilon}^{\dual}}_{\Omega_{i}}-\varepsilon_{i}\rbr{\rho_{\varepsilon,i}^{\prime}U_{\varepsilon,i},f^{\dual}}_{\Omega_{i}}.
\end{aligned}
\]
The⋅argument in⋅\cite[Subsection~8.1]{Senik:2023-1} then carries
over to each of these, yielding the representation~$\op I_{\mu,i}^{\varepsilon}+\op D_{\mu,i}^{\varepsilon}+\op B_{\mu,i}^{\varepsilon}$,
where
\begin{equation}
\begin{aligned}\rbr{\op I_{\mu,i}^{\varepsilon}f,f^{\dual}}_{\Omega} & =\rbr{\tau_{i}^{\varepsilon}\op T_{i}^{\varepsilon}\eta_{\varepsilon,i}D_{1}^{*}\op J_{i}^{\varepsilon}A_{i}\rbr{\op e_{i}D_{1}u_{0,i}+D_{2}U_{i}},\rbr{\op T_{i}^{\varepsilon}-\op I}u_{\varepsilon}^{\dual}}_{\Omega_{i}\times\cell_{i}}\\
 & \quad-\rbr{\tau_{i}^{\varepsilon}\op T_{i}^{\varepsilon}\eta_{\varepsilon,i}\sbr{A_{i},\op J_{i}^{\varepsilon}}\rbr{\op e_{i}D_{1}u_{0,i}+D_{2}U_{i}},D_{1}u_{\varepsilon}^{\dual}}_{\Omega_{i}\times\cell_{i}}\\
 & \quad-\rbr{\tau_{i}^{\varepsilon}\eta_{\varepsilon,i}\sbr{A_{i},\op T_{i}^{\varepsilon}}\op J_{i}^{\varepsilon}\rbr{\op e_{i}D_{1}u_{0,i}+D_{2}U_{i}},D_{1}u_{\varepsilon}^{\dual}}_{\Omega_{i}\times\cell_{i}}\\
 & \quad-\varepsilon_{i}\rbr{\tau_{i}^{\varepsilon}\eta_{\varepsilon,i}A_{i}\op T_{i}^{\varepsilon}D_{1}\op J_{i}^{\varepsilon}U_{i},D_{1}u_{\varepsilon}^{\dual}}_{\Omega_{i}\times\cell_{i}}\\
 & \quad-\rbr{\eta_{\varepsilon,i}A_{i}^{\varepsilon}\rbr{\op I-\op S_{i}^{\varepsilon}}\op J_{i}^{\varepsilon}\op e_{i}Du_{0,i},Du_{\varepsilon}^{\dual}}_{\Omega_{i}}\\
 & \quad+\varepsilon_{i}\mu\rbr{\eta_{\varepsilon,i}U_{\varepsilon,i},u_{\varepsilon}^{\dual}}_{\Omega_{i}},
\end{aligned}
\label{def: Ps | Operator I_(mu,i)}
\end{equation}
\begin{equation}
\begin{aligned}\rbr{\op D_{\mu,i}^{\varepsilon}f,f^{\dual}}_{\Omega} & =\rbr{\eta_{\varepsilon,i}\rbr{\op I-\op J_{i}^{\varepsilon}}A_{i}\rbr{\op e_{i}D_{1}u_{0,i}+D_{2}U_{i}},D_{1}u_{\varepsilon}^{\dual}}_{\Omega_{i}\times\cell_{i}}\\
 & \quad-\rbr{\eta_{\varepsilon,i}A_{i}^{\varepsilon}\rbr{\op I-\op J_{i}^{\varepsilon}}\op e_{i}Du_{0,i},Du_{\varepsilon}^{\dual}}_{\Omega_{i}}
\end{aligned}
\label{def: Ps | Operator D_(mu,i)}
\end{equation}
and
\begin{equation}
\begin{aligned}\rbr{\op B_{\mu,i}^{\varepsilon}f,f^{\dual}}_{\Omega} & =\rbr{\tau_{i}^{\varepsilon}\op T_{i}^{\varepsilon}\op J_{i}^{\varepsilon}A_{i}\rbr{\op e_{i}D_{1}u_{0,i}+D_{2}U_{i}}\cdot\op T_{i}^{\varepsilon}\rbr{D_{1}\xi_{\varepsilon,i}}^{*},\rbr{\op T_{i}^{\varepsilon}-\op I}u_{\varepsilon}^{\dual}}_{\Omega_{i}\times\cell_{i}}\\
 & \quad+\rbr{\tau_{i}^{\varepsilon}\sbr{\xi_{\varepsilon,i},\op T_{i}^{\varepsilon}}A_{i}\op J_{i}^{\varepsilon}\rbr{\op e_{i}D_{1}u_{0,i}+D_{2}U_{i}},D_{1}u_{\varepsilon}^{\dual}}_{\Omega_{i}\times\cell_{i}}\\
 & \quad+\rbr{\xi_{\varepsilon,i}\rbr{A_{i}^{0}-A_{i}^{\varepsilon}}Du_{0,i},Du_{\varepsilon}^{\dual}}_{\Omega_{i}}\\
 & \quad+\varepsilon_{i}\rbr{A_{i}^{\varepsilon}U_{\varepsilon,i}D\xi_{\varepsilon,i},Du_{\varepsilon}^{\dual}}_{\Omega_{i}}\\
 & \quad-\varepsilon_{i}\rbr{\rho_{\varepsilon,i}^{\prime}U_{\varepsilon,i},f^{\dual}}_{\Omega_{i}}
\end{aligned}
\label{def: Ps | Operator B_(mu,i)}
\end{equation}
(\cf~(8.1) and⋅(8.11)–(8.14) in~\cite{Senik:2023-1}). The⋅difference
from⋅\cite{Senik:2023-1} is that  here we extend the derivative
of the effective solution rather than the effective solution itself.

\subsection{A localized resolvent identity}

We now localize the above resolvent identity so that only interior
terms remain.

Take an arbitrary function~$\eta\in C^{0,1}\rbr{\smash{\wb{\Omega}}}$
that vanishes in the⋅$\delta$-neighborhood of⋅$\partial\Omega\cup\Gamma$
for some~$\delta>0$. Let⋅$\smash[t]{\mathring{\op A}^{\varepsilon}}\colon W_{p}^{1}\rbr{\Omega}^{n}\to\smash{\rbr{W_{p^{\prime}}^{1}\rbr{\Omega}^{n}}^{*}}$
be the operator associated with the form~$\rbr{A^{\varepsilon}Du,Du^{\prime}}_{\Omega}$,
and set⋅$\smash{\mathring{\op A}_{\mu}^{\varepsilon}}=\smash[t]{\mathring{\op A}^{\varepsilon}}-\mu$
(here and below, the ring marks objects related to the localized identity).
Observe that
\[
\rbr{\wb{\eta}\mathring{\op A}_{\mu}^{\varepsilon}u_{\varepsilon},u_{\varepsilon}^{\dual}}_{\Omega}=\rbr{f,\eta u_{\varepsilon}^{\dual}}_{\Omega}=\rbr{\op A_{\mu}^{0}u_{0,\#},\eta u_{\varepsilon}^{\dual}}_{\Omega}.
\]
Then we have
\[
\begin{aligned}\rbr{\wb{\eta}\mathring{\op A}_{\mu}^{\varepsilon}\rbr{u_{\varepsilon}-u_{0,\#}-\varepsilon_{\#}\chi_{\varepsilon}U_{\varepsilon,\#}},u_{\varepsilon}^{\dual}}_{\Omega} & =\rbr{A^{0}Du_{0,\#},D\eta u_{\varepsilon}^{\dual}}_{\Omega}\\
 & \quad-\rbr{A^{\varepsilon}D\rbr{u_{0,\#}+\varepsilon_{\#}\chi_{\varepsilon}U_{\varepsilon,\#}},D\eta u_{\varepsilon}^{\dual}}_{\Omega}\\
 & \quad+\mu\rbr{\varepsilon_{\#}\chi_{\varepsilon}U_{\varepsilon,\#},\eta u_{\varepsilon}^{\dual}}_{\Omega}.
\end{aligned}
\]
Note that the functional~$\smash[t]{\wb{\eta}}\smash{\mathring{\op A}_{\mu}^{\varepsilon}}\rbr{u_{\varepsilon}-u_{0,\#}-\varepsilon_{\#}\chi_{\varepsilon}U_{\varepsilon,\#}}$
belongs to⋅$\smash{\rbr{W_{p^{\hconj}}^{1}\rbr{\Omega}^{n}}^{*}}$
and hence, by the canonical projection~(\ref{def: OO | Canonical projection for the dual}),
to~$\smash{\set W_{p}^{-1}\rbr{\Omega;\C^{n}}}$. Next, the first
two terms on the right are the same as those in⋅(\ref{eq: Ps | First step}),
except that $\eta_{\varepsilon}$ is removed and $u_{\varepsilon}^{\dual}$
is replaced by⋅$\eta\smash{u_{\varepsilon}^{\dual}}$. Suppose that
$\varepsilon$ is small enough that⋅$\max_{i\in\set I}5\varepsilon_{i}r_{i}\le\delta$.
Then $\eta_{\varepsilon}=1$ on⋅$\supp\eta$, and, repeating the
argument  that led to⋅(\ref{eq: Ps | Resolvent identity}) (\cf~(8.21)
in~\cite{Senik:2023-1}), we obtain, for all such~$\varepsilon$,
that
\begin{equation}
\rbr{\op A_{\mu}^{\varepsilon}}^{-1}\wb{\eta}\mathring{\op A}_{\mu}^{\varepsilon}\rbr{\rbr{\op A_{\mu}^{\varepsilon}}^{-1}-\rbr{\op A_{\mu}^{0}}^{-1}-\varepsilon_{\#}\op K_{\mu}^{\varepsilon}}\restrict{L_{p}\rbr{\Omega}^{n}}=\sum_{i\in\set I}\bigrbr{\mathring{\op I}_{\mu,i}^{\varepsilon}+\mathring{\op D}_{\mu,i}^{\varepsilon}},\label{eq: Localized resolvent identity}
\end{equation}
where
\begin{equation}
\begin{aligned}\rbr{\mathring{\op I}_{\mu,i}^{\varepsilon}f,f^{\dual}}_{\Omega} & =\rbr{\tau_{i}^{\varepsilon}\op T_{i}^{\varepsilon}D_{1}^{*}\op J_{i}^{\varepsilon}A_{i}\rbr{\op e_{i}D_{1}u_{0,i}+D_{2}U_{i}},\rbr{\op T_{i}^{\varepsilon}-\op I}\eta u_{\varepsilon}^{\dual}}_{\Omega_{i}\times\cell_{i}}\\
 & \quad-\rbr{\tau_{i}^{\varepsilon}\op T_{i}^{\varepsilon}\sbr{A_{i},\op J_{i}^{\varepsilon}}\rbr{\op e_{i}D_{1}u_{0,i}+D_{2}U_{i}},D_{1}\eta u_{\varepsilon}^{\dual}}_{\Omega_{i}\times\cell_{i}}\\
 & \quad-\rbr{\tau_{i}^{\varepsilon}\sbr{A_{i},\op T_{i}^{\varepsilon}}\op J_{i}^{\varepsilon}\rbr{\op e_{i}D_{1}u_{0,i}+D_{2}U_{i}},D_{1}\eta u_{\varepsilon}^{\dual}}_{\Omega_{i}\times\cell_{i}}\\
 & \quad-\varepsilon_{i}\rbr{\tau_{i}^{\varepsilon}A_{i}\op T_{i}^{\varepsilon}D_{1}\op J_{i}^{\varepsilon}U_{i},D_{1}\eta u_{\varepsilon}^{\dual}}_{\Omega_{i}\times\cell_{i}}\\
 & \quad-\rbr{A_{i}^{\varepsilon}\rbr{\op I-\op S_{i}^{\varepsilon}}\op J_{i}^{\varepsilon}\op e_{i}Du_{0,i},D\eta u_{\varepsilon}^{\dual}}_{\Omega_{i}}\\
 & \quad+\varepsilon_{i}\mu\rbr{U_{\varepsilon,i},\eta u_{\varepsilon}^{\dual}}_{\Omega_{i}}
\end{aligned}
\label{def: Ps | Operator I0_(mu,i)}
\end{equation}
and
\begin{equation}
\begin{aligned}\rbr{\mathring{\op D}_{\mu,i}^{\varepsilon}f,f^{\dual}}_{\Omega} & =\rbr{\rbr{\op I-\op J_{i}^{\varepsilon}}A_{i}\rbr{\op e_{i}D_{1}u_{0,i}+D_{2}U_{i}},D_{1}\eta u_{\varepsilon}^{\dual}}_{\Omega_{i}\times\cell_{i}}\\
 & \quad-\rbr{A_{i}^{\varepsilon}\rbr{\op I-\op J_{i}^{\varepsilon}}\op e_{i}Du_{0,i},D\eta u_{\varepsilon}^{\dual}}_{\Omega_{i}}.
\end{aligned}
\label{def: Ps | Operator D0_(mu,i)}
\end{equation}

\subsection{Proofs of the main results}

The⋅proofs below closely follow⋅\cite{Senik:2023-1}, but we give
the details for the reader's convenience. Bounds for the terms~$\smash{\op I_{\mu,i}^{\varepsilon}}$,
$\smash{\op D_{\mu,i}^{\varepsilon}}$ and⋅$\smash{\op B_{\mu,i}^{\varepsilon}}$
are derived in a sequence of lemmas. Throughout, we keep in mind that
$\eta_{\varepsilon,i}$ vanishes in⋅$\smash[b]{\rbr{\partial\Omega_{i}}_{2\varepsilon_{i}\cellradius_{i}}}\cap\smash{\wb{\Omega}_{i}}$
and⋅$\xi_{\varepsilon,i}$ is supported in~$\smash[b]{\rbr{\partial\Omega_{i}}_{3\varepsilon_{i}\cellradius_{i}}}\cap\smash{\wb{\Omega}_{i}}$.
For⋅brevity, we write⋅$H_{\varepsilon,i}$ for⋅$\supp\eta_{\varepsilon,i}$,
and $X_{\varepsilon,i}$ and⋅$\smash{X_{\varepsilon,i}^{\prime}}$
for⋅$\supp\xi_{\varepsilon,i}$ and⋅$\supp D\xi_{\varepsilon,i}$,
respectively.  Recall also the notation for anisotropic norms from
Remark~\ref{rem: Cor | Anisotropic norms}.
\begin{lem}
\label{lem: Ps | Operator Ie_(i,mu) bounded}For⋅any⋅$\varepsilon\in\set E_{\mu}$
and⋅$f\in L_{p}\rbr{\Omega}^{n}$,
\begin{equation}
\norm{\op I_{\mu,i}^{\varepsilon}f}_{1,p,\Omega}\lesssim\omega_{i}\rbr{\varepsilon_{i}}\hairspace\norm f_{p,\Omega},\label{est: Ps | Operator Ie_(i,mu) bounded}
\end{equation}
where⋅$\op I_{\mu,i}^{\varepsilon}$ is given by~\textup{(\ref{def: Ps | Operator I_(mu,i)})}.
\end{lem}
\begin{proof}
We estimate each term in⋅(\ref{def: Ps | Operator I_(mu,i)}) using
Lemmas~\ref{lem: C | Operator =0003C4e Te bounded on domains},⋅\ref{lem: C | Operator =0003C4e Se bounded on domains}
and⋅\ref{lem: C | Bound for (Te - I)},⋅\ref{lem: C | Bound for (Se - I)}.
For⋅the⋅first, we have
\[
\begin{aligned}\hspace{1em} & \hspace{-1em}\bigabs{\rbr{\tau_{i}^{\varepsilon}\op T_{i}^{\varepsilon}\eta_{\varepsilon,i}D_{1}^{*}\op J_{i}^{\varepsilon}A_{i}\rbr{\op e_{i}D_{1}u_{0,i}+D_{2}U_{i}},\rbr{\op T_{i}^{\varepsilon}-\op I}u_{\varepsilon}^{\dual}}_{\Omega_{i}\times\cell_{i}}}\\
 & \le\norm{\tau_{i}^{\varepsilon}\op T_{i}^{\varepsilon}\eta_{\varepsilon,i}D_{1}^{*}\op J_{i}^{\varepsilon}A_{i}\rbr{\op e_{i}D_{1}u_{0,i}\!+\!D_{2}U_{i}}}_{p,\rbr{H_{\varepsilon,i}}_{\varepsilon_{i}\cellradius_{i}}\!\times\cell_{i}}\norm{\rbr{\op T_{i}^{\varepsilon}\!-\!\op I}u_{\varepsilon}^{\dual}}_{p^{\dual},\rbr{H_{\varepsilon,i}}_{\varepsilon_{i}\cellradius_{i}}\!\times\cell_{i}}\\
 & \lesssim\varepsilon_{i}\norm{D_{1}\op J_{i}^{\varepsilon}A_{i}\rbr{\op e_{i}D_{1}u_{0,i}+D_{2}U_{i}}}_{p,\Omega_{i}\times\cell_{i}}\norm{Du_{\varepsilon}^{\dual}}_{p^{\dual},\Omega_{i}}.
\end{aligned}
\]
Next, expressing⋅$\op J_{i}^{\varepsilon}$ in the second term 
as in Remark~\ref{rem: C | Mollification Je_i other expression}
and switching the order of integration yield
\[
\begin{aligned}\hspace{2em} & \hspace{-2em}\rbr{\tau_{i}^{\varepsilon}\op T_{i}^{\varepsilon}\eta_{\varepsilon,i}\sbr{A_{i},\op J_{i}^{\varepsilon}}\rbr{\op e_{i}D_{1}u_{0,i}+D_{2}U_{i}},D_{1}u_{\varepsilon}^{\dual}}_{\Omega_{i}\times\cell_{i}}\\
 & =\rbr{\tau_{i}^{\varepsilon}\op S_{i}^{\varepsilon}\eta_{\varepsilon,i}J_{i}\sbr{A_{i},\op T_{i}^{\varepsilon}}\rbr{\op e_{i}D_{1}u_{0,i}+D_{2}U_{i}},D_{1}u_{\varepsilon}^{\dual}}_{\Omega_{i}\times\cell_{i}}.
\end{aligned}
\]
Now notice that⋅$\sbr{A_{i},\op T_{i}^{\varepsilon}}=\rbr{\op I-\op T_{i}^{\varepsilon}}A_{i}\cdot\op T_{i}^{\varepsilon}$.
Then
\[
\begin{aligned}\hspace{1em} & \hspace{-1em}\bigabs{\rbr{\tau_{i}^{\varepsilon}\op T_{i}^{\varepsilon}\eta_{\varepsilon,i}\sbr{A_{i},\op J_{i}^{\varepsilon}}\rbr{\op e_{i}D_{1}u_{0,i}+D_{2}U_{i}},D_{1}u_{\varepsilon}^{\dual}}_{\Omega_{i}\times\cell_{i}}}\\
 & \le\norm{J_{i}}_{\infty,\cell_{i}}\norm{\rbr{\op I-\op T_{i}^{\varepsilon}}A_{i}}_{L_{\infty}}\norm{\tau_{i}^{\varepsilon}\op T_{i}^{\varepsilon}\rbr{\op e_{i}D_{1}u_{0,i}+D_{2}U_{i}}}_{p,H_{\varepsilon,i}\times\cell_{i}}\norm{D_{1}u_{\varepsilon}^{\dual}}_{p^{\dual},\Omega_{i}\times\cell_{i}}\\
 & \lesssim\omega_{i}\rbr{\varepsilon_{i}}\hairspace\bigrbr{\norm{Du_{0,i}}_{p,\Omega_{i}}+\norm{D_{2}U_{i}}_{p,\Omega_{i}\times\cell_{i}}}\norm{Du_{\varepsilon}^{\dual}}_{p^{\dual},\Omega_{i}}
\end{aligned}
\]
according to~(\ref{eq: C | Functions Ai and Ni regular}). Similarly,
by the identity~$\tau_{i}^{\varepsilon}\sbr{A_{i},\op T_{i}^{\varepsilon}}=\tau_{i}^{\varepsilon}\rbr{\op I-\op T_{i}^{\varepsilon}}A_{i}\cdot\tau_{i}^{\varepsilon}\op T_{i}^{\varepsilon}$,
\[
\begin{aligned}\hspace{2em} & \hspace{-2em}\bigabs{\rbr{\tau_{i}^{\varepsilon}\eta_{\varepsilon,i}\sbr{A_{i},\op T_{i}^{\varepsilon}}\op J_{i}^{\varepsilon}\rbr{\op e_{i}D_{1}u_{0,i}+D_{2}U_{i}},D_{1}u_{\varepsilon}^{\dual}}_{\Omega_{i}\times\cell_{i}}}\\
 & \le\norm{\rbr{\op I-\op T_{i}^{\varepsilon}}A_{i}}_{L_{\infty}}\norm{\tau_{i}^{\varepsilon}\op T_{i}^{\varepsilon}\op J_{i}^{\varepsilon}\rbr{\op e_{i}D_{1}u_{0,i}+D_{2}U_{i}}}_{p,H_{\varepsilon,i}\times\cell_{i}}\norm{D_{1}u_{\varepsilon}^{\dual}}_{p^{\dual},\Omega_{i}\times\cell_{i}}\\
 & \lesssim\omega_{i}\rbr{\varepsilon_{i}}\hairspace\bigrbr{\norm{Du_{0,i}}_{p,\Omega_{i}}+\norm{D_{2}U_{i}}_{p,\Omega_{i}\times\cell_{i}}}\norm{Du_{\varepsilon}^{\dual}}_{p^{\dual},\Omega_{i}}.
\end{aligned}
\]
The⋅rest is routine. We see that
\[
\begin{aligned}\hspace{2em} & \hspace{-2em}\varepsilon_{i}\bigabs{\rbr{\tau_{i}^{\varepsilon}\eta_{\varepsilon,i}A_{i}\op T_{i}^{\varepsilon}D_{1}\op J_{i}^{\varepsilon}U_{i},D_{1}u_{\varepsilon}^{\dual}}_{\Omega_{i}\times\cell_{i}}}\\
 & \le\varepsilon_{i}\norm{A_{i}}_{L_{\infty}}\norm{\tau_{i}^{\varepsilon}\op T_{i}^{\varepsilon}D_{1}\op J_{i}^{\varepsilon}U_{i}}_{p,H_{\varepsilon,i}\times\cell_{i}}\norm{D_{1}u_{\varepsilon}^{\dual}}_{p^{\dual},\Omega_{i}\times\cell_{i}}\\
 & \lesssim\varepsilon_{i}\norm{D_{1}\op J_{i}^{\varepsilon}U_{i}}_{p,\Omega_{i}\times\cell_{i}}\norm{Du_{\varepsilon}^{\dual}}_{p^{\dual},\Omega_{i}},
\end{aligned}
\]
while
\[
\begin{aligned}\hspace{2em} & \hspace{-2em}\bigabs{\rbr{\eta_{\varepsilon,i}A_{i}^{\varepsilon}\rbr{\op I-\op S_{i}^{\varepsilon}}\op J_{i}^{\varepsilon}\op e_{i}Du_{0,i},Du_{\varepsilon}^{\dual}}_{\Omega_{i}}}\\
 & \le\norm{A_{i}}_{L_{\infty}}\norm{\rbr{\op I-\op S_{i}^{\varepsilon}}\op J_{i}^{\varepsilon}\op e_{i}Du_{0,i}}_{p,H_{\varepsilon,i}}\norm{Du_{\varepsilon}^{\dual}}_{p^{\dual},\Omega_{i}}\\
 & \lesssim\varepsilon_{i}\norm{D\op J_{i}^{\varepsilon}\op e_{i}Du_{0,i}}_{p,\Omega_{i}}\norm{Du_{\varepsilon}^{\dual}}_{p^{\dual},\Omega_{i}}
\end{aligned}
\]
and
\[
\varepsilon_{i}\bigabs{\rbr{\eta_{\varepsilon,i}U_{\varepsilon,i},u_{\varepsilon}^{\dual}}_{\Omega_{i}}}\le\varepsilon_{i}\norm{U_{\varepsilon,i}}_{p,H_{\varepsilon,i}}\norm{u_{\varepsilon}^{\dual}}_{p^{\dual},\Omega_{i}}\lesssim\varepsilon_{i}\norm{U_{i}}_{p,\Omega_{i}\times\cell_{i}}\norm{u_{\varepsilon}^{\dual}}_{p^{\dual},\Omega_{i}}.
\]
Combining these bounds with Lemmas~\ref{lem: C | Operators R0_(mu,i) and K_(mu,i) bounded}
and⋅\ref{lem: C | Operators Je_i R0_(mu,i) and Je_i K_(mu,i) bounded}
and the estimate~(\dualref{assu: OO | Resolvent of Ae bounded}),
we get
\[
\bigabs{\rbr{\op I_{\mu,i}^{\varepsilon}f,f^{\dual}}_{\Omega}}\lesssim\omega_{i}\rbr{\varepsilon_{i}}\hairspace\norm f_{p,\Omega}\norm{f^{\dual}}_{-1,p^{\dual},\Omega}^{*},
\]
as⋅desired.
\end{proof}
\begin{lem}
\label{lem: Ps | Operator De_(i,mu) bounded}For⋅any⋅$\varepsilon\in\set E_{\mu}$
and⋅$f\in L_{p}\rbr{\Omega}^{n}$,
\begin{equation}
\norm{\op D_{\mu,i}^{\varepsilon}f}_{1,p,\Omega}\lesssim\omega_{i}\rbr{\varepsilon_{i}}\hairspace\norm f_{p,\Omega},\label{est: Ps | Operator De_(i,mu) bounded}
\end{equation}
where⋅$\op D_{\mu,i}^{\varepsilon}$ is given by~\textup{(\ref{def: Ps | Operator D_(mu,i)})}.
\end{lem}
\begin{proof}
The⋅bound follows from  Lemma~\ref{lem: C | Operators Je_i R0_(i,mu) and Je_i K_(i,mu) converge}
and the estimate~(\dualref{assu: OO | Resolvent of Ae bounded})
applied to the right-hand sides⋅of
\[
\begin{aligned}\hspace{2em} & \hspace{-2em}\bigabs{\rbr{\eta_{\varepsilon,i}\rbr{\op I-\op J_{i}^{\varepsilon}}A_{i}\rbr{\op e_{i}D_{1}u_{0,i}+D_{2}U_{i}},D_{1}u_{\varepsilon}^{\dual}}_{\Omega_{i}\times\cell_{i}}}\\
 & \le\norm{\rbr{\op I-\op J_{i}^{\varepsilon}}A_{i}\rbr{\op e_{i}D_{1}u_{0,i}+D_{2}U_{i}}}_{p,\Omega_{i}\times\cell_{i}}\norm{D_{1}u_{\varepsilon}^{\dual}}_{p^{\dual},\Omega_{i}\times\cell_{i}}
\end{aligned}
\]
and
\begin{align*}
\hspace{2em} & \hspace{-2em}\bigabs{\rbr{\eta_{\varepsilon,i}A_{i}^{\varepsilon}\rbr{\op I-\op J_{i}^{\varepsilon}}\op e_{i}Du_{0,i},Du_{\varepsilon}^{\dual}}_{\Omega_{i}}}\\
 & \le\norm{A_{i}}_{L_{\infty}}\norm{\rbr{\op I-\op J_{i}^{\varepsilon}}\op e_{i}Du_{0,i}}_{p,\Omega_{i}}\norm{Du_{\varepsilon}^{\dual}}_{p^{\dual},\Omega_{i}}.\qedhere
\end{align*}
\end{proof}
The⋅boundary term~$\op B_{\mu,i}^{\varepsilon}$ admits different
bounds in the two operator norms, and the⋅$L_{p}$\nobreakdash-bound
improves under additional assumptions (\cf~Theorems~\ref{thm: MR | Conv in L_q and approx in W_p^1}
and~\ref{thm: MR | Conv in L_p with better rate}). We therefore
treat the two cases separately. In⋅each⋅case, we need to control norms
over thin neighborhoods of the boundaries,  which is done by the
following lemma.

\begin{lem}
\label{lem: Ps | Norm over a thin layer}Let⋅$\Sigma$ be a uniformly
weakly Lipschitz domain in⋅$\R^{d}$ and⋅$\omega$ be a regular modulus
of continuity. Then, for any⋅$q\in[1,\infty)$, $\lambda>0$, $\delta\in(0,1]$
and any⋅$u\in C_{c}^{\infty}\rbr{\wb{\Sigma}}$,
\begin{equation}
\norm u_{q,\rbr{\partial\Sigma}_{\lambda\delta}\cap\Sigma}\lesssim\omega\rbr{\delta}^{1/q}\bigrbr{\delta\omega\rbr{\delta}^{-1}\norm{Du}_{q,\Sigma}+\norm u_{q,\Sigma}}.\label{est: Ps | Norm over a thin layer}
\end{equation}
The⋅implicit constant in the estimate depends only on⋅$q$, $d$,
$\lambda$ and~$\Sigma$.
\end{lem}
\begin{proof}
By~\cite[Lemma~A.1]{Senik:2023-1}, 
\[
\norm u_{q,\rbr{\partial\Sigma}_{\lambda\delta}\cap\Sigma}\lesssim\lambda^{1/q}\delta^{1/q}\norm u_{1,q,\Sigma}^{1/q}\norm u_{q,\Sigma}^{1/q^{\hconj}}
\]
with⋅a⋅constant depending on⋅$q$, $d$ and~$\Sigma$. Young's inequality
then gives
\[
\norm u_{q,\rbr{\partial\Sigma}_{\lambda\delta}\cap\Sigma}\lesssim\omega\rbr{\delta}^{1/q}\bigrbr{\delta\omega\rbr{\delta}^{-1}\norm u_{1,q,\Sigma}+\norm u_{q,\Sigma}}.
\]
It⋅remains to note that $\delta\omega\rbr{\delta}^{-1}$ is nondecreasing
and~$\omega\rbr 1=1$.
\end{proof}
We now derive an estimate involving⋅$\op B_{\mu,i}^{\varepsilon}$
that serves both cases. It⋅controls the boundary term, up to a correction,
through the gradient of the dual solution on a thin neighborhood
of~$\partial\Omega_{i}$. The⋅two cases then differ only in how the
corresponding norm is estimated: in the first, it is trivially bounded
by the same norm over all of⋅$\Omega_{i}$, and in the second, the
approximation for the adjoint problem, available under the additional
assumptions, yields an extra small factor.
\begin{lem}
\label{lem: Ps | Operator Be_(i,mu) wo last term bounded}For⋅any⋅$\varepsilon\in\set E_{\mu}$,
$f\in L_{p}\rbr{\Omega}^{n}$ and⋅$f^{\dual}\in\rbr{W_{p}^{1}\rbr{\Omega}^{n}}^{*}$,
\begin{equation}
\bigabs{\rbr{\op B_{\mu,i}^{\varepsilon}f+\varepsilon_{i}\rho_{\varepsilon,i}^{\prime}U_{\varepsilon,i},f^{\dual}}_{\Omega}}\lesssim\omega_{i}\rbr{\varepsilon_{i}}^{1/p}\norm f_{p,\Omega}\norm{Du_{\varepsilon}^{\dual}}_{p^{\dual},\rbr{\partial\Omega_{i}}_{5\varepsilon_{i}\cellradius_{i}}\cap\Omega_{i}},\label{est: Ps | Operator Be_i wo last term bounded}
\end{equation}
where⋅$\op B_{\mu,i}^{\varepsilon}$ is given by~\textup{(\ref{def: Ps | Operator B_(mu,i)})}.
\end{lem}
\begin{proof}
The⋅argument is similar to that of Lemma~\ref{lem: Ps | Operator Ie_(i,mu) bounded}.
We have
\[
\begin{aligned}\hspace{2em} & \hspace{-2em}\bigabs{\rbr{\tau_{i}^{\varepsilon}\op T_{i}^{\varepsilon}\op J_{i}^{\varepsilon}A_{i}\rbr{\op e_{i}D_{1}u_{0,i}+D_{2}U_{i}}\cdot\op T_{i}^{\varepsilon}\rbr{D_{1}\xi_{\varepsilon,i}}^{*},\rbr{\op T_{i}^{\varepsilon}-\op I}u_{\varepsilon}^{\dual}}_{\Omega_{i}\times\cell_{i}}}\\
 & \lesssim\norm{\op J_{i}^{\varepsilon}A_{i}\rbr{\op e_{i}D_{1}u_{0,i}+D_{2}U_{i}}}_{p,X_{\varepsilon,i}^{\prime}\times\cell_{i}}\norm{Du_{\varepsilon}^{\dual}}_{p^{\dual},\rbr{X_{\varepsilon,i}^{\prime}}_{2\varepsilon_{i}\cellradius_{i}}}
\end{aligned}
\]
and
\[
\begin{aligned}\hspace{2em} & \hspace{-2em}\bigabs{\rbr{\tau_{i}^{\varepsilon}\sbr{\xi_{\varepsilon,i},\op T_{i}^{\varepsilon}}A_{i}\op J_{i}^{\varepsilon}\rbr{\op e_{i}D_{1}u_{0,i}+D_{2}U_{i}},D_{1}u_{\varepsilon}^{\dual}}_{\Omega_{i}\times\cell_{i}}}\\
 & \lesssim\bigrbr{\norm{\op J_{i}^{\varepsilon}\op e_{i}Du_{0,i}}_{p,\rbr{X_{\varepsilon,i}^{\prime}}_{2\varepsilon_{i}\cellradius_{i}}}+\norm{D_{2}\op J_{i}^{\varepsilon}U_{i}}_{p,\rbr{X_{\varepsilon,i}^{\prime}}_{2\varepsilon_{i}\cellradius_{i}}\times\cell_{i}}}\norm{Du_{\varepsilon}^{\dual}}_{p^{\dual},\rbr{X_{\varepsilon,i}^{\prime}}_{\varepsilon_{i}\cellradius_{i}}}.
\end{aligned}
\]
The⋅remaining terms satisfy
\[
\begin{aligned}\hspace{2em} & \hspace{-2em}\bigabs{\rbr{\xi_{\varepsilon,i}\rbr{A_{i}^{0}-A_{i}^{\varepsilon}}Du_{0,i},Du_{\varepsilon}^{\dual}}_{\Omega_{i}}}\\
 & \lesssim\bigrbr{\norm{\rbr{\op I-\op J_{i}^{\varepsilon}}\op e_{i}Du_{0,i}}_{p,\Omega_{i}}+\norm{\op J_{i}^{\varepsilon}\op e_{i}Du_{0,i}}_{p,X_{\varepsilon,i}}}\norm{Du_{\varepsilon}^{\dual}}_{p^{\dual},X_{\varepsilon,i}}
\end{aligned}
\]
and
\[
\varepsilon_{i}\bigabs{\rbr{A_{i}^{\varepsilon}U_{\varepsilon,i}D\xi_{\varepsilon,i},Du_{\varepsilon}^{\dual}}_{\Omega_{i}}}\lesssim\norm{\op J_{i}^{\varepsilon}U_{i}}_{p,\rbr{X_{\varepsilon,i}^{\prime}}_{\varepsilon_{i}\cellradius_{i}}\times\cell_{i}}\norm{Du_{\varepsilon}^{\dual}}_{p^{\dual},X_{\varepsilon,i}^{\prime}}.
\]
These estimates together give
\[
\begin{aligned}\hspace{1em} & \hspace{-1em}\bigabs{\rbr{\op B_{\mu,i}^{\varepsilon}f+\varepsilon_{i}\rho_{\varepsilon,i}^{\prime}U_{\varepsilon,i},f^{\dual}}_{\Omega}}\\
 & \lesssim\!\bigl(\norm{\op J_{i}^{\varepsilon}A_{i}\rbr{\op e_{i}D_{1}u_{0,i}+D_{2}U_{i}}}_{p,\rbr{\partial\Omega_{i}}_{3\varepsilon_{i}\cellradius_{i}}\cap\Omega_{i}\times\cell_{i}}\\
 & \hphantom{{}\lesssim\!\bigl(}{}\!+\!\norm{\op J_{i}^{\varepsilon}\op e_{i}Du_{0,i}}_{p,\rbr{\partial\Omega_{i}}_{5\varepsilon_{i}\cellradius_{i}}\cap\Omega_{i}}\!+\!\norm{\rbr{\op I\!-\!\op J_{i}^{\varepsilon}}\op e_{i}Du_{0,i}}_{p,\Omega_{i}}\\
 & \hphantom{{}\lesssim\bigl({}\!+\!\norm{\op J_{i}^{\varepsilon}\op e_{i}Du_{0,i}}_{p,\rbr{\partial\Omega_{i}}_{5\varepsilon_{i}\cellradius_{i}}\cap\Omega_{i}}\!+\!\norm{\rbr{\op I\!-\!\op J_{i}^{\varepsilon}}\op e_{i}Du_{0,i}}_{p,\Omega_{i}}}\mathllap{{}+\norm{\op J_{i}^{\varepsilon}U_{i}}_{1_{2},p,\rbr{\partial\Omega_{i}}_{5\varepsilon_{i}\cellradius_{i}}\cap\Omega_{i}\times\cell_{i}}}\bigr)\norm{Du_{\varepsilon}^{\dual}}_{p^{\dual},\rbr{\partial\Omega_{i}}_{5\varepsilon_{i}\cellradius_{i}}\cap\Omega_{i}}.
\end{aligned}
\]
Using Lemma~\ref{lem: Ps | Norm over a thin layer} where needed
and then applying Lemmas~$\ref{lem: C | Operators R0_(mu,i) and K_(mu,i) bounded}$,
\ref{lem: C | Operators Je_i R0_(mu,i) and Je_i K_(mu,i) bounded}
and~\ref{lem: C | Operators Je_i R0_(i,mu) and Je_i K_(i,mu) converge},
we obtain~(\ref{est: Ps | Operator Be_i wo last term bounded}).
\end{proof}
As⋅a⋅consequence of this lemma, we have the following bound for~$\op B_{\mu,i}^{\varepsilon}$.
\begin{lem}
\label{lem: Ps | Operator Be_(i,mu) bounded in W_p^1}For⋅any⋅$\varepsilon\in\set E_{\mu}$
and⋅$f\in L_{p}\rbr{\Omega}^{n}$,
\begin{equation}
\norm{\op B_{\mu,i}^{\varepsilon}f}_{1,p,\Omega}\lesssim\omega_{i}\rbr{\varepsilon_{i}}^{1/p}\norm f_{p,\Omega},\label{est: Ps | Operator Be_(i,mu) bounded in W_p^1}
\end{equation}
where⋅$\op B_{\mu,i}^{\varepsilon}$ is given by~\textup{(\ref{def: Ps | Operator B_(mu,i)})}.
\end{lem}
\begin{proof}
We see from⋅Lemma~\ref{lem: Ps | Operator Be_(i,mu) wo last term bounded}
and the estimate~(\dualref{assu: OO | Resolvent of Ae bounded})
that
\[
\norm{\op B_{\mu,i}^{\varepsilon}f}_{1,p,\Omega}\lesssim\varepsilon_{i}\norm{\rho_{\varepsilon,i}^{\prime}U_{\varepsilon,i}}_{1,p,\Omega_{i}}+\omega_{i}\rbr{\varepsilon_{i}}^{1/p}\norm f_{p,\Omega}.
\]
For⋅the⋅first term on the right, recall that $\varepsilon_{i}D\tau_{i}^{\varepsilon}\op S_{i}^{\varepsilon}=\varepsilon_{i}\tau_{i}^{\varepsilon}\op S_{i}^{\varepsilon}D_{1}+\tau_{i}^{\varepsilon}\op S_{i}^{\varepsilon}D_{2}$.
Therefore, by the definition of⋅$U_{\varepsilon,i}$ and Lemma~\ref{lem: C | Operator =0003C4e Se bounded on domains},
\[
\begin{aligned}\varepsilon_{i}\norm{\rho_{\varepsilon,i}^{\prime}U_{\varepsilon,i}}_{1,p,\Omega_{i}} & \lesssim\varepsilon_{i}\norm{U_{\varepsilon,i}}_{1,p,\supp\rho_{\varepsilon,i}^{\prime}}+\norm{U_{\varepsilon,i}}_{p,\supp D\rho_{\varepsilon,i}^{\prime}}\\
 & \lesssim\varepsilon_{i}\norm{\op J_{i}^{\varepsilon}U_{i}}_{1_{1},p,\rbr{\Omega_{i}}_{\varepsilon_{i}\cellradius_{i}}\times\cell_{i}}+\norm{\op J_{i}^{\varepsilon}U_{i}}_{1_{2},p,\rbr{\supp\rho_{\varepsilon,i}^{\prime}}_{\varepsilon_{i}\cellradius_{i}}\times\cell_{i}}.
\end{aligned}
\]
Note that $\rho_{\varepsilon,i}^{\prime}$ is supported in the closure
of⋅$\rbr{\partial\Omega_{i}}_{3\varepsilon_{i}\cellradius_{i}}\cap\Omega_{i}$.
We then split the $\varepsilon_{i}\cellradius_{i}$\nobreakdash-neighborhood
of this closure into the parts inside and outside~$\Omega_{i}$,
the sets⋅$\rbr{\partial\Omega_{i}}_{4\varepsilon_{i}\cellradius_{i}}\cap\Omega_{i}$
and⋅$\rbr{\Omega_{i}}_{\varepsilon_{i}\cellradius_{i}}\negmedspace\setminus\Omega_{i}$.
Lemma~\ref{lem: Ps | Norm over a thin layer}, with⋅$\Sigma=\Omega_{i}$
and with⋅$\Sigma=\R^{d}\negmedspace\setminus\wb{\Omega}_{i}$, now
shows for the second term on the far right that
\begin{equation}
\varepsilon_{i}\norm{\rho_{\varepsilon,i}^{\prime}U_{\varepsilon,i}}_{1,p,\Omega_{i}}\lesssim\omega_{i}\rbr{\varepsilon_{i}}^{1/p}\bigrbr{\varepsilon_{i}\omega_{i}\rbr{\varepsilon_{i}}^{-1}\norm{D_{1}\op J_{i}^{\varepsilon}U_{i}}_{1_{2},p,\R^{d}\times\cell_{i}}+\norm{\op J_{i}^{\varepsilon}U_{i}}_{1_{2},p,\R^{d}\times\cell_{i}}}.\label{est: Norm of DUe near boundary}
\end{equation}
With⋅Lemmas~\ref{lem: C | Operators R0_(mu,i) and K_(mu,i) bounded}
and⋅\ref{lem: C | Operators Je_i R0_(mu,i) and Je_i K_(mu,i) bounded},
this becomes
\[
\varepsilon_{i}\norm{\rho_{\varepsilon,i}^{\prime}U_{\varepsilon,i}}_{1,p,\Omega_{i}}\lesssim\omega_{i}\rbr{\varepsilon_{i}}^{1/p}\norm f_{p,\Omega},
\]
and the proof is complete.
\end{proof}
The⋅above lemma is the last ingredient needed for the proof of
Theorem~\ref{thm: MR | Conv in L_q and approx in W_p^1}.
\begin{proof}[Proof of Theorem~\textup{\ref{thm: MR | Conv in L_q and approx in W_p^1}}]
From⋅the⋅identity~(\ref{eq: Ps | Resolvent identity}) and Lemmas~\ref{lem: Ps | Operator Ie_(i,mu) bounded},
\ref{lem: Ps | Operator De_(i,mu) bounded} and~\ref{lem: Ps | Operator Be_(i,mu) bounded in W_p^1},
we conclude that
\[
\norm{\rbr{\op A_{\mu}^{\varepsilon}}^{-1}f-\rbr{\op A_{\mu}^{0}}^{-1}f-\varepsilon_{\#}\op K_{\mu}^{\varepsilon}f}_{1,p,\Omega}\lesssim\varpi_{\join}\rbr{\varepsilon}^{1/p}\norm f_{p,\Omega}.
\]
This immediately gives~(\ref{est: MR | Approx in W_p^1}). The⋅bound⋅(\ref{est: MR | Conv in L_q})
also follows, since, by the Sobolev⋅embedding theorem and the triangle
inequality,
\[
\norm{\rbr{\op A_{\mu}^{\varepsilon}}^{-1}f-\rbr{\op A_{\mu}^{0}}^{-1}f}_{q,\Omega}\lesssim\norm{\rbr{\op A_{\mu}^{\varepsilon}}^{-1}f-\rbr{\op A_{\mu}^{0}}^{-1}f-\varepsilon_{\#}\op K_{\mu}^{\varepsilon}f}_{1,p,\Omega}+\norm{\varepsilon_{\#}\op K_{\mu}^{\varepsilon}f}_{q,\Omega},
\]
where  the last term on the right is estimated using⋅(\ref{est: C | Corrector Ke bounded in L_q}),
 since~$\varpi_{\join}\rbr{\varepsilon}\le1$.
\end{proof}
For⋅Theorem~\ref{thm: MR | Conv in L_p with better rate}, the
key observation is that the norm of⋅$Du_{\varepsilon}^{\dual}$ on
the right-hand side of⋅(\ref{est: Ps | Operator Be_i wo last term bounded})
is taken only over a thin neighborhood of~$\partial\Omega_{i}$.
Lemma~\ref{lem: Ps | Norm over a thin layer} cannot be applied here
directly: $u_{\varepsilon}^{\dual}$ need not be regular enough, and
even if it were, its higher norms would blow up as~$\varepsilon\to0$.
Under the hypotheses of Theorem~\dualref{thm: MR | Conv in L_q and approx in W_p^1},
however, $u_{\varepsilon}^{\dual}$ can be replaced near the boundary
by the sum of the more regular functions~$u_{0,\#}^{\dual}$ and⋅$\varepsilon_{\#}\chi_{\varepsilon}\smash[t]{U_{\varepsilon,\#}^{\dual}}$
with an error of the right order (this idea is due to Griso~\cite[Lemma~3.1]{Gri:2006}).
We carry this out in the next lemma, which gives a bound for⋅$\op B_{\mu,i}^{\varepsilon}$
in the operator norm on~$L_{p}$.
\begin{lem}
\label{lem: Ps | Operator Be_(i,mu) bounded in L_p}Suppose that
the hypotheses of Theorem~\textup{\dualref{thm: MR | Conv in L_q and approx in W_p^1}}
hold. Then, for any⋅$\varepsilon\in\set E_{\mu}$ and⋅$f\in L_{p}\rbr{\Omega}^{n}$,
\begin{equation}
\norm{\op B_{\mu,i}^{\varepsilon}f}_{p,\Omega}\lesssim\omega_{i}\rbr{\varepsilon_{i}}^{1/p}\varpi_{\join}^{\dual}\rbr{\varepsilon}^{1/p^{\dual}}\norm f_{p,\Omega},
\end{equation}
where $\op B_{\mu,i}^{\varepsilon}$ is given by~\textup{(\ref{def: Ps | Operator B_(mu,i)})}.
\end{lem}
\begin{proof}
Taking⋅$f^{\dual}\in L_{p^{\dual}}\rbr{\Omega}^{n}$ in⋅(\ref{est: Ps | Operator Be_i wo last term bounded})
and using⋅(\ref{est: C | Corrector Ke bounded}), we get
\[
\begin{aligned}\bigabs{\rbr{\op B_{\mu,i}^{\varepsilon}f,f^{\dual}}_{\Omega}} & \lesssim\norm f_{p,\Omega}\bigrbr{\omega_{i}\rbr{\varepsilon_{i}}^{1/p}\norm{Du_{\varepsilon}^{\dual}}_{p^{\dual},\rbr{\partial\Omega_{i}}_{5\varepsilon_{i}\cellradius_{i}}\cap\Omega_{i}}+\varepsilon_{i}\norm{f^{\dual}}_{p^{\dual},\Omega}}.\end{aligned}
\]
It⋅remains to show that
\begin{equation}
\norm{Du_{\varepsilon}^{\dual}}_{p^{\dual},\rbr{\partial\Omega_{i}}_{5\varepsilon_{i}\cellradius_{i}}\cap\Omega_{i}}\lesssim\varpi_{\join}^{\dual}\rbr{\varepsilon}^{1/p^{\dual}}\norm{f^{\dual}}_{p^{\dual},\Omega}.\label{est: Ps | Du_=00005Cvarepsilon^+ in a thin layer}
\end{equation}
To⋅this end, we approximate⋅$Du_{\varepsilon}^{\dual}$ by the
sum of⋅$\op J_{i}^{\varepsilon}\op e_{i}Du_{0,i}^{\dual}$ and⋅$\varepsilon_{i}D\chi_{\varepsilon}U_{\varepsilon,i}^{\dual}$,
estimating the error by⋅Theorem~\dualref{thm: MR | Conv in L_q and approx in W_p^1}:
\[
\begin{aligned}\norm{Du_{\varepsilon}^{\dual}}_{p^{\dual},\rbr{\partial\Omega_{i}}_{5\varepsilon_{i}\cellradius_{i}}\cap\Omega_{i}} & \lesssim\norm{\op J_{i}^{\varepsilon}\op e_{i}Du_{0,i}^{\dual}}_{p^{\dual},\rbr{\partial\Omega_{i}}_{5\varepsilon_{i}\cellradius_{i}}\cap\Omega_{i}}+\varepsilon_{i}\norm{D\chi_{\varepsilon}U_{\varepsilon,i}^{\dual}}_{p^{\dual},\rbr{\partial\Omega_{i}}_{5\varepsilon_{i}\cellradius_{i}}\cap\Omega_{i}}\\
 & \quad+\norm{\rbr{\op I-\op J_{i}^{\varepsilon}}\op e_{i}Du_{0,i}^{\dual}}_{p^{\dual},\Omega_{i}}+\varpi_{\join}^{\dual}\rbr{\varepsilon}^{1/p^{\dual}}\norm{f^{\dual}}_{p^{\dual},\Omega}.
\end{aligned}
\]
By⋅Lemma~\ref{lem: Ps | Norm over a thin layer}, the first two
terms on the right satisfy
\[
\begin{aligned}\hspace{2em} & \hspace{-2em}\norm{\op J_{i}^{\varepsilon}\op e_{i}Du_{0,i}^{\dual}}_{p^{\dual},\rbr{\partial\Omega_{i}}_{5\varepsilon_{i}\cellradius_{i}}\cap\Omega_{i}}\\
 & \lesssim\omega_{i}^{\dual}\rbr{\varepsilon_{i}}^{1/p^{\dual}}\bigrbr{\varepsilon_{i}\omega_{i}^{\dual}\rbr{\varepsilon_{i}}^{-1}\norm{D\op J_{i}^{\varepsilon}\op e_{i}Du_{0,i}^{\dual}}_{p^{\dual},\Omega_{i}}+\norm{Du_{0,i}^{\dual}}_{p^{\dual},\Omega_{i}}}
\end{aligned}
\]
and
\[
\begin{aligned}\hspace{2em} & \hspace{-2em}\varepsilon_{i}\norm{D\chi_{\varepsilon}U_{\varepsilon,i}^{\dual}}_{p^{\dual},\rbr{\partial\Omega_{i}}_{5\varepsilon_{i}\cellradius_{i}}\cap\Omega_{i}}\\
 & \lesssim\omega_{i}^{\dual}\rbr{\varepsilon_{i}}^{1/p^{\dual}}\bigrbr{\varepsilon_{i}\omega_{i}^{\dual}\rbr{\varepsilon_{i}}^{-1}\norm{D_{1}\op J_{i}^{\varepsilon}U_{i}^{\dual}}_{1_{2},p^{\dual},\R^{d}\times\cell_{i}}+\norm{\op J_{i}^{\varepsilon}U_{i}^{\dual}}_{1_{2},p^{\dual},\R^{d}\times\cell_{i}}}
\end{aligned}
\]
(\cf~(\ref{est: Norm of DUe near boundary}), where a neighborhood
of⋅$\partial\Omega_{i}$ is split in the same way). These bounds,
together with Lemmas~\dualref{lem: C | Operators R0_(mu,i) and K_(mu,i) bounded},
\dualref{lem: C | Operators Je_i R0_(mu,i) and Je_i K_(mu,i) bounded}
and~\dualref{lem: C | Operators Je_i R0_(i,mu) and Je_i K_(i,mu) converge},
yield the desired estimate~(\ref{est: Ps | Du_=00005Cvarepsilon^+ in a thin layer}).
\end{proof}
\begin{proof}[Proof of Theorem~\textup{\ref{thm: MR | Conv in L_p with better rate}}]
The theorem follows at once from the identity~(\ref{eq: Ps | Resolvent identity})
and the estimate~(\ref{est: C | Corrector Ke bounded}), combined
with Lemmas~\ref{lem: Ps | Operator Ie_(i,mu) bounded}, \ref{lem: Ps | Operator De_(i,mu) bounded}
and~\ref{lem: Ps | Operator Be_(i,mu) bounded in L_p}, since $\varpi_{\join}\le\varpi_{\join}^{\dual}$
by our convention.
\end{proof}
We turn to the interior estimate. Its proof is based on the localized
resolvent identity~(\ref{eq: Localized resolvent identity}), whose
terms are estimated by the  analogs of Lemmas~\ref{lem: Ps | Operator Ie_(i,mu) bounded}
and~\ref{lem: Ps | Operator De_(i,mu) bounded} stated below, with
implicit constants depending also on~$\eta$. The⋅arguments are the
same and are omitted.
\begin{lem}
\label{lem: Ps | Operator I0e_(i,mu) bounded}For⋅any⋅$\varepsilon\in\set E_{\mu}$
and⋅$f\in L_{p}\rbr{\Omega}^{n}$,
\begin{equation}
\norm{\mathring{\op I}_{\mu,i}^{\varepsilon}f}_{1,p,\Omega}\lesssim\omega_{i}\rbr{\varepsilon_{i}}\hairspace\norm f_{p,\Omega},\label{est: Ps | Operator I0e_(i,mu) bounded}
\end{equation}
where⋅$\mathring{\op I}_{\mu,i}^{\varepsilon}$ is given by~\textup{(\ref{def: Ps | Operator I0_(mu,i)})}.
\end{lem}
\begin{lem}
\label{lem: Ps | Operator D0e_(mu,i) bounded}For⋅any⋅$\varepsilon\in\set E_{\mu}$
and⋅$f\in L_{p}\rbr{\Omega}^{n}$,
\begin{equation}
\norm{\mathring{\op D}_{\mu,i}^{\varepsilon}f}_{1,p,\Omega}\lesssim\omega_{i}\rbr{\varepsilon_{i}}\hairspace\norm f_{p,\Omega},\label{est: Ps | Operator D0e_(i,mu) bounded}
\end{equation}
where⋅$\mathring{\op D}_{\mu,i}^{\varepsilon}$ is given by~\textup{(\ref{def: Ps | Operator D0_(mu,i)})}.
\end{lem}
We are now ready to prove Theorem~\ref{thm: MR | Interior approx in W_p^1}.
\begin{proof}[Proof of Theorem~\textup{\ref{thm: MR | Interior approx in W_p^1}}]
Set⋅$\eta=\wb{\chi}^{\prime}$. By⋅assumption, $\eta$ vanishes
in the⋅$\delta$\nobreakdash-neighborhood of⋅$\partial\Omega\cup\Gamma$
for some~$\delta>0$. As⋅the⋅estimate is trivial for⋅$\max_{i\in\set I}5\varepsilon_{i}r_{i}>\delta$,
we may assume that⋅$\max_{i\in\set I}5\varepsilon_{i}r_{i}\le\delta$,
in which case the localized resolvent identity~(\ref{eq: Localized resolvent identity})
applies. Fix⋅$\eta^{\prime}\in C^{0,1}\rbr{\wb{\Omega}}$ that vanishes
on⋅$\partial\Omega\cup\Gamma$  and equals~$1$ near⋅$\supp\eta$.
Let⋅$v_{\varepsilon}=u_{\varepsilon}-u_{0,\#}-\varepsilon_{\#}\chi_{\varepsilon}U_{\varepsilon,\#}$
and⋅$g_{\varepsilon}=\smash{\wb{\eta}\mathring{\op A}_{\mu}^{\varepsilon}v_{\varepsilon}}$.
Inserting⋅$\eta^{\prime}$ before⋅$v_{\varepsilon}$, we can also
write⋅$g_{\varepsilon}=\smash{\wb{\eta}\op A_{\mu}^{\varepsilon}\eta^{\prime}v_{\varepsilon}}$.
Then, by the interior energy estimate~(\ref{assu: OO | Interior energy estimate})
with⋅$u=\eta^{\prime}v_{\varepsilon}$,
\[
\norm{\chi Dv_{\varepsilon}}_{p,\Omega}\lesssim\norm{v_{\varepsilon}}_{p,\Omega}+\norm{g_{\varepsilon}}_{-1,p,\Omega}^{*}.
\]
The⋅first term on the right is estimated by⋅(\ref{est: C | Corrector Ke bounded})
and Theorem~\ref{thm: MR | Conv in L_p with better rate}:
\[
\norm{v_{\varepsilon}}_{p,\Omega}\lesssim\varpi_{\join}\rbr{\varepsilon}^{1/p}\varpi_{\join}^{\dual}\rbr{\varepsilon}^{1/p^{\dual}}\norm f_{p,\Omega}.
\]
For⋅the⋅other term, the localized resolvent identity~(\ref{eq: Localized resolvent identity})
and Lemmas~\ref{lem: Ps | Operator I0e_(i,mu) bounded} and⋅\ref{lem: Ps | Operator D0e_(mu,i) bounded}
show that $w_{\varepsilon}=\smash[b]{\rbr{\op A_{\mu}^{\varepsilon}}^{-1}}g_{\varepsilon}$
satisfies
\[
\norm{w_{\varepsilon}}_{1,p,\Omega}\lesssim\varpi_{\join}\rbr{\varepsilon}\hairspace\norm f_{p,\Omega},
\]
whence, using⋅(\ref{est: OO | Operator Ae bounded}) and the identity~$g_{\varepsilon}=\eta^{\prime}g_{\varepsilon}$,
we get
\[
\norm{g_{\varepsilon}}_{-1,p,\Omega}^{*}\lesssim\Norm{g_{\varepsilon}}_{-1,p,\Omega}\lesssim\varpi_{\join}\rbr{\varepsilon}\hairspace\norm f_{p,\Omega}.
\]
It⋅remains to note that~$\varpi_{\join}\le\varpi_{\join}^{\dual}$.
\end{proof}

\appendix

\section{\label{sec: Appendix - Elliptic regularity}Elliptic regularity in
\texorpdfstring{$\Lambda_{p}^{1+\omega}$}{Λ\textunderscore p\textasciicircum (1+ω)}-spaces}

In⋅this appendix, we prove an elliptic regularity result in piecewise
$\Lambda_{p}^{1+\omega}$-spaces for a special case of the effective
transmission problem studied in the main text. We begin by introducing
our assumptions and formulating the problem.

Let⋅$\Omega$ be a bounded domain divided into a finite number of⋅$C^{1,\omega_{i}}$-regular
subdomains~$\Omega_{i}$,⋅$i\in\set I$. Then the interface~$\Gamma_{ij}$
between two adjacent subdomains~$\Omega_{i}$ and⋅$\Omega_{j}$,
\ie, the common part of the boundaries of the two, is  $C^{1,\omega_{ij}}$-regular
with~$\omega_{ij}=\omega_{i}\meet\omega_{j}$.   Denote by⋅$\smash{\set I_{\Gamma}^{2}}$
the set of   pairs~$\rbr{i,j}\in\set I\times\set I$ with⋅$i\ne j$
for which⋅$\Gamma_{ij}$ is nonempty and  by⋅$\Gamma$ the union
of the interfaces~$\Gamma_{ij}$ over⋅$\rbr{i,j}\in\smash{\set I_{\Gamma}^{2}}$.
We suppose that there are neither interior nor boundary cross-points.
(Recall that this means that no two interfaces meet and none of them
touches the boundary; see⋅Figure~\ref{subfig: OO | Smooth domains}
for a typical geometry.) This assumption was discussed in Remark~\ref{rem: Ex | Cross-points}.

To⋅simplify the notation, we drop the superscript from the effective
coefficient and consider a function~$A\colon\Omega\to\B\rbr{\C^{n\times d}}$
satisfying, for each⋅$i\in\set I$,
\begin{equation}
A\in C^{0,\omega_{i}}\rbr{\wb{\Omega}_{i}}\label{assu: Appendix| Coefficient A regular in Omega_i}
\end{equation}
and,⋅moreover, for each⋅$\rbr{i,j}\in\smash{\set I_{\Gamma}^{2}}$,
\begin{equation}
A\in C^{0,\omega_{ij}}\rbr{\overline{\Omega_{i}\cap G_{ij}}}\label{assu: Appendix| Coefficient A regular near Gamma_ij}
\end{equation}
with⋅some open neighborhood~$G_{ij}=G_{ji}$ of~$\Gamma_{ij}$.
Without loss of generality, each coordinate neighborhood~$W$
of⋅$\Gamma_{ij}$  may be chosen so that~$\smash[t]{\overline{W}}\subset G_{ij}$.
As⋅noted in Remark~\ref{rem: Ex | On coefficient A additional regularity},
the assumption~(\ref{assu: Appendix| Coefficient A regular near Gamma_ij})
can always be met by enlarging the moduli of continuity, at the cost
of weaker regularity.

We will discuss the Dirichlet and  Neumann problems, which are all
we need  in Section~\ref{sec: Examples}.  Set⋅$\smash{\set W_{q}^{1}}\rbr{\Omega;\C^{n}}$
to be⋅$\smash{\mathring{W}_{q}^{1}}\rbr{\Omega}^{n}$ for the former
and⋅$\smash{W_{q}^{1}}\rbr{\Omega}^{n}$ for the latter.  Fix⋅$p\in\rbr{1,\infty}$
and define the operator~$\mathbb{A}\colon\smash{\set W_{p}^{1}}\rbr{\Omega;\C^{n}}\to\smash{\set W_{p}^{-1}}\rbr{\Omega;\C^{n}}$⋅by
\begin{equation}
\rbr{\mathbb{A}u,u^{\prime}}_{\Omega}=\rbr{ADu,Du^{\prime}}_{\Omega},\qquad u\in\set W_{p}^{1}\rbr{\Omega;\C^{n}},u^{\prime}\in\set W_{p^{\hconj}}^{1}\rbr{\Omega;\C^{n}}.\label{def: Appendix | Operator A}
\end{equation}
 We assume that⋅$\smash[b]{\mathbb{A}_{\mu}}=\mathbb{A}-\mu$
 is an isomorphism for some~$\mu\in\C$. 

Our main result here is the following.
\begin{thm}
\label{thm: Appendix | Elliptic regularity}Under⋅the⋅above hypotheses,
suppose that $f\in L_{p}\rbr{\Omega}^{n}$ and⋅$F\in L_{p}\rbr{\Omega}^{n\times d}$
are such that⋅$F\in\Lambda_{p}^{\omega_{i}}\rbr{\Omega_{i}}$ for
each⋅$i\in\set I$ and⋅$F\in\Lambda_{p}^{\smash[b]{\omega_{ij}}}\rbr{\Omega_{i}\cap G_{ij}}$
for each⋅$\rbr{i,j}\in\smash{\set I_{\Gamma}^{2}}$. Let⋅$u$ be the
solution of the problem
\begin{equation}
\mathbb{A}_{\mu}u=f+D^{*}F.\label{eq: Appendix | Elliptic regularity problem}
\end{equation}
Then, for each~$i\in\set I$,
\begin{equation}
\seminorm{Du}_{\omega_{i},p,\Omega_{i}}\lesssim\norm f_{p,\Omega}+\norm F_{p,\Omega}+\seminorm F_{\omega_{i},p,\Omega_{i}}+\sum_{\rbr{i,j}\in\set I_{\Gamma}^{2}}\bigrbr{\seminorm F_{\omega_{ij},p,\Omega_{i}\cap G_{ij}}+\seminorm F_{\omega_{ij},p,\Omega_{j}\cap G_{ij}}}.\label{est: Appendix | Elliptic regularity}
\end{equation}
Moreover, for each⋅$\rbr{i,j}\in\smash{\set I_{\Gamma}^{2}}$ and
every open neighborhood~$G_{ij}^{\prime}$ of⋅$\Gamma_{ij}$ compactly
contained in~$G_{ij}$,
\begin{equation}
\seminorm{Du}_{\omega_{ij},p,\Omega_{i}\cap G_{ij}^{\prime}}\lesssim\norm f_{p,\Omega}+\norm F_{p,\Omega}+\seminorm F_{\omega_{ij},p,\Omega_{i}\cap G_{ij}}+\seminorm F_{\omega_{ij},p,\Omega_{j}\cap G_{ij}}.\label{est: Appendix | Elliptic regularity near the interfaces}
\end{equation}
The⋅implicit constants  depend only on the parameters~$d$, $p$,
$n$, $\mu$, the $C^{1,\omega_{i}}$\nobreakdash-structures of⋅$\Omega_{i}$
 and the $C^{1,\omega_{ij}}$\nobreakdash-structures of⋅$\Gamma_{ij}$,
the $C^{0,\omega_{i}}$\nobreakdash-norms of⋅$A\restrict{\Omega_{i}}$
 and the $\smash[t]{C^{0,\omega_{ij}}}$\nobreakdash-norms of⋅$\smash[b]{A\restrict{\Omega_{i}\cap G_{ij}}}$,
the sets~$G_{ij}$ and, in⋅\textup{(\ref{est: Appendix | Elliptic regularity near the interfaces})},⋅$G_{ij}^{\prime}$,
 and the norm of~$\smash{\mathbb{A}_{\mu}^{-1}}$.
\end{thm}
\begin{rem}
Recall that, given⋅$u\in\set W_{p}^{1}\rbr{\Omega;\C^{n}}$ and a
functional~$\ell$ on⋅$\smash[b]{W_{p^{\prime}}^{1}}\rbr{\Omega}^{n}$
that coincides with⋅$\smash[b]{\mathbb{A}_{\mu}}u$ on⋅$\smash{\mathring{W}_{p^{\prime}}^{1}}\rbr{\Omega}^{n}$,
the generalized conormal derivative~$\smash[b]{\partial_{\nu_{\mathbb{A}}}}u\in\smash{W_{p}^{-1/p}}\rbr{\partial\Omega}^{n}$
of⋅$u$ relative to⋅$\ell$ is defined by the first Green identity
\begin{equation}
\rbr{ADu,Du^{\prime}}_{\Omega}-\mu\rbr{u,u^{\prime}}_{\Omega}=\rbr{\ell,u^{\prime}}_{\Omega}+\rbr{\partial_{\nu_{\mathbb{A}}}u,\gamma u^{\prime}}_{\partial\Omega},\qquad u^{\prime}\in W_{p^{\prime}}^{1}\rbr{\Omega}^{n}.\label{eq: Appendix | Green's identity}
\end{equation}
Note that, despite the notation, $\smash[b]{\partial_{\nu_{\mathbb{A}}}}u$
depends not only on⋅$u$ but also on the choice of⋅$\ell$, which
may be any extension of⋅$\smash[b]{\mathbb{A}_{\mu}}u$ up to the
boundary.  If⋅$\smash[b]{\mathbb{A}_{\mu}}u$ is an⋅$L_{p}$\nobreakdash-function
in⋅$\Omega$, we always take~$\ell$ to be its extension to the boundary
by zero.
\end{rem}
\begin{rem}
Both Dirichlet and Neumann boundary data can be encoded in the right-hand
side of~(\ref{eq: Appendix | Elliptic regularity problem}). Consider
the Dirichlet problem with⋅$\gamma u=\gamma u_{0}$, where⋅$u_{0}\in\smash{\Lambda_{p}^{1+\smash[b]{\omega_{\#}}}}\rbr{\Omega_{\#}}^{n}$.
Since the interfaces do not touch the boundary, multiplying⋅$u_{0}$
by a suitable cutoff function, we may assume that $u_{0}$ vanishes
near the interfaces. Then $u-u_{0}\in\smash{\mathring{W}_{p}^{1}}\rbr{\Omega}^{n}$
satisfies an equation of the form~(\ref{eq: Appendix | Elliptic regularity problem}),
because $\mathbb{A}_{\mu}u_{0}$ is the sum of an⋅$L_{p}$\nobreakdash-function
and the divergence of a $\smash{\Lambda_{p}^{\omega_{\#}}}$\nobreakdash-field.
As⋅for⋅the⋅Neumann problem, the equation~(\ref{eq: Appendix | Elliptic regularity problem})
imposes the boundary condition~$\gamma_{n}\rbr{ADu-F}=0$ with⋅$\gamma_{n}$
denoting the normal trace operator for~$\partial\Omega$. The⋅normal
trace is well defined as an element of⋅$W_{p}^{-\smash[b]{1/p}}\rbr{\partial\Omega}^{n}$,
since $D^{*}\rbr{ADu-F}=f+\mu u$ is an⋅$L_{p}$\nobreakdash-function
in~$\Omega$. If⋅so is⋅$D^{*}F$, then the right-hand side of⋅(\ref{eq: Appendix | Elliptic regularity problem})
is also an⋅$L_{p}$\nobreakdash-function in⋅$\Omega$,   which
fixes the generalized conormal derivative~$\smash[b]{\partial_{\nu_{\mathbb{A}}}}u=i\gamma_{n}ADu$,
and the condition becomes~$\smash[b]{\partial_{\nu_{\mathbb{A}}}}u=i\gamma_{n}F$.
Neumann data in⋅$W_{p}^{1-\smash[b]{1/p}}\rbr{\partial\Omega}^{n}$,
which are regular enough for every modulus of continuity~$\omega$,
can be represented by a right-hand side as in~(\ref{eq: Appendix | Elliptic regularity problem}),
see⋅Lemma~\ref{lem: Appendix | Regular Neumann data fit} below.
We use this form of the data, rather than boundary values, to avoid
trace and extension theorems for~$\Lambda_{p}^{1+\omega}$\nobreakdash-spaces.
\end{rem}
\begin{lem}
\label{lem: Appendix | Regular Neumann data fit}For every⋅$\varrho\in\smash{W_{p}^{1-\smash[b]{1/p}}}\rbr{\partial\Omega}^{n}$,
there are functions~$h\in L_{p}\rbr{\Omega}^{n}$ and⋅$H\in\Lambda_{p}^{\omega_{\#}}\rbr{\Omega_{\#}}^{n\times d}$
vanishing near⋅$\Gamma$ such that
\[
\gamma^{*}\varrho=h+D^{*}H
\]
as⋅functionals on⋅$W_{p^{\prime}}^{1}\rbr{\Omega}^{n}$,⋅and
\[
\norm h_{p,\Omega}+\abs H_{\omega_{\#},p,\Omega_{\#}}\lesssim\norm{\varrho}_{1-1/p,p,\partial\Omega}.
\]
\end{lem}
We postpone the proof to the end of the appendix.
\begin{rem}
The⋅class of data in Theorem~\ref{thm: Appendix | Elliptic regularity}
is optimal:  for a right-hand side in⋅$\smash{\set W_{p}^{-1}}\rbr{\Omega;\C^{n}}$,
the solution~$u$ satisfies the conclusion of the theorem if and
only if this right-hand side can be written as⋅$f+D^{*}F$, where⋅$f$
and⋅$F$ satisfy its assumptions (possibly with smaller neighborhoods~$G_{ij}$).
For⋅the ``only if'' part, one can take⋅$f=\nobreak-\mu u$ and⋅$F=ADu$.
\end{rem}
We point out that the operator~$\mathbb{A}$ in Theorem~\ref{thm: Appendix | Elliptic regularity}
is not assumed to be strongly elliptic. Nevertheless, the coefficients
of operators that are merely invertible satisfy a close analog of
the strong ellipticity condition, as in~(\ref{est: Ex | Legendre-Hadamard}).
This property plays a role in the proof of Theorem~\ref{thm: Appendix | Elliptic regularity},
and we state it in greater generality than needed here.
\begin{lem}
\label{lem: Appendix | Weaker ellipticity condition}Let⋅$\Sigma$
be an open set in⋅$\R^{d}$ and⋅$A\in L_{\infty}\rbr{\Sigma;\B\rbr{\C^{n\times d}}}$.
Suppose that
\begin{equation}
\norm{D^{*}ADu}_{-1,p,\Sigma}^{*}\ge c_{A}\norm{Du}_{p,\Sigma}-C_{A}\norm u_{p,\Sigma},\qquad u\in C_{c}^{\infty}\rbr{\Sigma}^{n},\label{est: Appendix | A_mu bounded below}
\end{equation}
where⋅$D^{*}ADu$ is regarded as a functional on~$W_{p^{\prime}}^{1}\rbr{\Sigma}^{n}$.
Then, for⋅\aev~$x\in\Sigma$, we have
\begin{equation}
\abs{\abr{A\rbr x\hairspace\eta\otimes\xi,\xi}}\ge c_{A}\abs{\eta}\hairspace\abs{\xi}^{2},\qquad\eta\in\C^{n},\xi\in\R^{d}.\label{est: Appendix | Weaker ellipticity condition}
\end{equation}
\end{lem}
\begin{proof}
First assume that $A$ is constant, denoted by~$A_{0}$.  Let⋅$\varepsilon>0$
be small, and let⋅$e_{\xi}^{\varepsilon}\rbr x=e^{i\abr{\xi,x/\varepsilon}}$
with a unit vector~$\xi\in\R^{d}$. Take⋅$u_{\varepsilon}=\varepsilon\varphi e_{\xi}^{\varepsilon}\eta$,
where⋅$\varphi\in C_{c}^{\infty}\rbr{\Sigma}$ and~$\eta\in\C^{n}$.
Then, by a simple calculation,
\[
D^{*}A_{0}Du_{\varepsilon}=\rbr{\varepsilon^{-1}\abr{A_{0}\xi,\xi}\varphi+\psi_{0}+\varepsilon\psi_{1}}e_{\xi}^{\varepsilon}\eta.
\]
Here $\abr{A_{0}\xi,\xi}=\sum_{k,l}\rbr{A_{0}}_{kl}\xi_{k}\xi_{l}$,
so that $\abr{A_{0}\eta\otimes\xi,\xi}=\abr{A_{0}\xi,\xi}\eta$,
and $\psi_{0}$ and~$\psi_{1}$ are smooth functions. We apply the
identity~$e_{\xi}^{\varepsilon}=\varepsilon D^{*}e_{\xi}^{\varepsilon}\xi$
to convert oscillating terms of the form~$\psi e_{\xi}^{\varepsilon}$
into the sum of a scaled divergence term and a small remainder of
the same form:
\[
\psi e_{\xi}^{\varepsilon}=\varepsilon D^{*}\psi e_{\xi}^{\varepsilon}\xi+\varepsilon\psi^{\prime}e_{\xi}^{\varepsilon},
\]
where⋅$\psi^{\prime}=-D^{*}\psi\xi$. Then
\[
D^{*}A_{0}Du_{\varepsilon}=D^{*}\abr{A_{0}\xi,\xi}\varphi e_{\xi}^{\varepsilon}\eta\otimes\xi+\varepsilon D^{*}\wt{\psi}_{0}e_{\xi}^{\varepsilon}\eta\otimes\xi+\varepsilon\wt{\psi}_{1}e_{\xi}^{\varepsilon}\eta
\]
with⋅smooth⋅$\smash{\wt{\psi}_{0}}$ and~$\smash{\wt{\psi}_{1}}$,
and⋅so
\begin{equation}
\limsup_{\varepsilon\to0}\norm{D^{*}A_{0}Du_{\varepsilon}}_{-1,p,\Sigma}^{*}\le\abs{\abr{A_{0}\xi,\xi}\eta}\hairspace\norm{\varphi}_{p,\Sigma}.\label{est: Appendix | lim sup ||A^0_mu u_e||}
\end{equation}
On⋅the⋅other hand,
\begin{equation}
\lim_{\varepsilon\to0}\bigrbr{c_{A}\norm{Du_{\varepsilon}}_{p,\Sigma}-C_{A}\norm{u_{\varepsilon}}_{p,\Sigma}}=c_{A}\abs{\eta}\hairspace\norm{\varphi}_{p,\Sigma}.\label{eq: Appendix | lim ||u_e||}
\end{equation}
The⋅desired result now follows from the assumption~(\ref{est: Appendix | A_mu bounded below}).

To⋅handle the general case, fix a Lebesgue point~$x_{0}\in\Sigma$
of⋅$A$ and let⋅$A_{0}=A\rbr{x_{0}}$.   In⋅the definition
of⋅$u_{\varepsilon}$, we also replace⋅$\varphi$ by⋅$\varphi_{\delta}\rbr x=\varphi\rbr{\delta^{-1}\rbr{x-x_{0}}}$
with⋅$\supp\varphi\subset B$ and⋅$\delta>0$ sufficiently small,
obtaining~$u_{\varepsilon,\delta}$. We write
\begin{equation}
D^{*}ADu_{\varepsilon,\delta}=D^{*}A_{0}Du_{\varepsilon,\delta}+D^{*}\rbr{A-A_{0}}Du_{\varepsilon,\delta}.\label{est: Appendix | A =00003D A^0 + Delta}
\end{equation}
The⋅first term on the right-hand side is dealt with by~(\ref{est: Appendix | lim sup ||A^0_mu u_e||}).
The⋅other one  is estimated as follows:
\[
\norm{D^{*}\rbr{A-A_{0}}Du_{\varepsilon,\delta}}_{-1,p,\Sigma}^{*}\le\abs{B_{\delta}\rbr{x_{0}}}^{1/p}\Delta_{\delta}\norm{Du_{\varepsilon,\delta}}_{\infty,\Sigma},
\]
where
\[
\Delta_{\delta}=\abs{B_{\delta}\rbr{x_{0}}}^{-1/p}\norm{A-A_{0}}_{p,B_{\delta}\rbr{x_{0}}}\to0
\]
as⋅$\delta\to0$, by the Lebesgue differentiation theorem. Therefore,
\begin{equation}
\limsup_{\varepsilon\to0}\norm{D^{*}\rbr{A-A_{0}}Du_{\varepsilon,\delta}}_{-1,p,\Sigma}^{*}\le\abs{B_{\delta}\rbr{x_{0}}}^{1/p}\Delta_{\delta}\abs{\eta}\hairspace\norm{\varphi}_{\infty,\Sigma}.\label{est: Appendix | ||Delta||}
\end{equation}
We plug  the estimates~(\ref{est: Appendix | ||Delta||}) and⋅(\ref{est: Appendix | lim sup ||A^0_mu u_e||})
into⋅(\ref{est: Appendix | A =00003D A^0 + Delta}) and then use (\ref{est: Appendix | A_mu bounded below})
together with⋅(\ref{eq: Appendix | lim ||u_e||}) to get
\[
c_{A}\abs{\eta}\hairspace\norm{\varphi_{\delta}}_{p,\Sigma}\le\abs{\abr{A_{0}\xi,\xi}\eta}\hairspace\norm{\varphi_{\delta}}_{p,\Sigma}+\abs{B_{\delta}\rbr{x_{0}}}^{1/p}\Delta_{\delta}\abs{\eta}\hairspace\norm{\varphi}_{\infty,\Sigma}.
\]
Observe here that both⋅$\abs{B_{\delta}\rbr{x_{0}}}^{1/p}$ and⋅$\norm{\varphi_{\delta}}_{p,\Sigma}$
scale as~$\delta^{d/p}$. Canceling the factor~$\delta^{d/p}$
and letting⋅$\delta\to0$ completes the proof.
\end{proof}
\begin{cor}
\label{cor:  Appendix | Diagonal coefficients invertible}Under the
hypotheses of Lemma~\textup{\ref{lem: Appendix | Weaker ellipticity condition}},
we have, for any $1\le k\le d$ and⋅\aev~$x\in\Sigma$,
\begin{equation}
\abs{A_{kk}\rbr x\hairspace\eta}\ge c_{A}\abs{\eta},\qquad\eta\in\C^{n}.\label{est: Appendix | Diagonal coefficients invertible}
\end{equation}
\end{cor}
Let us also  discuss some properties  of⋅$C^{0,\omega}$\nobreakdash-maps
and their interactions with⋅$\Lambda_{p}^{\omega}$\nobreakdash-spaces.

Obviously, the⋅$C^{0,\omega}$\nobreakdash-spaces form an algebra,
which is not closed under composition unless~$\omega\rbr t=t$. However,
they are invariant under composition with $C^{0,1}$\nobreakdash-functions.
Indeed, it follows from monotonicity of the function~$t\mapsto t^{-1}\omega\rbr t$
that $\omega\rbr{\lambda t}\le\lambda\omega\rbr t$ for any⋅$t\ge0$
and⋅$\lambda\ge1$, which proves the claim for pre-compositions with
Lipschitz functions. The⋅claim for post-compositions with Lipschitz
functions is obvious.

Functions of the class $C^{0,\omega}$ are easily seen to be
multipliers on $\Lambda_{p}^{\omega}$\nobreakdash-spaces. As⋅for⋅composition,
we provide a result  eventually showing that a change-of-vari\-ables
operator associated with a bi-Lipschitz⋅$C^{1,\omega}$\nobreakdash-map
is an isomorphism between $\Lambda_{q}^{1+\omega}$\nobreakdash-spaces.
In⋅what follows, $\partial\varkappa=iD\varkappa$ denotes the Jacobian
matrix of a map~$\varkappa$.

\begin{lem}
\label{lm: Ex | d(w^=00007B-1=00007D)}Let⋅$\varkappa\colon U\to V$
be a bi-Lipschitz map of an open~$U$ set onto~$V$. If⋅$\partial\varkappa\in C^{0,\omega}\rbr{\wb U}^{d\times d}$,
then⋅$\partial\rbr{\varkappa^{-1}}\in C^{0,\omega}\rbr{\wb V}^{d\times d}$.
\end{lem}
\begin{proof}
Notice that $\rbr{\partial\varkappa}^{-1}\in C^{0,\omega}\rbr{\wb U}^{d\times d}$
as a product of the functions~$\rbr{\det\partial\varkappa}^{-1}$
and⋅$\adj\partial\varkappa$ (the⋅adjugate of~$\partial\varkappa$),
both of which are~$C^{0,\omega}$. Therefore, so is $\partial\rbr{\varkappa^{-1}}=\rbr{\partial\varkappa}^{-1}\circ\varkappa^{-1}$.
\end{proof}

We set $T_{\kappa}u=u\circ\kappa$ for any map~$\kappa$ and any
function~$u$ for which the right-hand side makes sense. In⋅the⋅next
two lemmas, $\varkappa\colon\Sigma\to\Xi$ is a bi-Lipschitz bijection
between uniformly weakly Lipschitz domains~$\Sigma$ and~$\Xi$.
\begin{lem}
\label{lem: Appendix | Change of variables and Lambda_q^omega}Let⋅$q\in(1,\infty]$.
Then $u\in\Lambda_{q}^{\omega}\rbr{\Xi}$ if and only if $T_{\varkappa}u\in\Lambda_{q}^{\omega}\rbr{\Sigma}$,⋅with⋅$\abs{T_{\varkappa}u}_{\omega,q,\Sigma}\sim\abs u_{\omega,q,\Xi}$.
\end{lem}
\begin{proof}
For⋅$q=\infty$, we have
\[
\osc_{\infty}\rbr{T_{\varkappa}u;t}\le\omega\rbr{\lambda t}\hairspace\seminorm u_{\omega,\infty,\Xi}\le\lambda\omega\rbr t\hairspace\seminorm u_{\omega,\infty,\Xi},
\]
where⋅$\lambda$ is the larger of~$1$ and the Lipschitz constant
of⋅$\varkappa$, as discussed earlier.

Suppose now that~$q<\infty$. Using an operator  as in⋅(\ref{eq: Notation | E and osc_p}),
we extend a function~$u\in L_{q}\rbr{\Xi}$ to⋅$\wt u\in L_{q}\rbr{\R^{d}}$,
so that, for⋅$t\in\sbr{0,1}$,
\[
\osc_{q}\rbr{\wt u;t}\lesssim\osc_{q}\rbr{u;t}+t\norm u_{q,\Xi}.
\]
Let⋅$\wt u_{\delta}$ be a standard mollification of~$\wt u$. Since
$\Sigma$ need not be convex, we estimate the $L_{q}$\nobreakdash-modulus
of continuity  of the pull-back through its extension to⋅$\R^{d}$
by an operator~$\op E$ for⋅$\Sigma$ as in~(\ref{eq: Notation | E and osc_p}):
\[
\osc_{q}\rbr{T_{\varkappa}\wt u_{\delta};t}\le\osc_{q}\rbr{\op ET_{\varkappa}\wt u_{\delta};t}\le t\norm{D\op ET_{\varkappa}\wt u_{\delta}}_{q,\R^{d}}\lesssim t\norm{\wt u_{\delta}}_{1,q,\Xi},
\]
whence
\[
\osc_{q}\rbr{T_{\varkappa}\wt u;\delta}\lesssim\osc_{q}\rbr{T_{\varkappa}\wt u_{\delta};\delta}+\norm{T_{\varkappa}\rbr{\wt u_{\delta}-\wt u}}_{q,\Sigma}\lesssim\delta\norm{\wt u_{\delta}}_{1,q,\Xi}+\norm{\wt u_{\delta}-\wt u}_{q,\Xi}.
\]
Combining this with Lemmas~\ref{lem: C | Mollification Je and regularity}
and~\ref{lem: C | Bound for (Je - I)} yields
\[
\osc_{q}\rbr{T_{\varkappa}\wt u;\delta}\lesssim\omega\rbr{\delta}\hairspace\abs u_{\omega,q,\Xi}.
\]

The⋅reverse inequalities are proved similarly.
\end{proof}

\begin{lem}
\label{lem: Appendix | Change of variables and Lambda_q^(1+omega)}Let⋅$q\in(1,\infty]$
and⋅$\partial\varkappa\in C^{0,\omega}\rbr{\wb{\Sigma}}^{d\times d}$.
Then $u\in\Lambda_{q}^{1+\omega}\rbr{\Xi}$ if and only if $T_{\varkappa}u\in\Lambda_{q}^{1+\omega}\rbr{\Sigma}$,⋅with⋅$\abs{T_{\varkappa}u}_{1+\omega,q,\Sigma}\sim\abs u_{1+\omega,q,\Xi}$.
\end{lem}
\begin{proof}
It⋅suffices to estimate the $\smash[b]{L_{q}}$\nobreakdash-modulus
of continuity of⋅$DT_{\varkappa}u=T_{\varkappa}Du\cdot\nobreak\partial\varkappa$.
But⋅$Du\in\smash[b]{\Lambda_{q}^{\omega}}\rbr{\Xi}^{d}$ if and only
if⋅$T_{\varkappa}Du\in\smash[b]{\Lambda_{q}^{\omega}}\rbr{\Sigma}^{d}$,
with equivalent norms, by Lemma~\ref{lem: Appendix | Change of variables and Lambda_q^omega},
and both⋅$\partial\varkappa$ and⋅$\rbr{\partial\varkappa}^{-1}$
belong to⋅$\smash[t]{C^{0,\omega}\rbr{\wb{\Sigma}}^{d\times d}}$,
as observed in the proof of Lemma~\ref{lm: Ex | d(w^=00007B-1=00007D)},
so the conclusion follows.
\end{proof}

We are now in a position to prove our regularity result.
\begin{proof}[Proof of Theorem~\textup{\ref{thm: Appendix | Elliptic regularity}}]
We follow the standard difference quotient argument, introducing
appropriate regularizations of the  coefficients and data along the
way. As⋅usual, the proof naturally  divides into two parts: regularity
near the boundary or the interfaces and regularity away from them.
Let us concentrate on the interfaces, the other cases being similar
or simpler.

Our first step is to localize the operator~$\mathbb{A}_{\mu}$
to a coordinate neighborhood of the interface~$\Gamma$ and make
a change of variables to locally flatten the interface.

Choosing a fixed cutoff function on⋅$\ball$ and pulling it back
to each coordinate neighborhood of⋅$\Gamma_{ij}$, we obtain a finite
family of cutoff functions whose element~$\varphi$ corresponding
to a  chart~$\rbr{W,w}$ satisfies $\varphi\in C_{c}^{1,\smash[b]{\omega_{ij}}}\rbr W$
and⋅$\smash[t]{\wh{\varphi}}=T_{w^{-1}}\varphi\in C_{c}^{\infty}\rbr{\ball}$,
and whose sum~$\Phi$ is  at least~$1$ near~$\Gamma_{ij}$, so
that $1/\Phi$ is $C^{1,\omega_{ij}}$ there. Since $C^{1,\omega_{ij}}$\nobreakdash-func\-tions
are multipliers on⋅$\smash{\Lambda_{p}^{1+\omega_{ij}}}$\nobreakdash-spaces,
it suffices to estimate the $L_{p}$\nobreakdash-modulus of continuity
of⋅$D\varphi u$ on each side of~$\Gamma_{ij}$.

Fix a chart~$\rbr{W,w}$ for⋅$\Gamma_{ij}$, and let $\varphi$ be
the corresponding element of the family.  Then $z=\varphi u\in\smash{\mathring{W}_{p}^{1}}\rbr W^{n}$
solves
\begin{equation}
\mathbb{A}_{\mu}z=g+D^{*}G,\label{eq: Appendix | z}
\end{equation}
where⋅$g=\varphi f+\rbr{F-ADu}\cdot\rbr{D\varphi}^{*}$ and⋅$G=AuD\varphi+\varphi F$,
with⋅$D^{*}G$ being viewed as a functional on~$\smash{\set W_{p^{\prime}}^{1}}\rbr{\Omega;\C^{n}}$.

Writing $\mathbb{A}_{\mu}=\mathbb{A}-\mu\alpha$ with⋅$\alpha=1$,
define the operator~$\smash{\wh{\mathbb{A}}_{\mu}}\colon\smash{\mathring{W}_{p}^{1}}\rbr{\ball}^{n}\to W_{p}^{-1}\rbr{\ball}^{n}$⋅by
\[
\wh{\mathbb{A}}_{\mu}=T_{w}^{*}\mathbb{A}_{\mu}T_{w},
\]
where functions in the range of⋅$T_{w}$ are extended by zero outside~$W$.
Let $J_{w^{-1}}=\det\partial\rbr{w^{-1}}$, assuming, as we may,
that⋅$\det\partial w>0$. Then $DT_{w}\smash[t]{\wh v}=T_{w}D\smash[t]{\wh v}\cdot\partial w$
for every⋅$\smash[t]{\wh v}\in W_{p}^{1}\rbr{\ball}^{n}$ and $T_{w}^{*}=J_{w^{-1}}T_{w^{-1}}$.
Therefore,
\[
\wh{\mathbb{A}}_{\mu}=D^{*}\wh AD-\mu\wh{\alpha},
\]
with⋅$\smash[t]{\wh A}\zeta=T_{w}^{*}\rbr{A\zeta\partial w\cdot\rbr{\partial w}^{t}}$
and⋅$\smash[t]{\wh{\alpha}}=T_{w}^{*}\alpha$. Here and below, a
hat marks the image of an object under the change of variables induced
by⋅$w$, and a check marks the inverse operation, that is, the hat
operation associated with~$w^{-1}$. Namely, functions transform
by⋅$\smash[t]{\wh v}=\smash[b]{T_{w^{-1}}}v$, densities by⋅$\smash[t]{\wh h}=T_{w}^{*}h$,
fields under the divergence by⋅$\smash[t]{\wh H}=T_{w}^{*}H\rbr{\partial w}^{t}$
and coefficients as⋅$\smash[t]{\wh A}$ and⋅$\smash[t]{\wh{\alpha}}$
above. These rules preserve the pairings, so $v$ satisfies⋅$\smash[b]{\mathbb{A}_{\mu}}v=h+D^{*}H$
 in⋅$W_{p}^{-1}\rbr W^{n}$ if and only if $\smash{\wh{\mathbb{A}}_{\mu}\wh v}=\smash[t]{\wh h}+D^{*}\smash[t]{\wh H}$
in~$W_{p}^{-1}\rbr{\ball}^{n}$. It⋅follows that (\ref{eq: Appendix | z})
can be rewritten⋅as
\begin{equation}
\wh{\mathbb{A}}_{\mu}\wh z=\wh g+D^{*}\wh G.\label{eq: Appendix | z^}
\end{equation}
By⋅Lemma~\ref{lem: Appendix | Change of variables and Lambda_q^(1+omega)}
applied to⋅$w$ on⋅$\Omega_{i}\cap W$ and⋅$\Omega_{j}\cap W$, we
need only estimate the $L_{p}$\nobreakdash-modulus of continuity
of⋅$D\wh z$ on⋅$\ball_{+}$ and~$\ball_{-}$.

Our⋅next step is to regularize the equation~(\ref{eq: Appendix | z^}).
Since $\smash{\wh{\mathbb{A}}_{\mu}}$ need not be an isomorphism,
we cannot do this on⋅$\ball$ directly. Instead, we modify the coefficient~$A$
of⋅$\smash[b]{\mathbb{A}_{\mu}}$ so that the new coefficient~$A_{\delta}$
coincides with⋅$A$ outside⋅$W$ and the transformed coefficient~$\smash{\wh A_{\delta}}$
is a mollification of⋅$\smash[t]{\wh A}$ on each half-ball~$\ball_{\pm}$
near⋅$\supp\smash[t]{\wh{\varphi}}$. Mollifying $A$ itself would
not be enough, because $\smash[t]{\wh A}$ also contains the factors~$\partial w$
and⋅$J_{w^{-1}}$, which are not Lipschitz. The⋅field~$F$ on the
right is modified in the same way.

For⋅notational convenience, we use the subscript~``$\pm$'' for
objects associated with⋅$\ball_{\pm}$; \eg,  $\smash{\wh A_{\pm}}$
and⋅$\smash[t]{\wh z_{\pm}}$ are the restrictions of⋅$\smash[t]{\wh A}$
and⋅$\smash[t]{\wh z}$ to⋅$\ball_{\pm}$, and $\Omega_{\pm}$ denotes
whichever of⋅$\Omega_{i}$ and⋅$\Omega_{j}$ satisfies⋅$w\rbr{\Omega_{\pm}\cap W}=\ball_{\pm}$.
According to⋅(\ref{assu: Appendix| Coefficient A regular near Gamma_ij}),
$\smash{\wh A_{\pm}}\in C^{0,\omega}\rbr{\smash{\wb{\ball}_{\pm}}}$,
where⋅$\omega=\omega_{ij}$. Indeed, $\smash[t]{\wh A}\zeta=\smash[b]{J_{w^{-1}}}\rbr{A\zeta\partial w\cdot\rbr{\partial w}^{t}}\circ w^{-1}$,
and $C^{0,\omega}$\nobreakdash-functions form an algebra invariant
under compositions with Lipschitz maps. By⋅reflection across the
flat part of⋅$\partial\ball_{\pm}$ and then across the unit sphere,
we  extend $\smash{\wh A_{\pm}}$ to a⋅$C^{0,\omega}$-function
on⋅$\R^{d}$, keeping the old notation. Let⋅$\smash{\wh A_{\pm,\delta}}$,⋅$\delta>0$,
be a  mollification of⋅$\smash{\wh A_{\pm}}$ obtained by applying
an operator~$\op J^{\delta}$ as in~(\ref{def: C | Mollification Je}).
Fix  the characteristic function~$\chi$ of an open neighborhood
of⋅$\supp\varphi$ compactly contained in⋅$W$ and introduce $A_{\delta}=A+\chi\rbr{\smash{\wh A_{\delta}}-\smash[t]{\wh A}}\smash[t]{\spcheck}$,
where⋅$\smash{\wh A_{\delta}}=\smash{\wh A_{\#,\delta}}$. Replacing
the coefficient~$A$ by⋅$A_{\delta}$, we define the regularized
version of the operator~$\mathbb{A}_{\mu}$, denoted~$\mathbb{A}_{\delta,\mu}$.
Since $\norm{A_{\delta}-A}_{\infty,\Omega}\lesssim\omega\rbr{\delta}$
(see⋅Lemma~\ref{lem: C | Bound for (Je - I)}),  the operator~$\mathbb{A}_{\delta,\mu}$
is an isomorphism with uniformly bounded inverse for all⋅$\delta\in\set D_{0}=(0,\delta_{0}]$,
 provided that $\delta_{0}\le1$ is sufficiently small:
\begin{equation}
\norm{\mathbb{A}_{\delta,\mu}^{-1}h}_{1,p,\Omega}\lesssim\Norm h_{-1,p,\Omega},\qquad h\in\set W_{p}^{-1}\rbr{\Omega;\C^{n}}.\label{est: Appendix | Resolvent of A_=00007Bdelta,=00005Cmu=00007D bounded}
\end{equation}
 Similarly, we extend⋅$\smash[t]{\wh F_{\pm}}$ to a⋅$\Lambda_{p}^{\omega}$-function
on⋅$\R^{d}$, mollify it with⋅$\op J^{\delta}$ and set⋅$F_{\delta}=F+\chi\rbr{\smash[t]{\wh F_{\delta}}-\smash[t]{\wh F}}\smash[t]{\spcheck}$,
where⋅$\smash{\wh F_{\delta}}=\smash{\wh F_{\#,\delta}}$. Let⋅$u_{\delta}\in\set W_{p}^{1}\rbr{\Omega;\C^{n}}$
be the solution⋅of
\begin{equation}
\mathbb{A}_{\delta,\mu}u_{\delta}=f+D^{*}F_{\delta}.\label{eq: Appendix | u_delta}
\end{equation}
We have
\[
\mathbb{A}_{\delta,\mu}\rbr{u_{\delta}-u}=-\rbr{\mathbb{A}_{\delta,\mu}-\mathbb{A}_{\mu}}u+D^{*}\rbr{F_{\delta}-F},
\]
so⋅that, by⋅(\ref{est: Appendix | Resolvent of A_=00007Bdelta,=00005Cmu=00007D bounded})
and Lemma~\ref{lem: C | Bound for (Je - I)},
\begin{equation}
\begin{aligned}\norm{u_{\delta}-u}_{1,p,\Omega} & \lesssim\norm{A_{\delta}-A}_{\infty,\Omega}\norm u_{1,p,\Omega}+\norm{F_{\delta}-F}_{p,\Omega}\\
 & \lesssim\omega\rbr{\delta}\hairspace\bigrbr{\norm f_{p,\Omega}+\norm F_{p,\Omega}+\abs F_{\omega,p,\Omega_{+}\cap W}+\abs F_{\omega,p,\Omega_{-}\cap W}}.
\end{aligned}
\label{est: Appendix | ||u_delta - u|| in W_p^1}
\end{equation}

We now localize⋅$u_{\delta}$ in the same way as~$u$: repeating
the argument that led to⋅(\ref{eq: Appendix | z}) with⋅$A_{\delta}$,
$F_{\delta}$ and⋅$u_{\delta}$ in place of⋅$A$, $F$ and⋅$u$, we
see that $z_{\delta}=\varphi u_{\delta}\in\smash{\mathring{W}_{p}^{1}}\rbr W^{n}$
satisfies
\begin{equation}
\mathbb{A}_{\delta,\mu}z_{\delta}=g_{\delta}+D^{*}G_{\delta},\label{eq: Appendix | z_delta}
\end{equation}
where⋅$g_{\delta}$ and⋅$G_{\delta}$ are defined as⋅$g$ and⋅$G$
with the same replacements. Changing the variables in this equation
yields
\begin{equation}
\wh{\mathbb{A}}_{\delta,\mu}\wh z_{\delta}=\wh g_{\delta}+D^{*}\wh G_{\delta}\label{eq: Appendix | z_delta^}
\end{equation}
in~$W_{p}^{-1}\rbr{\ball}^{n}$. Note that, on⋅$\supp\smash[t]{\wh{\chi}}$,
which contains the supports of⋅$\smash[t]{\wh z_{\delta}}$, $\smash[t]{\wh g_{\delta}}$
and~$\smash[t]{\wh G_{\delta}}$, the transforms of⋅$A_{\delta}$
and⋅$F_{\delta}$ coincide with⋅$\smash{\wh A_{\delta}}$ and⋅$\smash{\wh F_{\delta}}$,
respectively. For⋅every⋅$r\in\rbr{0,1}$, the restriction of⋅$\smash{\wh{\mathbb{A}}_{\delta,\mu}}$
to functions supported in the closure of⋅$\ball_{r}=r\ball$ has a
left inverse bounded uniformly in⋅$\delta\in\set D_{0}$, that⋅is,
\begin{equation}
\norm{\wh v}_{1,p,\ball}\lesssim\norm{\wh{\mathbb{A}}_{\delta,\mu}\wh v}_{-1,p,\ball},\qquad\wh v\in\mathring{W}_{p}^{1}\rbr{\ball_{r}}^{n}.\label{est: Appendix | Left inverse of ^A_=00007Bdelta,mu=00007Dbounded}
\end{equation}
To⋅see this, choose a cutoff function~$\smash[t]{\wh{\varphi}^{\prime}}\in C_{c}^{0,1}\rbr{\ball}$
equal to~$1$ on⋅$\ball_{r}$, and set~$v=T_{w}\smash[t]{\wh v}$.
Then $\smash{\wh{\mathbb{A}}_{\delta,\mu}\wh v}=\smash{\wh{\varphi}^{\prime}\wh{\mathbb{A}}_{\delta,\mu}\wh v}$,
where
\[
\Norm{\varphi^{\prime}\mathbb{A}_{\delta,\mu}v}_{-1,p,\Omega}\lesssim\norm{\mathbb{A}_{\delta,\mu}v}_{-1,p,W}\sim\norm{\wh{\mathbb{A}}_{\delta,\mu}\wh v}_{-1,p,\ball},
\]
while, in view of⋅(\ref{est: Appendix | Resolvent of A_=00007Bdelta,=00005Cmu=00007D bounded}),
the norm of⋅$\mathbb{A}_{\delta,\mu}v$ is bounded below by the $W_{p}^{1}$\nobreakdash-norm
of~$\wh v$.

Using⋅(\ref{est: Appendix | ||u_delta - u|| in W_p^1}), we see that
\begin{equation}
\norm{z_{\delta}-z}_{1,p,\Omega}\lesssim\omega\rbr{\delta}\hairspace\bigrbr{\norm f_{p,\Omega}+\norm F_{p,\Omega}+\abs F_{\omega,p,\Omega_{+}\cap W}+\abs F_{\omega,p,\Omega_{-}\cap W}}.\label{est: Appendix | ||z_delta - z|| in W_p^1}
\end{equation}
Combined with the obvious estimate
\[
\osc_{p}\rbr{D\wh z_{\pm,\delta}-D\wh z_{\pm};\delta}\lesssim\norm{D\wh z_{\delta}-D\wh z}_{p,\ball_{\pm}},
\]
this shows that it suffices to prove that
\begin{equation}
\osc_{p}\rbr{D\wh z_{\pm,\delta};\delta}\lesssim\omega\rbr{\delta}\hairspace\bigrbr{\norm f_{p,\Omega}+\norm F_{p,\Omega}+\abs F_{\omega,p,\Omega_{+}\cap W}+\abs F_{\omega,p,\Omega_{-}\cap W}}\label{est: Appendix | Reduced to z_=00007Bpm,delta=00007D}
\end{equation}
for⋅all⋅$\delta\in\set D_{0}$, since the $L_{p}$\nobreakdash-modulus
of continuity of $D\wh z_{\pm}$ is trivially controlled at large
scales~$\delta>\delta_{0}$.

Let us look at the regularized problem~(\ref{eq: Appendix | z_delta^}).
Denote the flat portion of the boundary of⋅$\ball_{\pm}$ by~$S$.
We  introduce the  trace operator~$\gamma_{S}\colon\smash{W_{p^{\prime}}^{1}}\rbr{\ball}^{n}\to\smash{W_{p^{\prime}}^{1-1/p^{\prime}}}\!\rbr S^{n}$
 and the one-sided normal trace operators~$\gamma_{S_{\pm},n}\colon\smash{W_{p}^{1}\rbr{\ball_{\pm}}^{n\times d}}\to\smash{W_{p}^{1-1/p}\rbr S^{n}}$
 associated with the normal vector~$e_{d}$.  Then, by the Green
theorem,
\[
\rbr{\wh G_{\pm,\delta},D\wh u^{\prime}}_{\ball_{\pm}}=\rbr{D^{*}\wh G_{\pm,\delta},\wh u^{\prime}}_{\ball_{\pm}}\mp i\rbr{\gamma_{S_{\pm},n}\wh G_{\pm,\delta},\gamma_{S}\wh u^{\prime}}_{S},\qquad\wh u^{\prime}\in W_{p^{\prime}}^{1}\rbr{\ball}^{n},
\]
$D^{*}\smash{\wh G_{\pm,\delta}}$  now being a function.
Putting this into the right-hand side of⋅(\ref{eq: Appendix | z_delta^}),
we obtain
\begin{equation}
\wh{\mathbb{A}}_{\delta,\mu}\wh z_{\delta}=\wh g_{\delta}+D^{*}\wh G_{\#,\delta}-i\gamma_{S}^{*}\wh{\psi}_{\delta},\label{eq: Appendix | z_delta^as a transmisison problem}
\end{equation}
where⋅$\smash{\wh{\psi}_{\delta}}=\gamma_{S_{+},n}\smash{\wh G_{+,\delta}}-\gamma_{S_{-},n}\smash{\wh G_{-,\delta}}$.
This equation resembles a model transmission problem with right-hand
sides~$\smash{\wh g_{\pm,\delta}}+D^{*}\smash{\wh G_{\pm,\delta}}$
and conormal derivative jump~$i\smash{\wh{\psi}_{\delta}}$⋅(\cf~(\ref{eq: Appendix | Green's identity})).
 Such problems are treated in standard elliptic theory (see,⋅\eg,~\cite[Theorem~4.20]{McL:2000}),
yet ours is not directly covered. Nevertheless, the same proof carries
over, yielding
\begin{equation}
\begin{aligned}\norm{DD\wh z_{\pm,\delta}}_{p,\ball_{\pm}} & \lesssim\bigrbr{\norm{D\wh A_{\#,\delta}}_{\infty,\ball}+1}\norm{\wh z_{\delta}}_{1,p,\ball}\\
 & \quad+\norm{\wh g_{\delta}}_{p,\ball}+\norm{D^{*}\wh G_{\#,\delta}}_{p,\ball}+\norm{\wh{\psi}_{\delta}}_{1-1/p,p,S}
\end{aligned}
\label{est: Appendix | DD ^u_delta in L_p}
\end{equation}
uniformly in~$\delta\in\set D_{0}$.  Two comments are in order.
The⋅role of the coercivity condition, which drives tangential regularity,
is now played by the left-inverse bound~(\ref{est: Appendix | Left inverse of ^A_=00007Bdelta,mu=00007Dbounded})
for~$\smash{\wh{\mathbb{A}}_{\delta,\mu}}$. This bound also implies,
via Corollary~\ref{cor:  Appendix | Diagonal coefficients invertible},
that the diagonal  coefficient~$\rbr{\smash{\wh A_{\delta}}}_{dd}$
has a uniformly bounded inverse on⋅$\ball_{r}$, which is what the
proof of normal regularity requires.

Once we have (\ref{est: Appendix | DD ^u_delta in L_p}), we can estimate
the $L_{p}$\nobreakdash-modulus of continuity of⋅$D\smash{\wh z_{\pm,\delta}}$
using the relation
\[
\sup_{t>0}t^{-1}\osc_{p}\rbr{D\wh z_{\pm,\delta};t}\sim\norm{DD\wh z_{\pm,\delta}}_{p,\ball_{\pm}}.
\]
Applying Lemma~\ref{lem: C | Mollification Je and regularity}
to the right-hand side of the resulting bound gives
\[
\begin{aligned}\osc_{p}\rbr{D\wh z_{\pm,\delta};\delta} & \lesssim\omega\rbr{\delta}\hairspace\bigl(\norm{z_{\delta}}_{1,p,W}+\norm{u_{\delta}}_{1,p,W}\\
 & \hphantom{{}\lesssim\omega\rbr{\delta}\hairspace\bigl(}{}+\norm f_{p,W}+\norm F_{p,W}+\abs F_{\omega,p,\Omega_{+}\cap W}+\abs F_{\omega,p,\Omega_{-}\cap W}\bigr)
\end{aligned}
\]
(we have used that $\delta\le\omega\rbr{\delta}$ for any⋅$\delta\le1$,
by regularity of~$\omega$).  Now⋅(\ref{est: Appendix | Reduced to z_=00007Bpm,delta=00007D})
follows from⋅(\ref{est: Appendix | Resolvent of A_=00007Bdelta,=00005Cmu=00007D bounded})
and~(\ref{eq: Appendix | u_delta}).
\end{proof}
\begin{proof}[Proof of Lemma~\textup{\ref{lem: Appendix | Regular Neumann data fit}}]
In⋅a⋅boundary chart~$\rbr{W,w}$, with⋅$\Omega_{W}=\Omega\cap W$
and⋅$\varphi$ the corresponding element of a partition of unity,
we have, for every⋅$u^{\prime}\in W_{p^{\prime}}^{1}\rbr{\Omega_{W}}^{n}$,
\[
\rbr{\varphi\varrho,\gamma u^{\prime}}_{\partial\Omega\cap W}=\rbr{J_{\smash[t]{w^{-1}}}^{S}\wh{\varphi}\wh{\varrho},\gamma_{S}\wh u^{\prime}}_{S},
\]
where⋅$\smash{J_{\smash[t]{w^{-1}}}^{S}}\in C^{0,\omega}\rbr S$ is
the surface Jacobian. Let⋅$\smash{\wh v}\in W_{p}^{1}\rbr{\ball_{+}}^{n}$
be a  lifting of⋅$\smash{\wh{\varphi}\wh{\varrho}}$ that vanishes
near⋅$\partial\ball\cap\partial\ball_{+}$, and set⋅$\smash{\wh H}=i\smash{J_{\smash[t]{w^{-1}}}^{S}}\smash{\wh v}\otimes e_{d}$.
Then $\smash{\wh H}\in\smash{\Lambda_{p}^{\omega}\rbr{\ball_{+}}^{n\times d}}$,
and its divergence involves only⋅$D_{d}$, so  $D^{*}\smash{\wh H}=-\smash[t]{\wh h}$
with⋅$\smash[t]{\wh h}=-i\smash{J_{\smash[t]{w^{-1}}}^{S}}D_{d}\smash{\wh v}\in L_{p}\rbr{\ball_{+}}^{n}$.
By⋅Green's theorem,
\[
\rbr{J_{\smash[t]{w^{-1}}}^{S}\wh{\varphi}\wh{\varrho},\gamma_{S}\wh u^{\prime}}_{S}=\rbr{\wh H,D\wh u^{\prime}}_{\ball_{+}}+\rbr{\wh h,\wh u^{\prime}}_{\ball_{+}}.
\]
Pulling this back to⋅$\Omega_{W}$, we get
\[
\rbr{\varphi\varrho,\gamma u^{\prime}}_{\partial\Omega\cap W}=\rbr{H,Du^{\prime}}_{\Omega_{W}}+\rbr{h,u^{\prime}}_{\Omega_{W}},
\]
where⋅$h=\smash{T_{w^{-1}}^{*}}\smash[t]{\wh h}$ and⋅$H=J_{w}T_{w}\smash[t]{\wh H}\cdot\rbr{\partial w}^{-t}$
belong to⋅$L_{p}\rbr{\Omega_{W}}^{n}$ and⋅$\Lambda_{p}^{\omega}\rbr{\Omega_{W}}^{n\times d}$
by Lemmas~\ref{lm: Ex | d(w^=00007B-1=00007D)} and~\ref{lem: Appendix | Change of variables and Lambda_q^omega}.

Summing over the charts gives the desired representation on⋅$\Omega$,
and the estimate follows from the boundedness of the lifting.
\end{proof}

\bibliographystyle{amsalpha}
\bibliography{bibliography}

\end{document}